\documentclass[11pt,a4paper]{article}
\usepackage[round,authoryear,comma]{natbib}
\usepackage{optimstyle}
\usepackage{tabularx,pifont,calc,xurl,fontawesome5,nicefrac}
\usepackage{placeins}
\usepackage{adjustbox}
\usepackage{tikz}
\usetikzlibrary{arrows.meta,positioning,shapes.geometric,calc}
\colorlet{cheapline}{optimBlue}
\colorlet{cheapfill}{optimBlue!5!white}
\newcommand{\R}{\mathbb R}
\newcommand{\E}{\mathbb E}
\newcommand{\ip}[2]{\langle #1,#2\rangle}
\newcommand{\norm}[1]{\lVert#1\rVert}

\DeclareMathOperator{\argmin}{arg\,min}
\DeclareMathOperator{\dist}{dist}
\newcommand{\by}[2]{\mathrel{\underset{#2}{#1}}}

\newcolumntype{P}[1]{>{\raggedright\arraybackslash}p{#1}}
\newcolumntype{Y}{>{\raggedright\arraybackslash}X}

\setlist[itemize]{label=$\vartriangleright$}
\BeforeBeginEnvironment{theorem}{\Needspace{6\baselineskip}}
\BeforeBeginEnvironment{lemma}{\Needspace{6\baselineskip}}
\BeforeBeginEnvironment{proposition}{\Needspace{6\baselineskip}}
\BeforeBeginEnvironment{definition}{\Needspace{6\baselineskip}}
\BeforeBeginEnvironment{assumption}{\Needspace{7\baselineskip}}
\tcbset{optimassumption/.style={optimbox,colback=optimGreen!5!white,
  colframe=optimGreen!45!white,borderline west={1.2pt}{0pt}{optimGreen}}}
\pretocmd{\subsection}{\Needspace{8\baselineskip}}{}{}
\renewcommand{\optimtitlefont}{\fontsize{16}{21}\selectfont\bfseries\boldmath}
\title{Descending Everest, Part II\texorpdfstring{\\[3pt]}{: }
Oracle Complexity of Minimizing a Sum\texorpdfstring{\\}{ }
of Convex \texorpdfstring{$(L_0,L_1)$}{(L0,L1)}-Smooth Functions\texorpdfstring{\thanks{The Russian version of this paper is available at \url{https://github.com/chitoraphelion/descending-everest-with-first-order-oracles}.}}{}}
\author{%
  Mikhail Rudomanov\thanks{MIRIAI, Moscow Independent Research Institute of Artificial Intelligence, Moscow, Russia. \email{rudomanov.ma@phystech.edu}.}\,\thanks{Equal contribution.}\quad
  Ivan Fomin\thanks{Yandex, Moscow, Russia. \email{ash-ivanfomin@yandex-team.ru}.}\,\footnotemark[3]\quad
  Alexander V. Gasnikov\thanks{Innopolis University, Innopolis, Russia; \email{gasnikov@yandex.ru}.}\quad
  Fedor S. Stonyakin\thanks{V.\,I.~Vernadsky Crimean Federal University, Simferopol, Russia. \email{fedyor@mail.ru}.}}
\hypersetup{pdfauthor={Mikhail Rudomanov, Ivan Fomin, Alexander V. Gasnikov, Fedor S. Stonyakin}}
\date{}
\begin{document}
\begingroup
\makeatletter
\renewcommand{\@fnsymbol}[1]{\ensuremath{\ifcase#1\or *\or\spadesuit\or\text{\tiny\faMountain}\or\clubsuit\or\dagger\or\ddagger\fi}}
\let\origmakefntext\@makefntext
\long\def\@makefntext#1{\rightskip\z@\@plus3em\relax\origmakefntext{#1}}
\makeatother
\maketitle
\endgroup
\begin{abstract}
We study the minimization of the average of $n$ convex $(L_0,L_1)$-smooth
functions: the Hessian norm of each function is bounded by an affine function
of its gradient norm. We assume known upper bounds $R$ on the distance
from the starting point to a minimizer and $\delta$ on the difference between the minimum of the average and the average
of the infima of the components. We derive upper and lower bounds on the oracle complexity
of reaching accuracy $\varepsilon$ in function value in expectation. Up to logarithmic factors, both bounds are of order
\[
 n\bigl(1+(L_1R)^{2/3}\bigr)
 +\sqrt{\frac{n(L_0+L_1^2\delta)R^2}{\varepsilon}+n(L_1R)^2}.
\]
The upper bound is attained by the proposed method YETI. The method combines
ideas of accelerated proximal methods and variance reduction methods
with a heuristic for choosing additional points.
The number of oracle calls made by the heuristic is bounded, which
preserves the upper bound on the total number of queries needed to reach the target accuracy.
In the experiments, YETI is compared with first-order methods
on synthetic problems and Poisson regression problems.
We also examine the stopping criterion and the effect of the heuristic separately.
The lower bound is obtained by an explicit construction of
families of instances that are hard for any first-order method.
\end{abstract}

\begin{tcolorbox}[
  optimbox,
  colback=optimBlue!3!white,
  colframe=optimBlue!40!white,
  borderline west={1.4pt}{0pt}{optimBlue},
  top=7pt,bottom=7pt
]
\small
\noindent\textbf{Use of language models.}\enspace
The models \textit{Codex Astra 6} and \textit{Claude Opus 5.5} were used
in obtaining the theoretical results of this work.
They were also used to prepare and typeset the manuscript and to translate
the original text into English.
Nevertheless, every line of the paper, including every statement, every
proof, and every computation, has been personally edited and checked
by the authors. The authors bear sole responsibility for the correctness
of all mathematical statements, proofs, and conclusions.
\end{tcolorbox}

\section{Introduction}\label{sec:intro}

We consider the problem of minimizing a finite sum
\[
 f(x)=\frac1n\sum_{i=1}^n f_i(x)\to\min_{x\in\R^d}.
\]
Problems of this form arise in empirical risk minimization
in machine learning, where $f_i$ is the loss on the $i$-th example
of the training set \citep{shalev2014understanding,bottou2018optimization}.
In federated and distributed learning, the component $f_i$ is the loss
on the data of the $i$-th device or computing node
\citep{mcmahan2017communication}. The gradient of a single component is often
much cheaper to compute than the full gradient. Therefore, we measure the complexity
of methods by the number of calls to an oracle that returns the value and the gradient of a single component.

For a sum of $L$-smooth functions, variance reduction methods
\citep{johnson2013svrg,xiao2014proxsvrg} make it possible to reduce the number of oracle calls.
They periodically compute the full gradient of $f$, and at the remaining
steps they use the gradients of individual components $f_i$.
Lower bounds for arbitrary randomized methods show
which dependencies on the number of components and on the accuracy cannot be improved
\citep{woodworth2016finite}.

However, the $L$-smoothness condition fails, for example, for exponential
loss functions. The $(L_0,L_1)$-smoothness condition was introduced by
\citet{zhang2019gradient} to explain the success of methods
with gradient clipping. When training an LSTM-based language model, the authors
observed that the local smoothness constant grows together with the gradient norm.
A similar dependence was also found when training transformers
\citep{crawshaw2022robustness}.
The $(L_0,L_1)$-smoothness condition admits precisely this kind of curvature growth.
Its generalizations and first-order methods for such functions were studied
in
\citet{li2023generalized,gorbunov2024methods,vankov2025optimizing,tyurin2025unified}.
For finite sums, one has to account simultaneously for the number of components and for the growth of curvature.
For example, in Poisson regression with the log link
\citep{stanPoisson}
$f_i(x)=\exp(a_i^{\mathsf T}x)-b_i a_i^{\mathsf T}x$, where $b_i\ge0$,
the Hessian is unbounded on the whole space, but
\[
 \norm{\nabla^2 f_i(x)}
 \le b_i\norm{a_i}^2+\norm{a_i}\norm{\nabla f_i(x)}.
\]
In phase retrieval from intensity measurements, one uses the losses
$f_i(x)=(\ip{a_i}x^2-b_i)^2$
\citep{candes2015wirtinger}. Here, too, one has to minimize a finite sum
of components with unbounded curvature.
Another example is the exponential loss in classification
\citep[Example~1.5]{gorbunov2024methods}.
In all these problems, when the number of observations is large, computing the full gradient
may be expensive, and the gradient itself is not globally Lipschitz.
Therefore, one needs methods that account for the growth of curvature and do not compute
the gradients of all components at every step. Below we assume
that $f$ attains its minimum, and that the distance to the minimizer, the initial optimality gap,
and the heterogeneity of the components are bounded.

\paragraph{The peculiarity of finite sums.}
In the problem under consideration, each component satisfies the inequality
\[
 0\preceq\nabla^2f_i(x)
 \preceq\bigl(L_{0,i}+L_1\norm{\nabla f_i(x)}\bigr)I.
\]
Even at a minimizer of the average, the gradients of individual components can be
large: the equality $\nabla f(x^\ast)=0$ only means that their sum
is zero. Therefore, a small norm of the full gradient does not imply small
smoothness constants of the components, on which the bounds of variance reduction
methods depend. What is the oracle complexity of minimizing a finite sum in this situation?
We answer this question in this work.

\paragraph{Notation.} We denote by $\varepsilon>0$ the required accuracy of the solution in function value: the method must output a point $\hat x$ such that $f(\hat x)-f^\ast\le\varepsilon$. Throughout, $\norm{\cdot}$ denotes the Euclidean norm of a vector and the spectral norm of a matrix, $[x,y]$ is the segment with endpoints $x,y\in\R^d$, and $B(x,r)$ is the closed ball of radius $r$ centered at $x$. We denote the sum of the first $k$ terms of a sequence by the corresponding capital letter: for example, $A_k=a_1+\dots+a_k$ and $A_0=0$.

For randomized methods, the stated accuracy is required
in expectation. All projections are Euclidean; $\Pi_C$
denotes the projection onto a nonempty closed convex set $C$.

\paragraph{Main result.}
Let the components be convex, twice continuously differentiable, and bounded
from below, and let $n^{-1}\sum_iL_{0,i}\le L_0$. We assume that the starting point $x_0$
is zero: this can be achieved by a shift, with respect to which the problem class
is invariant. We assume that $f$ attains its minimum at a point
$x^\ast$ with $\norm{x^\ast}\le R$ and $f(0)-f^\ast\le\Delta$.
Under the conditions of Theorem~\ref{thm:main}, the upper and lower complexity bounds
are, up to logarithmic factors, of order
\[
 n\bigl(1+(L_1R)^{2/3}\bigr)
 +\sqrt{\frac{n(L_0+L_1^2\delta)R^2}{\varepsilon}+n(L_1R)^2}.
\]
Here, complexity means the worst-case expected number of queries
of values and gradients of individual components that suffices to guarantee
$\E[f(\hat x)-f^\ast]\le\varepsilon$.
The lower bound holds for arbitrary randomized methods;
the dimension of the problem is not bounded.
The precise definition of the class and the logarithmic factors are given
in Section~\ref{sec:setting}.

\paragraph{Note on the logarithmic factors.}
In the upper bound, the hidden factor is $\log^2N$, where $N$
from~\eqref{eq:complexity-parameters} grows linearly in $\Delta/\varepsilon$.
In the worst case, when $\Delta$ is exponentially large in $L_1R$,
this logarithm can be of order $L_1R$, so the factor $\log^2N$ can be of order $(L_1R)^2$,
that is, it can give an additional polynomial factor.

The term $\sqrt{nL_0R^2/\varepsilon}$ is known from the theory of minimizing
a sum of $L$-smooth functions \citep{woodworth2016finite}.
The term $n(L_1R)^{2/3}$ reflects the number of oracle calls
required to approach the minimum from the region where the gradients
are exponentially large. The term
\[
 \sqrt{\frac{nL_1^2R^2\delta}{\varepsilon}+n(L_1R)^2}
\]
arises because the curvature of the components near the minimum
is determined by their gradients, which can remain large even at the
minimizer of the average.

\paragraph{Method and experiments.}
The upper bound is attained by the method YETI. When the local curvature is small,
it performs an accelerated step, solving the auxiliary subproblem by SVRG,
and when it is large, it projects the auxiliary point onto a halfspace
containing the minimizer. In addition, YETI uses a heuristic
with a bounded number of oracle queries: it makes it possible to take advantage of
good proposed points without worsening the theoretical bound.
In the experiments, YETI is compared with other methods, and
its stopping criterion is tested.

\paragraph{Relation to known methods.}
YETI uses the accelerated proximal scheme
of Monteiro--Svaiter \citep{ms2013} and the idea of adaptive choice of
regularization without solving an implicit equation at every step
\citep{carmon2022adaptive}. The local subproblem is solved by
a projected variant of SVRG. Acceleration via approximate solution of
regularized subproblems was also studied in the Catalyst scheme
\citep{lin2015catalyst,lin2018catalyst} and in
\citet{frostig2015}. In YETI, the choice of the step depends on the computed curvature
estimate of the sum: if this estimate is too large, the method constructs
a separating halfspace instead of a step. The complexity bound accounts for the queries
of both the local solver and the heuristic.

The lower bound relies on Yao's minimax principle \citep{yao1977}: a hard instance
is chosen at random from a finite family, after which we average
over the random variables of the algorithm. Hidden coordinates and the distribution
of information across different components continue the known approach to constructing
lower bounds for finite-sum minimization \citep{woodworth2016finite}.
For the class under consideration, we apply exponential
transformations of convex chains and choose their parameters so
that the constraints on $L_1$, $\delta$, and the initial optimality gap are satisfied.

The present paper continues our previous work \citet{rudomanov2026everest}, where
the oracle complexity of this class was found for a single function ($n=1$):
its order is $\sqrt{L_0R^2/\varepsilon}+L_1R$, and for a small
initial optimality gap it is $\sqrt{L_0R^2/\varepsilon}+(L_1R)^{2/3}$.
Here we consider a finite sum, and the oracle returns the value
and the gradient of only one component; therefore, to obtain the upper bound,
we use variance reduction.
The results for a single function
\citep{vankov2025optimizing,tyurin2026accelerated} concern
an oracle that computes the full gradient. In particular, the optimality
of the method of \citet{tyurin2026accelerated} for sufficiently small $\varepsilon$ is based
on the classical lower bound for the smooth subclass and concerns only
the leading term $\sqrt{L_0R^2/\varepsilon}$. The optimality of the term
$n(L_1R)^{2/3}$ for an oracle returning the gradient of a single component
does not follow from it.
Similarly,
the bounds for stochastic methods in the case of a common minimizer of all components
\citep{gorbunov2024methods}
do not cover an arbitrary positive heterogeneity $\delta$.

\paragraph{Organization of the paper.}
The problem setting and the main result are given in Section~\ref{sec:setting},
the method YETI and its analysis in Sections~\ref{sec:cheap}--\ref{sec:cheap-analysis},
the main experiments in Section~\ref{sec:experiments}
(further experiments and their details are in Appendix~\ref{app:experiments}), and the lower bound
in Section~\ref{sec:lower}. The main text contains
the results and a proof sketch of the lower bound, while auxiliary
proofs and variants of the method are deferred to the appendices. A summary table
of known bounds is given in Appendix~\ref{app:literature}.

\section{Problem setting and main bounds}\label{sec:setting}

We consider the function $f=n^{-1}\sum_{i=1}^n f_i$ on $\R^d$.
We are given $n\ge2$, $L_0,\delta\ge0$, and $L_1,R,\Delta,\varepsilon>0$.
We assume that the starting point is zero; this can be achieved by a shift of coordinates.
All general assumptions are collected below.

\begin{assumption}[components]\label{asm:components}
The functions $f_i:\R^d\to\R$, $1\le i\le n$, are convex,
twice continuously differentiable, and bounded from below.
\end{assumption}

\begin{assumption}[generalized smoothness]\label{asm:smoothness}
Numbers $L_{0,i}\ge0$ are known such that
\begin{equation}\label{eq:component-class}
 \begin{gathered}
 \norm{\nabla^2 f_i(x)}\le L_{0,i}+L_1\norm{\nabla f_i(x)}
 \quad(x\in\R^d,\ 1\le i\le n),\\
 \frac1n\sum_{i=1}^n L_{0,i}\le L_0.
 \end{gathered}
\end{equation}
\end{assumption}

\begin{assumption}[minimizer]\label{asm:minimizer}
There exists a minimizer $x^\ast\in\argmin f$ such that $\norm{x^\ast}\le R$.
We denote $f^\ast=f(x^\ast)$.
\end{assumption}

\begin{assumption}[initial optimality gap]\label{asm:initial-gap}
We have $f(0)-f^\ast\le\Delta$, and $\Delta\ge2\varepsilon$.
\end{assumption}

\begin{assumption}[heterogeneity]\label{asm:heterogeneity}
\begin{equation}\label{eq:problem-bounds}
 f^\ast-\frac1n\sum_{i=1}^n\inf_{x\in\R^d}f_i(x)\le\delta.
\end{equation}
\end{assumption}

\begin{assumption}[range of accuracies]\label{asm:accuracy}
\begin{equation}\label{eq:nondegenerate}
 \frac{L_1^2R^2\delta}{\varepsilon}\ge3.
\end{equation}
\end{assumption}

\begin{definition}[problem class]\label{def:class}
We denote by $\mathcal P_d$ the set of collections of functions
$f_1,\ldots,f_n:\R^d\to\R$ and known numbers $L_{0,1},\ldots,L_{0,n}$
satisfying Assumptions~\ref{asm:components}--\ref{asm:heterogeneity}.
\end{definition}

Assumption~\ref{asm:smoothness} bounds only the average
of the numbers $L_{0,i}$; individual $L_{0,i}$ may exceed $L_0$.
It suffices to construct the lower bound on the subclass $L_{0,i}=L_0$.
For $C^2$ functions, the Hessian bound from Assumption~\ref{asm:smoothness} is equivalent to
symmetric $(L_{0,i},L_1)$-smoothness
\citep[Assumption~1.3, Eq.~(5)]{gorbunov2024methods}:
\[
 \norm{\nabla f_i(x)-\nabla f_i(y)}
 \le\Bigl(L_{0,i}+L_1\sup_{z\in[x,y]}\norm{\nabla f_i(z)}\Bigr)\norm{x-y}.
\]
Assumption~\ref{asm:heterogeneity} requires knowing neither $f^\ast$
nor the infima of the components.
Assumption~\ref{asm:accuracy} guarantees that in the hard instance
with a hidden pair of functions from Section~\ref{sec:lower}, the initial optimality gap exceeds
the required accuracy (Lemma~\ref{lem:lower-pair}).

\paragraph{Oracle model.}
The algorithm knows the problem parameters and the numbers $L_{0,i}$.
It obtains information about the functions only from an oracle that returns
$(f_i(x),\nabla f_i(x))$ for a query $(i,x)$.
Queries are allowed anywhere in $\R^d$.
The query points, the stopping time, and the output are measurable
functions of the received oracle answers and the random variables of the algorithm.
Complexity is measured by the number of oracle calls
\citep[Section~1.1]{nesterov2018lectures}.
Computing the function value and the full gradient at one point requires $n$ queries.

\begin{definition}[oracle complexity]\label{def:oracle}
Let $\tau_{\mathsf M}(p)$ be the number of oracle calls made by the algorithm $\mathsf M$
on the input $p$, and let $\hat x_{\mathsf M}(p)$ be its output.
The algorithm is called valid if for any $d$ and $p\in\mathcal P_d$
it stops almost surely and
\[
 \E\bigl[f(\hat x_{\mathsf M}(p))-f^\ast\bigr]\le\varepsilon.
\]
The complexity of an algorithm on the class and the minimax complexity are defined as
\begin{equation}\label{eq:minimax}
 \begin{aligned}
 \mathcal T(\mathsf M)
 &=\sup_{d\ge1}\ \sup_{p\in\mathcal P_d}\E\tau_{\mathsf M}(p),\\
 \mathcal T^\ast&=\inf_{\mathsf M\ \text{valid}}\mathcal T(\mathsf M).
 \end{aligned}
\end{equation}
\end{definition}

To bound the number of queries, we set
\begin{equation}\label{eq:complexity-parameters}
 N_0=4+n+\frac{(L_0+L_1^2\delta)R^2}{\varepsilon}+(L_1R)^2,
 \qquad N=N_0-2+\frac\Delta\varepsilon.
\end{equation}

For nonnegative quantities, the notation $a\lesssim b$ means
$a\le Cb$ with an absolute constant $C$
independent of the problem parameters, the dimension, and the algorithm.
We write $a\sim b$ if $a\lesssim b\lesssim a$.
The notation $O$ and $\Omega$ is used in the same uniform sense,
and $\widetilde\Theta$ also hides powers of $\log N$.
The lower bound uses the quantity $N_0$, which does not depend on $\Delta$,
and the upper bound uses the quantity $N$, which by Assumption~\ref{asm:initial-gap}
is not smaller than $N_0$.
We set
\begin{equation}\label{eq:complexity-order}
 T_0=
 n\bigl(1+(L_1R)^{2/3}\bigr)
 +\sqrt{\frac{n(L_0+L_1^2\delta)R^2}{\varepsilon}+n(L_1R)^2}.
\end{equation}

\begin{theorem}[two-sided complexity bound]\label{thm:main}
Under Assumptions~\ref{asm:components}--\ref{asm:accuracy}, we have
\begin{equation}\label{eq:main-bounds}
 \frac{T_0}{\log^2N_0}\lesssim\mathcal T^\ast
 \lesssim T_0\log^2N.
\end{equation}
\end{theorem}

The upper bound is attained by the method YETI, including
the case of different $L_{0,i}$.
The lower bound holds already for $L_{0,i}=L_0$.

The upper bound is proved in Theorem~\ref{thm:cheap-upper},
and the lower bound in Theorem~\ref{thm:lower}. The proofs are given
in Appendices~\ref{app:upper-constants} and~\ref{app:lower-constants}.
The question of closing the gap in the logarithmic factors remains open.

\section{The method YETI}\label{sec:cheap}
\begin{figure}[H]
\centering\begin{tikzpicture}[
 x=1mm,y=1mm,font=\small,
 box/.style={draw=optimBlue!65,fill=optimBlue!5!white,rounded corners=2pt,
 line width=.55pt,align=center,inner sep=5pt},
 arr/.style={-{Stealth[length=2mm]},line width=.65pt,draw=optimBlue},
 lab/.style={font=\footnotesize,fill=white,inner sep=2pt,align=center}]
\node[box,text width=108mm] (fast) at (0,0) {
 \textbf{1. Heuristic block}\\[2pt]
 Proposes points based on the history of previous queries.\\
 Only points that do not increase the function value are accepted.\\
 The number of oracle calls in the block is bounded.};
\node[box,text width=85mm] (cert) at (0,-27) {
 \textbf{2. Accuracy check}\\[2pt]
 Is the gap between the best function value\\
 and the lower bound on its minimum small enough?};
\node[box,text width=20mm] (answer) at (66,-27) {Return\\the point};
\node[box,text width=94mm] (gate) at (0,-52) {
 \textbf{3. Main iteration}\\[2pt]
 Compute component gradients and the curvature estimate.\\
 Is the curvature small enough for the local subproblem?};
\node[box,text width=61mm,anchor=north] (cut) at (-35,-74) {
 \textbf{Build a cut}\\[3pt]
 Large curvature yields a lower\\
 bound on the gap. Convexity defines\\
 a halfspace containing the minimizer.\\
 Projection moves the auxiliary\\
 point closer to the minimizer.};
\node[box,text width=61mm,anchor=north] (local) at (35,-74) {
 \textbf{Take an accelerated step}\\[3pt]
 On a small ball, the smoothness\\
 constants of components are known.\\
 SVRG solves the local subproblem\\
 to a verifiable accuracy. The solution\\
 determines the update of the iterates.};
\draw[arr] (fast) -- (cert);
\draw[arr] (cert.east) -- node[lab,above]{yes} (answer.west);
\draw[arr] (cert) -- node[lab,right]{no} (gate);
\draw[draw=optimBlue,line width=.65pt] (gate.south) -- (0,-66);
\draw[arr] (0,-66) -- (-35,-66) -- node[lab,left,pos=.5]{no} (cut.north);
\draw[arr] (0,-66) -- (35,-66) -- node[lab,right,pos=.5]{yes} (local.north);
\coordinate (join) at (0,-114);
\draw[draw=optimBlue] (cut.south) |- (join);
\draw[draw=optimBlue] (local.south) |- (join);
\draw[arr] (join) -- (0,-120) -- (-76,-120) -- (-76,0) -- (fast.west);
\end{tikzpicture}
\caption{Flowchart of the method YETI.}\label{fig:cheap-flow}
\end{figure}
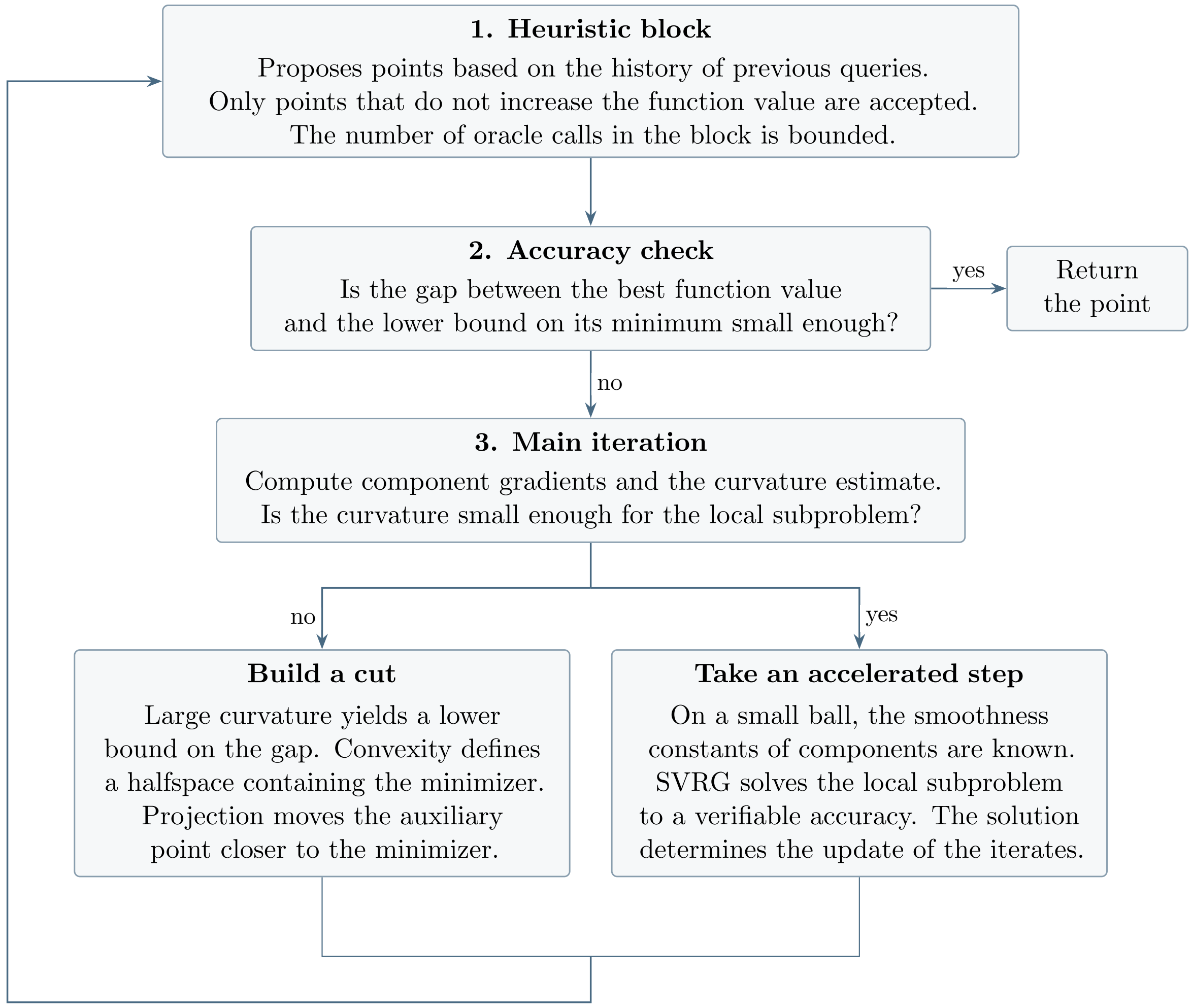

YETI combines an accelerated proximal method with a heuristic.
The heuristic proposes a point $z\in B(0,R)$ based on the query history.
The method accepts it if $f(z)\le f(x)$. A single heuristic block
is allowed at most $n\lceil\log N\rceil$ oracle queries, including
the verification of the proposed points. If the block has not terminated the method,
the main iteration is performed (Figure~\ref{fig:cheap-flow}).

\paragraph{Accelerated step and cut.}
From the gradients of the components at the point $y$, we compute the curvature estimate
\[
 \ell(y)=\frac1n\sum_i\bigl(L_{0,i}+L_1\norm{\nabla f_i(y)}\bigr).
\]
By Lemma~\ref{lem:observed-curvature}, the following inequalities hold:
\[
 \norm{\nabla^2f(y)}\le\ell(y)
 \le\tfrac32L_0+2L_1^2\bigl(\delta+f(y)-f^\ast\bigr),
 \qquad \ell(z)\le e^{L_1\norm{z-y}}\ell(y).
\]
If the quantity $\ell(y)$ is small, then in the ball of radius $r=\min\{R,1/(2L_1)\}$
one can solve the regularized subproblem by SVRG and perform
an accelerated step. If, however, the curvature check fails,
then the upper bound on $\ell(y)$ implies a lower bound on the optimality gap $f(y)-f^\ast$.
By convexity (Assumption~\ref{asm:components}), $\ip{\nabla f(y)}{x^\ast-y}\le-(f(y)-f^\ast)$,
so this bound defines a halfspace $\mathcal H_y$
containing the minimizer. The projection of the auxiliary point $v$
onto this halfspace brings it closer to $x^\ast$.

\paragraph{Stopping.}
By Assumptions~\ref{asm:components} and~\ref{asm:minimizer},
each computation of $f(p)$ and $\nabla f(p)$ yields the lower bound
\[
 \widetilde f(p)=f(p)-\ip{\nabla f(p)}p-R\norm{\nabla f(p)}\le f^\ast.
\]
The subroutine \textsc{Evaluate} stores the maximum $\widetilde f_{\rm best}$
of such bounds and the best point found so far, $w_{\rm best}\in B(0,R)$.
If $f(w_{\rm best})-\widetilde f_{\rm best}\le\varepsilon/2$,
the method stops and returns $w_{\rm best}$. Component answers that have already been obtained
are not queried again.

Algorithm~\ref{alg:cheap} stores the points $x,v$, the weight $\mathcal A$,
and the parameter $h$, from which the regularization is chosen. The bounds
$\lambda_\downarrow(\mathcal A)$ and $\lambda_\uparrow$ are the lower and upper limits
of the proximal parameter: at each iteration, the method takes
$\lambda=\min\{\max\{h,\lambda_\downarrow(\mathcal A)\},\lambda_\uparrow\}$.
The lower limit keeps the weight of the accelerated step bounded, makes the local
subproblem well conditioned, and ensures stopping; it does not increase as the weight
$\mathcal A$ grows. The upper limit $\lambda_\uparrow$ is the value of the lower limit at the
initial weight $\mathcal A=R^2/(2\Delta)$. Explicit formulas are given
in~\eqref{eq:outer-parameters}. The subroutines \textsc{Evaluate} and
\textsc{Block}, including the rules for counting queries, are given
in Appendix~\ref{app:implementation}. The requirements on the output of
\textsc{Solver} are stated in the next section.

\begin{algorithm}[H]
\caption{YETI}\label{alg:cheap}
\begin{algorithmic}[1]
\Require $n$, $R$, $\Delta$, $\delta$, $\varepsilon$, $L_0$, $L_1$,
constants $L_{0,i}$, a heuristic
\State set $r$, $\lambda_\downarrow$, and $\lambda_{\uparrow}$ by~\eqref{eq:outer-parameters}
\State $x=v=0$
\State $\mathcal A=R^2/(2\Delta)$
\State $h=\lambda_{\uparrow}$
\State the query history is empty
\Loop
 \State $x\gets\Call{Block}{x}$ \Comment{heuristic block}
 \State $\lambda=\min\bigl\{\max\{h,\lambda_\downarrow(\mathcal A)\},\lambda_{\uparrow}\bigr\}$
 \State $a_\lambda=\bigl(1+\sqrt{1+4\lambda\mathcal A}\bigr)/(2\lambda)$
 \State $\alpha=a_\lambda/(\mathcal A+a_\lambda)$
 \State $y=(1-\alpha)x+\alpha v$
 \State \Call{Evaluate}{$y$}
 \State $\ell(y)=\frac1n\sum_i\bigl(L_{0,i}+L_1\norm{\nabla f_i(y)}\bigr)$
 \If{$\ell(y)>e^{L_1r}\bigl(\frac32L_0+2L_1^2(\delta+R^2/\mathcal A)\bigr)$} \Comment{cut}
  \State $\mathcal H_y=\Bigl\{w:\ \ip{\nabla f(y)}{w-y}\le-\dfrac{\ell(y)-\frac32L_0-2L_1^2\delta}{2L_1^2}\Bigr\}$
  \State $v\gets\Pi_{B(0,R)}\Pi_{\mathcal H_y}v$
  \State $h\gets2\lambda$
 \Else \Comment{step based on the solution of the local subproblem}
  \State $(z,g,\beta)\gets\Call{Solver}{y,\lambda}$
  \If{$\norm g\le\varepsilon/(4R)$}
   \State \Return $z$
  \EndIf
  \State $\gamma=\lambda/(\lambda+\beta)$
  \State $x\gets\dfrac{(1-\gamma)\mathcal Ax+\gamma(\mathcal A+a_\lambda)z}{\mathcal A+\gamma a_\lambda}$
  \State $v\gets\Pi_{B(0,R)}(v-\gamma a_\lambda g)$
  \State $\mathcal A\gets\mathcal A+\gamma a_\lambda$
  \State $h\gets\begin{cases}2\lambda,&\beta>0,\\\lambda/2,&\beta=0\end{cases}$
  \If{$\mathcal A\ge2R^2/\varepsilon$}
   \State \Return $x$
  \EndIf
 \EndIf
\EndLoop
\end{algorithmic}
\end{algorithm}

\section{The local subproblem}\label{sec:local-solver}
At the point $y$, we solve the subproblem
\[
 \min_{z\in B(0,R)\cap B(y,r)}
 \left\{f(z)+\frac\lambda2\norm{z-y}^2\right\}.
\]
The solver returns $z\in B(0,R)\cap B(y,r)$,
$g\in\partial(f+I_{B(0,R)})(z)$, and $\beta\ge0$, where
$\beta>0$ is possible only if $\norm{y-z}=r$.
Here $I_{B(0,R)}$ is the indicator function of the ball $B(0,R)$. The output of the solver
must satisfy one of the two conditions
\[
 \norm g\le\frac\varepsilon{4R}
 \quad\text{or}\quad
 \norm{g-(\lambda+\beta)(y-z)}
 \le\frac45(\lambda+\beta)\norm{y-z}.
\]
The first one allows YETI to stop, and the second one ensures the accuracy
of the accelerated update. If $\beta=0$, the weight of this update does not change,
and if $\beta>0$, it is multiplied by the coefficient
$\gamma=\lambda/(\lambda+\beta)<1$.

To solve the subproblem, we use SVRG with component sampling probabilities
proportional to their local smoothness constants.
After each epoch, the above conditions are checked.
The solver terminates almost surely
(Lemma~\ref{lem:local-certificate}), and the expected number of its
queries per iteration of YETI with a fixed history is $O(n\log(n+2))$
(Appendix~\ref{app:upper-total}).
The full algorithm of the solver is given in Appendix~\ref{app:solver-implementation}.

\section{Complexity bound for YETI}\label{sec:cheap-analysis}
The analysis uses the Lyapunov function
\[
 \Phi=\mathcal A\bigl(f(x)-f^\ast\bigr)
       +\tfrac12\norm{v-x^\ast}^2\le R^2.
\]
Accepting a point proposed by the heuristic does not increase the first term.
The projection onto $\mathcal H_y$ decreases the second one. Under the accelerated update,
$\Phi$ does not increase owing to the accuracy condition of the local solver.
Hence $f(x)-f^\ast\le R^2/\mathcal A$.

The number of iterations at which the regularization increases is bounded
owing to the decrease of $\Phi$. At the remaining iterations, the regularization
is halved or the weight $\mathcal A$ grows. This yields a bound on
the number of iterations until the condition $\mathcal A\ge2R^2/\varepsilon$ holds.
The corresponding bounds are proved in Lemmas~\ref{lem:outer-potential}
and~\ref{lem:outer-count} of Appendix~\ref{app:cheap-proofs}.
Each iteration, together with the heuristic block preceding it,
requires $O(n\log N)$ queries in expectation for a fixed history
of the previous iterations.

\begin{theorem}[upper bound]\label{thm:cheap-upper}
Under Assumptions~\ref{asm:components}--\ref{asm:accuracy},
Algorithm~\ref{alg:cheap} stops almost surely and returns
a point $\hat x$ such that
\begin{equation}\label{eq:cheap-upper-final}
 \E[f(\hat x)-f^\ast]\le\varepsilon,
\end{equation}
and the expected total number $T$ of oracle calls
satisfies
\begin{equation}\label{eq:cheap-upper-cost}
 \E T\lesssim T_0\log^2N.
\end{equation}
The count $T$ includes all oracle calls in the heuristic blocks,
the main iterations, and the checks.
\end{theorem}

The proof is given
in Appendix~\ref{app:upper-total}. A variant in which the allowed number of queries
of the heuristic depends on the work already done is considered
in Appendix~\ref{sec:budget}.

\section{Numerical experiments}\label{sec:experiments}

We compare the methods by the number of oracle calls on chains S and X
(the hard families from Section~\ref{sec:lower-origin}),
on Poisson models, and on 65 instances from 13 families (five
instances each), including sums of quartic functions and logistic
regression. The number of queries is divided by $n$, that is, it is measured in passes
over the data; the queries spent on rejected trial points
during line search are also counted.
The relative optimality gap equals $(f-f_{\rm ref})/\widehat\Delta$,
where $\widehat\Delta=f(0)-f_{\rm ref}$ and $f_{\rm ref}$ is the reference
value of the minimum. On the chains, it coincides with $f^\ast$, which is known
from the construction. On the other instances, $f^\ast$ is unknown, and $f_{\rm ref}$ is
its approximation computed numerically in advance; therefore, we denote it by
$f_{\rm ref}$ rather than $f^\ast$. The queries spent on computing the reference value
are not counted in the comparison.
The problem formulations, method settings, and detailed results are given
in Appendix~\ref{app:experiments}.

\begin{figure}[p]
 \centering
 \input{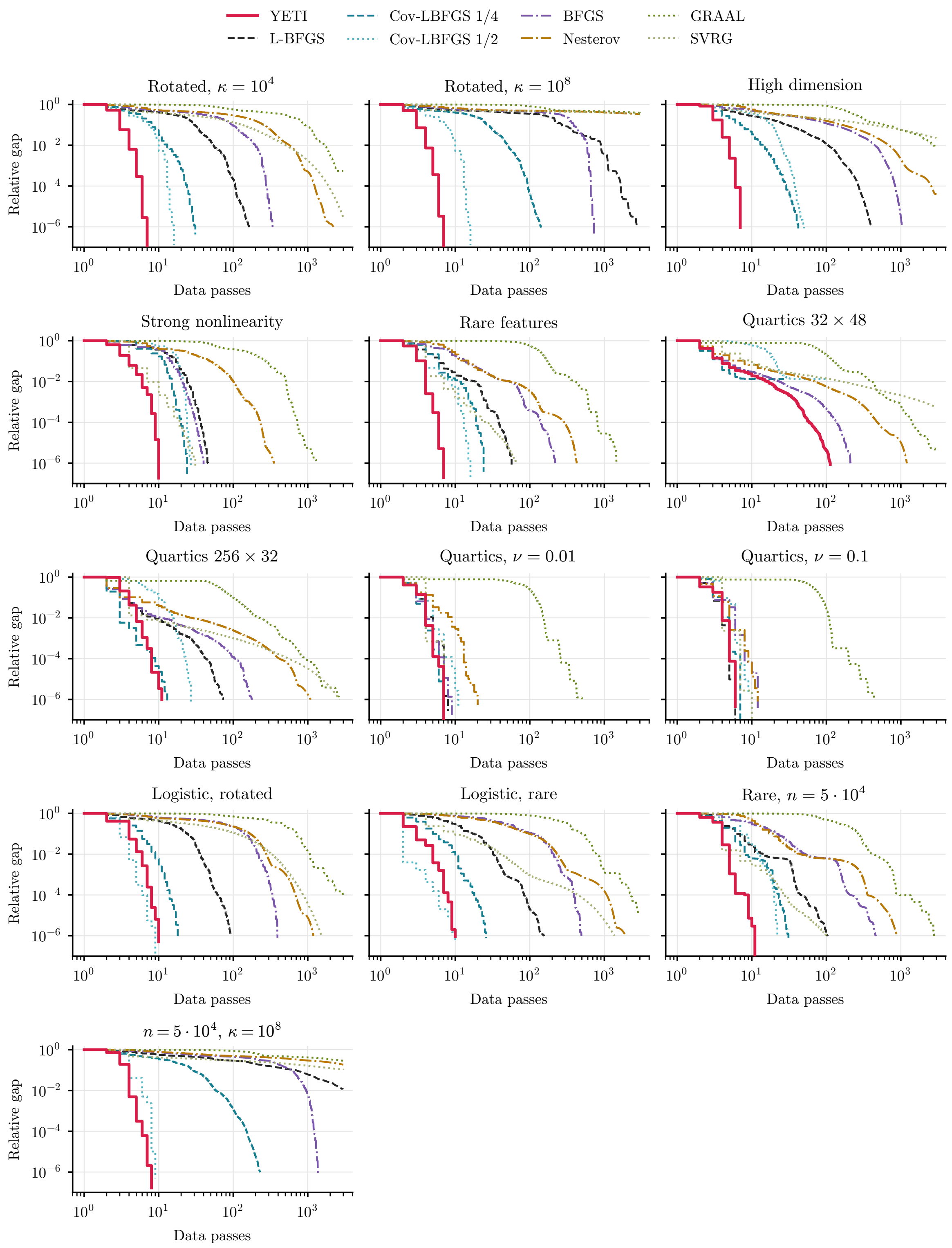}
 \caption{Best relative optimality gap so far on the first instance of each family
 (log--log scale; YETI in red; $\kappa$ is the condition number of the feature covariance).}
 \label{fig:exp-families}
\end{figure}

On the chains, YETI performs approximately as well as L-BFGS. On five
Poisson instances, it reaches the relative accuracy $10^{-4}$
within 4--9 passes; on freMTPL2freq, SAGA requires fewer passes:
5 versus 7 (Tables~\ref{tab:exp-chains}, \ref{tab:exp-glm}).

\paragraph{Stopping.}\label{sec:exp-certificate}
The history-based stopping criterion~\eqref{eq:history-stop} does not use
the reference value. For $\varepsilon=10^{-6}\widehat\Delta$, it is triggered
on all 65 instances, and all returned points satisfy the tolerance
with respect to the reference value. The median number of passes until stopping is 15,
and the maximum is 273. The criterion may be triggered later than the tolerance is actually
reached (Figure~\ref{fig:exp-certificate}).

\paragraph{Test problems.}
Chains S and X are numerical variants of the hard families from
Section~\ref{sec:lower-origin}. The Poisson models are synthetic instances
and instances on the Bike Sharing and freMTPL2freq data. The thirteen families
from which the 65 instances are composed (their names are shown in Figure~\ref{fig:exp-families}
and in Table~\ref{tab:exp-series}) are constructed as follows.
\begin{itemize}
\item \emph{Rotated, $\kappa$}: Poisson regression in which
the feature covariance matrix has condition number
$\kappa\in\{10^4,10^8\}$ and is not diagonal in the original basis (its eigenvectors
are rotated with respect to the coordinate axes).
\item \emph{High dimension}: Poisson regression with a large number of
features: $n=600$, $d=160$, $\kappa=10^4$.
\item \emph{Strong nonlinearity}: Poisson regression in which
the true parameter vector used to generate the responses is three times
longer than in the other Poisson families; this increases the spread
of the linear predictors $a_i^{\mathsf T}x$ and of the values of the exponential.
\item \emph{Rare features} and \emph{Rare, $n=5\cdot10^4$}: Poisson
regression with sparse features; in the second family, the number of components
equals $5\cdot10^4$.
\item \emph{$n=5\cdot10^4$, $\kappa=10^8$}: ill-conditioned
Poisson regression with $5\cdot10^4$ components.
\item \emph{Quartics}: sums of quartic functions
$f_i(x)=\frac12(x-c_i)^{\mathsf T}H_i(x-c_i)+\frac\nu4\norm{x-c_i}^4$,
where $H_i\succ0$ and $c_i$ are random centers. In the family names
$32\times48$ and $256\times32$, the number of components $n$ and the dimension $d$ are indicated;
$\nu$ is the weight of the fourth-degree term (in these two families,
$\nu=10^{-6}$).
\item \emph{Logistic, rotated} and \emph{Logistic, rare}: regularized
logistic regression with rotated and with sparse features,
respectively.
\end{itemize}
\FloatBarrier

\section{Lower bound for arbitrary queries}\label{sec:lower}

We bound from below the number of queries required by any method
that solves all instances of the class under consideration.

\subsection{The origin of the four terms}\label{sec:lower-origin}

We prove the lower bound separately for each of the four terms. Denote
\begin{equation}\label{eq:lower-four-terms}
 \begin{gathered}
 B_n=n,\\
 B_S=\sqrt{\frac{nL_0R^2}{\varepsilon}},\\
 B_Q=n(L_1R)^{2/3},\\
 B_X=\sqrt{\frac{nL_1^2R^2\delta}{\varepsilon}+n(L_1R)^2}.
 \end{gathered}
\end{equation}
Since $\sqrt{a^2+b^2}\sim a+b$ for $a,b\ge0$, we have
$T_0\sim B_n+B_S+B_X+B_Q$.

\begin{center}
\small
\renewcommand{\arraystretch}{1.18}
\begin{tabularx}{\linewidth}{@{}lYY@{}}
\toprule
Term & What makes the problem hard & Hard family\\
\midrule
$B_n$ & Search for an informative component & Hidden pair of functions\\
$B_S$ & Sequential revealing of smooth chains & Independent chains S\\
$B_X$ & Heterogeneity of components under generalized smoothness & Exponential chains X\\
$B_Q$ & Repeated search for the component for the next step & Chain Q with hidden placement of links\\
\bottomrule
\end{tabularx}
\end{center}

For each term, we construct its own family of instances.
Since the method must solve the instances of all four families,
its complexity is not smaller than the maximum of the obtained bounds.
Below, the dependencies are discussed up to numerical and logarithmic
factors. The precise conditions of the constructions are given in the appendix;
if a particular construction is not applicable, the corresponding term
is bounded via the others.

\paragraph{Search for an informative component: $B_n=n$.}
Among the $n$ components, two nonzero ones are chosen at random, and the rest
are zero. Until the first
call to one of them, the algorithm receives only zero answers.
These answers do not allow it to distinguish two instances whose minimizers lie on opposite
sides of the starting point. For the expected error to be small,
one has to find the hidden pair of functions with sufficient probability, and this
requires on the order of $n$ oracle calls
(Lemma~\ref{lem:lower-pair}).

\paragraph{Independent smooth chains: $B_S$.}
Each component contains its own chain of coordinates. By a chain
we mean a function in which the next coordinate affects the oracle answer
only after sufficient progress along the previous
coordinates. In the idealized model, a single query reveals at most
one new link of one chain. How this property carries over to
arbitrary queries is explained below.

In family S, the number of links of a single chain is of order
$\sqrt{L_0R^2/(n\varepsilon)}$. To reach accuracy $\varepsilon$, one needs
to reveal a sufficient fraction of the total number of links of all chains. Therefore,
the number of queries is of order
\[
 n\sqrt{\frac{L_0R^2}{n\varepsilon}}
 =\sqrt{\frac{nL_0R^2}{\varepsilon}}=B_S.
\]
The instances of this family are smooth: $\|\nabla^2f_i\|\le L_0$
(Lemma~\ref{lem:lower-S}). If the conditions for constructing the chain are not satisfied, then
$B_S\lesssim B_n$, and the bound for $B_n$ suffices.

\paragraph{Exponential chains: $B_X$.}
Under generalized smoothness, the curvature of a component can grow together
with its gradient. Even at a minimizer of the average, where
$\nabla f(x^\ast)=0$, the individual gradients $\nabla f_i(x^\ast)$
can be nonzero and cancel each other. Therefore, a small norm of the gradient of the average does not imply
small smoothness constants of the components.
The parameter $\delta$ bounds the discrepancy between the minimum of the average
and the average of the infima of the components.

In family X, the chains enter the exponents
of components grouped into pairs. The construction preserves the sequential
revealing of links and sets the heterogeneity equal to $\delta$.
The bound for this family, together with the bound for the hidden pair of functions, yields $B_X$
(Lemmas~\ref{lem:lower-X} and~\ref{lem:lower-pair}). Since
$B_X=\sqrt n\,L_1R\sqrt{1+\delta/\varepsilon}$, we have
\[
  B_X\sim\sqrt n\,L_1R\max\left\{1,\sqrt{\frac\delta\varepsilon}\right\}.
\]

\paragraph{Search for the next link: $B_Q$.}
Family Q uses a single chain, but its links are randomly
placed across pairs of components. After revealing the next coordinate,
it is not enough to know in which direction to continue: one also has
to find the pair containing the next link. A component that has already
been queried may become informative only at the next level of the
chain, so computing all components once does not eliminate
these searches.

Under the conditions of Lemma~\ref{lem:lower-Q}, the length of the chain is of order
$(L_1R)^{2/3}$, and the search for a pair at each level requires on the order of $n$
queries in expectation. Hence, up to logarithmic factors,
\[
 \underbrace{(L_1R)^{2/3}}_{\text{number of levels}}
 \;\underbrace{n}_{\text{search per level}}
 =B_Q.
\]
The exponent $2/3$ is explained by the choice of the chain length $T$.
At its minimizer, the $T$ active coordinates are equal to $R/\sqrt T$, so
that the norm of this point equals $R$. The bound on the optimality gap in terms of the number of revealed
links~\eqref{eq:lower-Q-gap} shows that to reach the accuracy
one has to pass on the order of $L_1R/\sqrt T$ levels (up to
logarithmic factors), while the total number of levels is $T$. The chain length
is chosen so that the quantities $T$ and $L_1R/\sqrt T$ coincide
up to absolute constants. From the equality $T=L_1R/\sqrt T$
we obtain $T^{3/2}=L_1R$, that is, $T=(L_1R)^{2/3}$; in
Lemma~\ref{lem:lower-Q}, an integer $T$ of this order is chosen.
The precise conditions of this choice also take into account the initial
optimality gap and the heterogeneity. The number of oracle calls required
to find the next link is bounded in Lemmas~\ref{lem:lower-search}
and~\ref{lem:lower-revelations}.

\subsection{The lower bound}

Let Assumptions~\ref{asm:components}--\ref{asm:accuracy} hold.
The number $N_0$ is defined in~\eqref{eq:complexity-parameters}.

\proofdown{thm:lower}
\begin{theorem}\label{thm:lower}
For an oracle returning $(f_i(x),\nabla f_i(x))$, we have
\[
 \mathcal T^\ast\gtrsim\frac{T_0}{\log^2N_0}.
\]
\end{theorem}

The right-hand side does not depend on $\Delta$ for $\Delta\ge2\varepsilon$.
The bound refers to the worst instance of the class in arbitrary dimension.
The method may be randomized and may choose query points anywhere in the
space, not only in the linear span of the previously obtained gradients.

Together with the upper bound for YETI (Theorem~\ref{thm:cheap-upper}),
this theorem shows that in this setting the exponents
of $n$, $L_0$, $L_1$, $R$, $\delta$, and $\varepsilon$ cannot be improved
(up to logarithmic factors). The logarithmic gap between the upper and lower bounds
is not closed here. Assumption~\ref{asm:accuracy} is essential for the scope
of the result: this lower bound cannot be automatically
transferred, for example, to the subclass with $\delta=0$. The equality $\delta=0$ means
that the point $x^\ast$ minimizes each component, that is, all $f_i$ have
a common minimizer. For this case, convergence bounds for the stochastic
methods $(L_0,L_1)$-SGD and SGD-PS were obtained in
\citet[Section~7]{gorbunov2024methods}.

\subsection{Proof of the lower bound}

\begin{proof}[Proof sketch of Theorem~\ref{thm:lower}]
The instance with a hidden pair of functions gives the bound for $B_n$.
Family S gives the bound for $B_S$, and if its conditions
are violated, we use $B_S\lesssim B_n$.
Family X, together with the hidden pair of functions, gives the bound for $B_X$
up to a logarithmic factor. Family Q gives the bound for
$B_Q$. When conditions~\eqref{eq:lower-Q-regime} are not satisfied,
this term is bounded via the bounds already obtained.

Consequently,
\[
 T_0\sim\max\{B_n,B_S,B_X,B_Q\}
 \lesssim\mathcal T^\ast\log^2N_0.
\]
Lemmas~\ref{lem:lower-finite}--\ref{lem:lower-revelations}
make it possible to apply these bounds to arbitrary randomized
methods. The analysis of all parameter regimes is given
in Appendix~\ref{app:lower-constants}.
\end{proof}

\section{Conclusion}\label{sec:conclusion}
For minimizing a sum of convex $(L_0,L_1)$-smooth functions, the method YETI
attains the upper bound $O(T_0\log^2N)$, and for any valid method
we prove the lower bound $\Omega(T_0/\log^2N_0)$. If $\Delta/\varepsilon\le N_0^c$
for some fixed $c\ge1$, then these bounds coincide up to
powers of $\log N_0$. Without a restriction on the growth of the initial optimality gap $\Delta$,
the optimality of the upper bound has not been established.

On 65 instances with relative accuracy $10^{-4}$, the tested variant of YETI
ranked first (including ties with other methods)
in 46 cases, and its history-based stopping criterion was triggered
on all 65 instances. In these experiments, the main gain comes from the heuristic
(Appendix~\ref{app:experiments}). However, its matrix operations may
be expensive, so a gain in the number of queries does not necessarily mean
the same gain in running time.

\subsection{Open questions}\label{sec:open-questions}
\begin{enumerate}
\item \emph{Large initial optimality gap.}
What is the optimal dependence of the complexity on $\Delta$?
If $\log(\Delta/\varepsilon)$ is of order $L_1R$, then the obtained upper
and lower bounds may differ polynomially in $L_1R$. One needs either
to improve the upper bound or to prove a lower bound that accounts for $\Delta$.
For moderate $\Delta$, it remains to close the logarithmic gap.

\item \emph{Individual smoothness constants.}
Can one obtain a sharper complexity bound for a given collection
$(L_{0,i},L_{1,i})$ without replacing all $L_{1,i}$ by their maximum?
Non-uniform sampling of components already makes it possible to use only the average
value of $L_{0,i}$ (Appendix~\ref{app:variants}). To improve the general
bound, one needs to account for the relation between $L_{1,i}$ and the optimality gaps of individual components
and to construct a matching lower bound: our lower bound is obtained
on the subclass with equal $L_{0,i}$.

\item \emph{Constraints and non-Euclidean geometry.}
How can the result be extended to minimization over a convex set
and to methods with a Bregman distance? Mirror descent under generalized
smoothness was studied in~\citet{yu2026mirror}. Here, one will additionally
need an analysis of variance reduction, a replacement of the Euclidean potential,
and a bound on the number of cuts.

\item \emph{Nonsmooth components.}
Can one obtain similar two-sided bounds for finite sums
whose subgradient norms are bounded in terms of the optimality gap?
Such a nonsmooth stochastic setting is considered in
\citet{kovalev2026nonsmooth}. For a finite sum, one needs to define the
oracle model separately and to determine what gain variance reduction provides.

\item \emph{Second- and higher-order methods.}
Can variance reduction for finite sums
\citep{tsedenov2026second} be combined with acceleration under generalized smoothness
of the Hessian \citep{agafonov2026inexact}? This requires error bounds for
stochastic derivatives and a new local solver. Conditions on higher-order
derivatives do not follow from Assumption~\ref{asm:smoothness}.

\item \emph{The dual problem.}
Can YETI be interpreted as a coordinate method for the dual
problem and related to the scheme of \citet{podlipnova2026sparse}?
Besides transforming the problem, one will need to match the available oracles
and the cost of a single step: a query of a component $f_i$ and a query of a coordinate
of the dual gradient are not equivalent in general.
\end{enumerate}

Separately, it remains to determine the complexity outside the scope of
Assumption~\ref{asm:accuracy}, in particular for $\delta=0$.
This question concerns the heterogeneity of the components and is distinct
from the question of a large initial optimality gap $\Delta$.

\begingroup
\selectlanguage{english}
\renewcommand{\refname}{References}
\renewcommand{\bibfont}{\small}
\bibliographystyle{plainnat}
\bibliography{refs}
\endgroup
\clearpage
\appendix
\begingroup
\small
\section{Analytic properties of the class}\label{sec:analytic-tools}

The auxiliary lemmas are stated for an arbitrary function $\varphi$.
They apply to the components $f_i$ by Assumptions~\ref{asm:components}
and~\ref{asm:smoothness}, with the constant $L_{0,i}$ in place of $L_0$.

\subsection{Local smoothness}

See also \citet[Lemma~A.1]{rudomanov2026everest}
and \citet[Lemma~2.5]{vankov2025optimizing}.

\begin{lemma}[growth of the gradient along a segment]\label{lem:local-transfer}
Let $\varphi\in C^2(\mathbb R^d)$ and
$\norm{\nabla^2\varphi(x)}\le L_0+L_1\norm{\nabla\varphi(x)}$
for all $x$. Then for any $x,y$,
\begin{equation}\label{eq:local-transfer}
 L_0+L_1\norm{\nabla\varphi(x)}
 \le e^{L_1\norm{x-y}}
       \bigl(L_0+L_1\norm{\nabla\varphi(y)}\bigr).
\end{equation}
In particular, on the ball $B(y,r)$ the gradient of $\varphi$ is Lipschitz with constant
$e^{L_1r}(L_0+L_1\norm{\nabla\varphi(y)})$.
\end{lemma}

\begin{proof}
By the Newton--Leibniz formula,
\[
 \nabla\varphi(y+t(x-y))-\nabla\varphi(y)
 =\int_0^t\nabla^2\varphi(y+u(x-y))(x-y)\,du.
\]
Taking norms and applying the Hessian condition, for $t\in[0,1]$ we obtain
\[
 L_0+L_1\norm{\nabla\varphi(y+t(x-y))}
 \le L_0+L_1\norm{\nabla\varphi(y)}
 +L_1\norm{x-y}\int_0^t\bigl(L_0+L_1\norm{\nabla\varphi(y+u(x-y))}\bigr)\,du.
\]
By Grönwall's inequality,
\[
 L_0+L_1\norm{\nabla\varphi(y+t(x-y))}
 \le e^{L_1\norm{x-y}t}\bigl(L_0+L_1\norm{\nabla\varphi(y)}\bigr).
\]
Setting $t=1$ gives the first claim. For $x\in B(y,r)$ it yields a bound
on the Hessian norm over the whole ball, and integrating the Hessian along segments
inside the convex ball proves that the gradient is Lipschitz.
\end{proof}

\citep[Baillon--Haddad theorem,][]{bauschke2010}

\begin{lemma}[cocoercivity on a convex domain]\label{lem:local-cocoercivity}
Let $\varphi$ be convex and twice continuously differentiable
in a neighborhood of a convex set $\mathcal D$. If
$\norm{\nabla^2\varphi}\le L$ on $\mathcal D$, where $L>0$, then for all
$x,z\in\mathcal D$,
\[
 \norm{\nabla\varphi(x)-\nabla\varphi(z)}^2
 \le L\ip{\nabla\varphi(x)-\nabla\varphi(z)}{x-z}.
\]
\end{lemma}

\begin{proof}
The averaged Hessian
\[
 H=\int_0^1\nabla^2\varphi(z+t(x-z))\,dt
\]
satisfies $0\preceq H\preceq LI$ by convexity; hence
$H^2\preceq LH$ and
\[
 \norm{\nabla\varphi(x)-\nabla\varphi(z)}^2
 =\ip{x-z}{H^2(x-z)}
 \le L\ip{x-z}{H(x-z)}
 =L\ip{\nabla\varphi(x)-\nabla\varphi(z)}{x-z}.
\]
\end{proof}

\begin{lemma}[bound on the gradient difference]\label{lem:local-bregman}
Let $\varphi$ be convex, $\varphi\in C^2(\mathbb R^d)$, and
$\norm{\nabla^2\varphi}\le L$ on $B(y,3r)$, where $L>0$.
Then for $x,w\in B(y,r)$,
\[
 \norm{\nabla\varphi(x)-\nabla\varphi(w)}^2
 \le2L\bigl(\varphi(x)-\varphi(w)
                     -\ip{\nabla\varphi(w)}{x-w}\bigr).
\]
\end{lemma}

\begin{proof}
Since the gradient of $\varphi$ is $L$-Lipschitz on $B(y,r)$
and $x,w\in B(y,r)$, we have
\[
 \norm{\nabla\varphi(x)-\nabla\varphi(w)}\le L\norm{x-w}\le2Lr.
\]
Therefore, the point
\[
 q=x-\frac1L\bigl(\nabla\varphi(x)-\nabla\varphi(w)\bigr)
\]
lies in $B(y,3r)$, and hence the whole segment $[x,q]$ is contained in $B(y,3r)$.
The function
\[
 z\mapsto\varphi(z)-\varphi(w)-\ip{\nabla\varphi(w)}{z-w}
\]
is convex, nonnegative on the whole space, and has on $B(y,3r)$
the $L$-Lipschitz gradient $\nabla\varphi(z)-\nabla\varphi(w)$.
The quadratic upper bound on the segment $[x,q]$ gives
\begin{align*}
 0&\le\varphi(q)-\varphi(w)-\ip{\nabla\varphi(w)}{q-w}
 \\&\le\varphi(x)-\varphi(w)-\ip{\nabla\varphi(w)}{x-w}
   +\ip{\nabla\varphi(x)-\nabla\varphi(w)}{q-x}+\frac L2\norm{q-x}^2
 \\&=\varphi(x)-\varphi(w)-\ip{\nabla\varphi(w)}{x-w}
   -\frac1{2L}\norm{\nabla\varphi(x)-\nabla\varphi(w)}^2.
\end{align*}
It remains to move the last term to the left-hand side.
\end{proof}

\subsection{Bounding the gradient norm by the optimality gap}

Cf.\ \citet{li2023generalized} and \citet[Lemma~3.4]{yu2026mirror}.

\begin{lemma}[integral bound on the gradient norm]\label{lem:self-bound}
Let $\varphi$ be convex, bounded from below, belong to
$C^2(\mathbb R^d)$, and satisfy the condition of Lemma~\ref{lem:local-transfer}.
Then for any $x$,
\begin{equation}\label{eq:integral-self-bound}
 \begin{gathered}
 \int_0^{\norm{\nabla\varphi(x)}}\frac{t}{L_0+L_1t}\,dt
 \le\varphi(x)-\inf\varphi,\\
 \norm{\nabla\varphi(x)}
 \le L_1\bigl(\varphi(x)-\inf\varphi\bigr)
 +\sqrt{2L_0\bigl(\varphi(x)-\inf\varphi\bigr)}.
 \end{gathered}
\end{equation}
For $L_0=0$, the integrand is extended to zero by continuity.
The infimum is not assumed to be attained.
\end{lemma}

\begin{proof}
If $\nabla\varphi(x)=0$, the claim is obvious. Otherwise, choose $h>0$ such that
\[
 \begin{gathered}
 L_1h\norm{\nabla\varphi(x)}\le\tfrac12,\\
 h\bigl(L_0+L_1\norm{\nabla\varphi(x)}\bigr)\le\tfrac12,
 \end{gathered}
\]
and consider gradient descent started from $x_0=x$:
\[
 x_{j+1}=x_j-h\nabla\varphi(x_j).
\]

\begin{itemize}
\item The gradient norm does not increase. Suppose
$\norm{\nabla\varphi(x_j)}\le\norm{\nabla\varphi(x)}$. The segment
$[x_j,x_{j+1}]$ has length $h\norm{\nabla\varphi(x_j)}$; hence,
by Lemma~\ref{lem:local-transfer}, the Hessian norm on it does not exceed
\begin{equation}\label{eq:self-bound-step-curvature}
 e^{L_1h\norm{\nabla\varphi(x_j)}}
 \bigl(L_0+L_1\norm{\nabla\varphi(x_j)}\bigr)
 \le2\bigl(L_0+L_1\norm{\nabla\varphi(x)}\bigr)
 \le\frac1h.
\end{equation}
Thus, the averaged Hessian $H_j$ on this segment satisfies
$0\preceq hH_j\preceq I$, and
\[
 \begin{gathered}
 \nabla\varphi(x_{j+1})=(I-hH_j)\nabla\varphi(x_j),\\
 \norm{\nabla\varphi(x_{j+1})}\by{\le}{(a)}\norm{\nabla\varphi(x_j)}.
 \end{gathered}
\]
Here (a) follows from $0\preceq I-hH_j\preceq I$.
By induction, these bounds hold for all $j$.

\item We bound two telescoping sums. The quadratic upper bound
together with the curvature estimate~\eqref{eq:self-bound-step-curvature} gives
\[
 \varphi(x_j)-\varphi(x_{j+1}) \ge h\norm{\nabla\varphi(x_j)}^2
 -\frac{h^2}{2}e^{L_1h\norm{\nabla\varphi(x_j)}}
 \bigl(L_0+L_1\norm{\nabla\varphi(x_j)}\bigr)\norm{\nabla\varphi(x_j)}^2 \ge\frac h2\norm{\nabla\varphi(x_j)}^2.
\]
The integrand $t/(L_0+L_1t)$ is nondecreasing, and
\[
 \norm{\nabla\varphi(x_j)}-\norm{\nabla\varphi(x_{j+1})}
 \le\norm{hH_j\nabla\varphi(x_j)}.
\]
Therefore,
\[
 \int_{\norm{\nabla\varphi(x_{j+1})}}^{\norm{\nabla\varphi(x_j)}}
 \frac{t\,dt}{L_0+L_1t} \le\frac{\norm{\nabla\varphi(x_j)}}{L_0+L_1\norm{\nabla\varphi(x_j)}}
 \norm{hH_j\nabla\varphi(x_j)} \by{\le}{(a)}
 he^{L_1h\norm{\nabla\varphi(x_j)}}\norm{\nabla\varphi(x_j)}^2.
\]
Here (a) follows from the bound~\eqref{eq:self-bound-step-curvature}
on $\norm{H_j}$. Since $e^s-1\le2s$ and $e^s\le2$ for $0\le s\le1/2$,
the function decrease bound and the last bound on the integral yield
\[
 \int_{\norm{\nabla\varphi(x_{j+1})}}^{\norm{\nabla\varphi(x_j)}}
 \frac{t\,dt}{L_0+L_1t}
 \le\varphi(x_j)-\varphi(x_{j+1})
 +h^2\bigl(L_0+3L_1\norm{\nabla\varphi(x)}\bigr)\norm{\nabla\varphi(x_j)}^2.
\]
Summing the inequality
\[
 \varphi(x_i)-\varphi(x_{i+1})\ge\frac h2\norm{\nabla\varphi(x_i)}^2
\]
over $i=0,\ldots,m-1$, we obtain
\[
 h\sum_{i<m}\norm{\nabla\varphi(x_i)}^2\le2\bigl(\varphi(x)-\inf\varphi\bigr).
\]
Summing the bound on the difference of integrals over $i=0,\ldots,j-1$, where $0\le j<m$,
and using the resulting bound on the sum of squared gradient norms, we obtain
\[
 \int_0^{\norm{\nabla\varphi(x)}}\frac{t\,dt}{L_0+L_1t}
 \le\bigl(\varphi(x)-\inf\varphi\bigr)
 \Bigl(1+2h\bigl(L_0+3L_1\norm{\nabla\varphi(x)}\bigr)\Bigr)
 +\int_0^{\norm{\nabla\varphi(x_j)}}\frac{t\,dt}{L_0+L_1t}.
\]
Take $m=\lceil h^{-2}\rceil$ and an index $j<m$ with the smallest norm
$\norm{\nabla\varphi(x_j)}$. Then
\[
 \norm{\nabla\varphi(x_j)}^2
 \le\frac{2(\varphi(x)-\inf\varphi)}{mh}\to0
 \quad\text{as }h\downarrow0.
\]
The last integral does not exceed $\norm{\nabla\varphi(x_j)}/L_1$;
hence, passing to the limit, we obtain the first
inequality in~\eqref{eq:integral-self-bound}.

\item For $L_0=0$, the integral equals $\norm{\nabla\varphi(x)}/L_1$, which immediately
gives the second inequality. For $L_0>0$, after the substitution $t=L_0s/L_1$,
the integral bound takes the form
\[
 \frac{L_1\norm{\nabla\varphi(x)}}{L_0}
 -\log\Bigl(1+\frac{L_1\norm{\nabla\varphi(x)}}{L_0}\Bigr)
 \le\frac{L_1^2(\varphi(x)-\inf\varphi)}{L_0}.
\]
The function $s\mapsto s-\log(1+s)$ is increasing for $s\ge0$. For
$v\ge0$, the power series expansion of the exponential gives $e^{\sqrt{2v}}\ge1+\sqrt{2v}+v$; hence
\[
 (v+\sqrt{2v})-\log(1+v+\sqrt{2v})\ge v.
\]
Applying this inequality with $v=L_1^2(\varphi(x)-\inf\varphi)/L_0$, we obtain
\[
 \frac{L_1\norm{\nabla\varphi(x)}}{L_0}
 \le\frac{L_1^2(\varphi(x)-\inf\varphi)}{L_0}
 +\sqrt{\frac{2L_1^2(\varphi(x)-\inf\varphi)}{L_0}},
\]
as required.
\end{itemize}
\end{proof}

\subsection{Bounding the local curvature of the sum}\label{app:observed-curvature}

\begin{proof}[Proof of Lemma~\ref{lem:observed-curvature}\proofup{lem:observed-curvature}]
For the average of the nonnegative quantities $f_i(y)-\inf f_i$ we have
\[
 \frac1n\sum_i\bigl(f_i(y)-\inf f_i\bigr)
 =f(y)-f^\ast+f^\ast-\frac1n\sum_i\inf f_i
 \by{\le}{\text{\ref{asm:heterogeneity}}}\delta+f(y)-f^\ast.
\]
By Lemma~\ref{lem:self-bound}, which applies by
Assumptions~\ref{asm:components} and~\ref{asm:smoothness},
\begin{align*}
 &\frac1n\sum_i\norm{\nabla f_i(y)}\le\frac{L_1}n\sum_i\bigl(f_i(y)-\inf f_i\bigr)
 +\frac1n\sum_i\sqrt{2L_{0,i}\bigl(f_i(y)-\inf f_i\bigr)}
 \\&\qquad\by{\le}{(a)}
 L_1\bigl(\delta+f(y)-f^\ast\bigr)+\sqrt{2L_0\bigl(\delta+f(y)-f^\ast\bigr)},
 \\
 &\ell(y)\le L_0+L_1^2\bigl(\delta+f(y)-f^\ast\bigr)
 +L_1\sqrt{2L_0\bigl(\delta+f(y)-f^\ast\bigr)}
 \\&\qquad\by{\le}{(b)}
 \tfrac32L_0+2L_1^2\bigl(\delta+f(y)-f^\ast\bigr).
\end{align*}
In (a) we used the Cauchy--Schwarz inequality
and the constraint on the average value of $L_{0,i}$ from Assumption~\ref{asm:smoothness};
(b) follows from $2ab\le a^2+b^2$.
The triangle inequality gives both bounds for $\norm{\nabla f(y)}$.
The remaining claims follow by summing the bounds of
Lemma~\ref{lem:local-transfer} and taking Assumption~\ref{asm:smoothness} into account.
\end{proof}

\Needspace{10\baselineskip}
\subsection{Degenerate class}\label{app:degenerate}

\begin{proposition}\label{prop:degenerate-class}
Under Assumptions~\ref{asm:components}--\ref{asm:minimizer}
and~\ref{asm:heterogeneity}, if $L_0=\delta=0$, then all components are constant.
Consequently, the minimax oracle complexity of this class is zero.
\end{proposition}

\begin{proof}
By Assumption~\ref{asm:smoothness}, $L_0=0$ implies
$L_{0,i}=0$. At the point $x^\ast$ from Assumption~\ref{asm:minimizer} we have
\[
 0\le\frac1n\sum_i\bigl(f_i(x^\ast)-\inf f_i\bigr)
 =f^\ast-\frac1n\sum_i\inf f_i\by{\le}{\text{\ref{asm:heterogeneity}}}\delta=0.
\]
Hence each summand is zero, whence $\nabla f_i(x^\ast)=0$.
Lemma~\ref{lem:local-transfer} with $L_{0,i}=0$ gives
\[
 \norm{\nabla f_i(x)}\le e^{L_1\norm{x-x^\ast}}
 \norm{\nabla f_i(x^\ast)}=0
\]
for any $x$. Consequently, all components are constant, and the answer $0$
can be output without a single oracle call.
\end{proof}

\section{Proofs of the results on YETI}\label{app:cheap-proofs}

\subsection{History lower bounds}\label{app:history-proofs}

\begin{proof}[Proof of Lemma~\ref{lem:history-small-gradient}\proofup{lem:history-small-gradient}]
Choose a point $w\in B(0,R)$ minimizing the linear function
$u\mapsto\ip{\nabla f(z)}u$. By the definition of the normal cone,
$\ip{v_z}{w-z}\le0$. Therefore,
\[
 f(z)-\widetilde f(z)=\ip{\nabla f(z)}{z-w}
 \le\ip g{z-w}
 \by{\le}{(a)}\norm g\norm{z-w}
 \le2R\norm g,
\]
where (a) follows from the Cauchy--Schwarz inequality.
After $z$ is recorded in the history, we have $f(w_{\rm best})\le f(z)$ and
$\widetilde f_{\rm best}\ge\widetilde f(z)$, which completes the proof.
\end{proof}

\subsection{A lower bound from the responses of individual components}\label{app:component-history}

The stopping criterion~\eqref{eq:history-stop} uses the lower
bounds~\eqref{eq:history-full-lower}, each of which requires the responses
of all $n$ components at the same point. The heuristic in the experiments
(Appendix~\ref{app:experiments}) queries individual components
at different points. We show that such responses can also be combined into a lower
bound on $f^\ast$ without making new queries. This subsection is needed only
for the experiments: it is not used in the proof of the complexity bound for YETI.

For each index $i$, choose one of the responses of this component
already obtained; let it be obtained at the point $p_i$. Set
\begin{equation}\label{eq:mixed-history-lower}
 \widetilde f_{\rm mix}
 =\frac1n\sum_{i=1}^n\bigl(f_i(p_i)-\ip{\nabla f_i(p_i)}{p_i}\bigr)
 -R\,\Bigl\|\frac1n\sum_{i=1}^n\nabla f_i(p_i)\Bigr\|.
\end{equation}
If all points $p_i$ coincide with some point $p$, then
$\widetilde f_{\rm mix}=\widetilde f(p)$; that is,~\eqref{eq:mixed-history-lower}
generalizes~\eqref{eq:history-full-lower}.

\begin{proposition}\label{prop:mixed-history}
We have $\widetilde f_{\rm mix}\le f^\ast$. Therefore, the quantity
$\widetilde f_{\rm mix}$ can be included in the maximum $\widetilde f_{\rm best}$
on a par with the bounds~\eqref{eq:history-full-lower}. When one
selected response is replaced, recomputing both sums in~\eqref{eq:mixed-history-lower}
requires $O(d)$ arithmetic operations and no new oracle calls.
\end{proposition}

\begin{proof}
By the convexity of the components (Assumption~\ref{asm:components}),
for any $w$,
\[
 f(w)\ge\frac1n\sum_i\bigl(f_i(p_i)+\ip{\nabla f_i(p_i)}{w-p_i}\bigr).
\]
For $w=x^\ast$, Assumption~\ref{asm:minimizer}
and the Cauchy--Schwarz inequality give
$f^\ast\ge\widetilde f_{\rm mix}$. Replacing the response with index $i$
changes only one term of the scalar sum and one term
of the vector sum.
\end{proof}

The bound~\eqref{eq:mixed-history-lower} can be computed only after
at least one response has been obtained for each component. No conditions
on the points $p_i$ or on the order in which they are obtained are required.
Under subsequent replacements of responses, the value $\widetilde f_{\rm mix}$ may
decrease; therefore, $\widetilde f_{\rm best}$ stores the maximum of all
lower bounds obtained so far. To recompute $\widetilde f_{\rm mix}$,
one needs to store one selected response per component; this is
additional memory compared with storing only the best lower
bound and the best point.

\Needspace{12\baselineskip}
\subsection{Local solver}\label{app:local-proofs}

\begin{proof}[Proof of Lemma~\ref{lem:local-projection}\proofup{lem:local-projection}]
\leavevmode
\begin{itemize}
\item For $\zeta\ge0$, consider the point
\[
 z_\zeta=\Pi_{B(0,R)}\Bigl(\frac{b+\zeta y}{1+\zeta}\Bigr).
\]
It depends continuously on $\zeta$ and tends to $y$ as $\zeta\to\infty$.
If $\norm{z_0-y}\le r$, it suffices to take $\zeta=0$. Otherwise, there exists
$\zeta>0$ for which $\norm{z_\zeta-y}=r$. The form of the normal cone
to the ball $B(0,R)$ gives~\eqref{eq:two-ball-normal}. For any $t\in\mathcal K$ we obtain
\[
 \ip{b-z}{t-z}
 =\xi\ip z{t-z}+\zeta\ip{z-y}{t-z}
 \by{\le}{(a)}0.
\]
Here (a) follows from the fact that $\xi z$ and $\zeta(z-y)$ are normals
to the balls $B(0,R)$ and $B(y,r)$ at the point $z$.
This inequality characterizes the projection onto $\mathcal K$ and proves
the claim, including the cases when the balls coincide or one of them is contained in the other.
\item Inequality \eqref{eq:local-bregman} follows from
Lemma~\ref{lem:local-bregman} with $\varphi=P_i$ and $L=L_i$.
The auxiliary point from the proof of
this lemma also lies in the ball of radius $3r$, so the smoothness bound applies to it.
The convexity of $P_i$, which follows from Assumption~\ref{asm:components},
gives the nonnegativity
of its deviation from the tangent. The cocoercivity of $P$ follows from
Lemma~\ref{lem:local-cocoercivity}: on the convex ball $B(y,r)$,
the averaged Hessian of $P$ lies between $0$ and $\overline LI$.
This is a local version of the classical
cocoercivity property~\citep{baillon1977}.
\end{itemize}
\end{proof}

\Needspace{14\baselineskip}
For the local subproblem and the SVRG epoch, we keep the notation $u$, $w$, $s$,
$h_s$, $M$, and $\overline L$ from Section~\ref{sec:local-solver}.

\begin{lemma}[decrease of the Lyapunov function over an epoch]\label{lem:svrg-contraction}
Set
\begin{equation}\label{eq:svrg-potential}
 \Psi(w,s)=\norm{w-u}^2
 +4h_s^2\sum_{j=0}^{M-1}(1-h_s\lambda)^j\,
 \sum_{i=1}^n\frac{\norm{\nabla P_i(s)-\nabla P_i(u)}^2}{nL_i/\overline L}.
\end{equation}
Then, conditionally on the entire history before the start of the epoch,
\begin{equation}\label{eq:svrg-contraction}
 \begin{gathered}
 \E\Psi(w^+,s^+)\le\tfrac12\Psi(w,s),\\
 \Psi(y,y)\le\left(1+\frac{2\overline L}{3\lambda}\right)\norm{y-u}^2.
 \end{gathered}
\end{equation}
\end{lemma}

\begin{proof}
\leavevmode
\begin{itemize}
\item Let the random index $I$ take the value $i$ with probability
$L_i/(n\overline L)$. Then
\[
 v_I-\nabla P(u)
 =\frac{\nabla P_I(x)-\nabla P_I(u)}{L_I/\overline L}
 +\Bigl(\nabla P(s)-\nabla P(u)
 -\frac{\nabla P_I(s)-\nabla P_I(u)}{L_I/\overline L}\Bigr).
\]
The second term is centered, and the expectation of its squared norm
does not exceed the sum from~\eqref{eq:svrg-potential}. For any point
$t\in B(y,r)$,
\begin{equation}\label{eq:svrg-variance-bound}
 \sum_{i=1}^n\frac{\norm{\nabla P_i(t)-\nabla P_i(u)}^2}{nL_i/\overline L}
 \by{\le}{(a)}2\overline L\bigl(P(t)-P(u)-\ip{\nabla P(u)}{t-u}\bigr),
\end{equation}
where (a) follows from the choice of the probabilities $L_i/(n\overline L)$, the bound~\eqref{eq:local-bregman},
and the identity $P=\frac1n\sum_iP_i$. Therefore,
\begin{align*}
 \E\norm{v_I-\nabla P(u)}^2&\by{\le}{(a)}2\sum_{i=1}^n\frac{\norm{\nabla P_i(x)-\nabla P_i(u)}^2}{nL_i/\overline L}+2\sum_{i=1}^n\frac{\norm{\nabla P_i(s)-\nabla P_i(u)}^2}{nL_i/\overline L}
 \\&\by{\le}{(b)}4\overline L\bigl(P(x)-P(u)-\ip{\nabla P(u)}{x-u}\bigr)
 +2\sum_{i=1}^n\frac{\norm{\nabla P_i(s)-\nabla P_i(u)}^2}{nL_i/\overline L}.
\end{align*}
Here (a) follows from the inequality
\[
 \norm{a+b}^2\le2\norm a^2+2\norm b^2,
\]
and (b) from~\eqref{eq:svrg-variance-bound} with $t=x$.
\item The optimality of $u$ means that
$u=\Pi_{\mathcal K}(u-h_s\nabla P(u))$. Therefore,
\begin{align*}
 \E\norm{x^+-u}^2&\by{\le}{(a)}
 \norm{x-u}^2-2h_s\ip{\nabla P(x)-\nabla P(u)}{x-u}
       +h_s^2\E\norm{v_I-\nabla P(u)}^2
 \\&\by{\le}{(b)}
 (1-h_s\lambda)\norm{x-u}^2
 -\frac43h_s\bigl(P(x)-P(u)-\ip{\nabla P(u)}{x-u}\bigr)
 \\&+2h_s^2\sum_{i=1}^n\frac{\norm{\nabla P_i(s)-\nabla P_i(u)}^2}{nL_i/\overline L}.
\end{align*}
Here (a) follows from the nonexpansiveness of the projection and the unbiasedness of $v_I$,
and (b) from the strong convexity of $P$, the previous item, and the identity
$2h_s-4h_s^2\overline L=\frac43h_s$.
\item By \eqref{eq:svrg-parameters}, we have $0<h_s\lambda\le1/6$
and $h_s\lambda M\ge7/10$. We unroll the obtained recursive inequality over the whole epoch.
Since $8h_s^2\overline L=\frac43h_s$, the terms with the linearization
remainders at the points $x_j$ cancel, and
\[
 \E\Psi(w^+,s^+) \by{\le}{(a)}
 (1-h_s\lambda)^M\norm{w-u}^2
 +2h_s^2\sum_{j=0}^{M-1}(1-h_s\lambda)^j\,\sum_{i=1}^n\frac{\norm{\nabla P_i(s)-\nabla P_i(u)}^2}{nL_i/\overline L} \by{\le}{(b)}\tfrac12\Psi(w,s).
\]
Here (a) follows from the recursive inequality of the previous item,
the bound~\eqref{eq:svrg-variance-bound} with $t=s^+$, and Jensen's inequality
for the convex function $t\mapsto P(t)-P(u)-\ip{\nabla P(u)}{t-u}$,
and (b) from $(1-h_s\lambda)^M\le e^{-7/10}<1/2$.
Finally, local smoothness gives
\[
 P(y)-P(u)-\ip{\nabla P(u)}{y-u}\le\frac{\overline L}2\norm{y-u}^2,
\]
so, by~\eqref{eq:svrg-variance-bound}, the sum
from~\eqref{eq:svrg-potential} with $s=y$ does not exceed $\overline L^2\norm{y-u}^2$.
Together with
\[
 \sum_{j=0}^{M-1}(1-h_s\lambda)^j\le\frac1{h_s\lambda}
\]
this gives
\[
 \Psi(y,y)\le\left(1+\frac{4h_s\overline L^2}{\lambda}\right)\norm{y-u}^2
 =\left(1+\frac{2\overline L}{3\lambda}\right)\norm{y-u}^2.
\]
\end{itemize}
\end{proof}

\begin{proof}[Proof of Lemma~\ref{lem:local-certificate}\proofup{lem:local-certificate}]
\leavevmode
\begin{itemize}
\item By Lemma~\ref{lem:local-projection}, the map
$I-\nabla P/\overline L$ is nonexpansive on $\mathcal K$. Comparing
the projections of the points $w-\nabla P(w)/\overline L$ and $u-\nabla P(u)/\overline L$,
we obtain the first bound in \eqref{eq:sharp-residual}. For the second bound, we have
\begin{align*}
 \norm e^2&=\overline L^2\norm{w-z}^2
 -2\overline L\ip{\nabla P(w)-\nabla P(z)}{w-z}
 +\norm{\nabla P(w)-\nabla P(z)}^2
 \\&\by{\le}{(a)}\overline L\ip e{w-z}
 \by{\le}{(b)}\overline L\ip e{w-u}
 \by{\le}{(c)}\overline L\norm e\norm{w-u}.
\end{align*}
Here (a) follows from the cocoercivity of $\nabla P$, (b) from the inequality
\[
 \ip e{z-u}\ge P(z)-P(u)\ge0,
\]
which holds since $e\in\nabla P(z)+N_{\mathcal K}(z)$
by the definition of the projection, and (c) from the
Cauchy--Schwarz inequality.
\item If $\norm{w-u}$ satisfies the condition of the lemma, then
\[
 \norm e\le\overline L\norm{w-u}
 \le\frac45\lambda\bigl(\norm{y-u}-\norm{w-u}\bigr)
 \le\frac45\lambda\norm{y-z}
 \le\frac45(\lambda+\beta)\norm{y-z}.
\]
The second-to-last step follows from $\norm{z-u}\le\norm{w-u}$
and the triangle inequality. Hence the second
condition in~\eqref{eq:certificate-test} holds.
If $y=u$, then $\Psi(y,y)=0$, after the first epoch $w=u$ almost surely,
and the triple is accepted. If $y\ne u$, then after $j$ epochs
Lemma~\ref{lem:svrg-contraction} and Markov's inequality give
\[
 \begin{aligned}
 \Pr(T_{\rm ep}>j)
 &\le\Pr\left(\norm{w_j-u}>\frac{4\lambda}{5\overline L+4\lambda}\norm{y-u}\right)\\
 &\le\left(1+\frac{5\overline L}{4\lambda}\right)^2
 \frac{\E\Psi(w_j,s_j)}{\norm{y-u}^2}
 \le\Gamma2^{-j}.
 \end{aligned}
\]
In deriving this bound, the epochs can be conceptually continued even after the first
acceptance, using independent samples. The sum of the right-hand sides is finite, which proves
both claims about the stopping time.
\end{itemize}
\end{proof}

\subsection{Decrease of the Lyapunov function}\label{app:potential-proof}

\begin{lemma}[decrease of the Lyapunov function]\label{lem:outer-potential}
In all states of the method, the following conditions hold:
\begin{equation}\label{eq:outer-invariant}
 \begin{gathered}
 x,v\in B(0,R),\\
 \mathcal A\ge R^2/(2\Delta),\\
 h\ge\tfrac12\lambda_\downarrow(\mathcal A),\\
 \mathcal A\bigl(f(x)-f^\ast\bigr)+\tfrac12\norm{v-x^\ast}^2\le R^2.
 \end{gathered}
\end{equation}
Denote the left-hand side of the last inequality by $\Phi$
(the Lyapunov function). We always have $\alpha\le\frac12(r/R)^{2/3}\le\frac12$ and $\lambda\le2h$.
Accepting a point proposed by the heuristic block does not change $v,\mathcal A,h$
and does not increase $\Phi$.
If \eqref{eq:simple-curvature-test} is violated, the halfspace
$\mathcal H_y$ contains $x^\ast$ and
\begin{equation}\label{eq:cut-descent}
 \Phi^+\le\Phi-\frac{9r^2}{128\alpha^2}.
\end{equation}
Let $y=(1-\alpha)x+\alpha v$ be the center of the local subproblem
from~\eqref{eq:outer-step-setup}, and let $z$ be the point returned by the solver.
If the triple $(z,g,\beta)$ is accepted by the second condition in~\eqref{eq:certificate-test}, then
\begin{equation}\label{eq:damped-descent}
 \Phi^+\le\Phi-\frac9{50}\frac{\norm{y-z}^2}{\alpha^2}.
\end{equation}
Under any mode of stopping, the optimality gap of the output does not exceed
$\varepsilon/2$.
\end{lemma}

\begin{proof}
\leavevmode
\begin{itemize}
\item At the initial moment,
\[
 \Phi_0=\frac{R^2}{2\Delta}\bigl(f(0)-f^\ast\bigr)
          +\tfrac12\norm{x^\ast}^2
 \by{\le}{\text{\ref{asm:minimizer}, \ref{asm:initial-gap}}}R^2.
\]
The function $\lambda_\downarrow$ is nonincreasing; hence
$\lambda_\downarrow(\mathcal A)\le\lambda_{\uparrow}$ and
$\lambda\ge\lambda_\downarrow(\mathcal A)$.
The identity $\lambda\mathcal A\alpha^2=1-\alpha$ implies
$\alpha\le\frac12(r/R)^{2/3}$, since the function $t\mapsto(1-t)/t^2$
is decreasing. The inequality $h\ge\lambda_\downarrow(\mathcal A)/2$ implies
$\lambda\le2h$.
\item By the rule for accepting a heuristic point,
$z\in B(0,R)$ and $f(z)\le f(x)$. Replacing $x$ by $z$ changes the Lyapunov function by
$\mathcal A(f(z)-f(x))\le0$, while $v,\mathcal A,h$ do not change.
Therefore, conditions~\eqref{eq:outer-invariant} are preserved.
\item From Assumption~\ref{asm:components}
and Corollary~\ref{cor:curvature-lower-gap} we obtain
\[
 \ip{\nabla f(y)}{x^\ast-y}\le f^\ast-f(y)
 \le-\frac{\ell(y)-\frac32L_0-2L_1^2\delta}{2L_1^2},
\]
that is, $x^\ast\in\mathcal H_y$. Suppose condition~\eqref{eq:simple-curvature-test}
is violated. Since $L_1r\le1/2$,
\[
 \frac{\ell(y)-\frac32L_0-2L_1^2(\delta+R^2/\mathcal A)}{2L_1^2}
 \by{>}{(a)}\frac{(1-e^{-L_1r})\ell(y)}{2L_1^2}
 \by{\ge}{(b)}\frac{3r\ell(y)}{8L_1}>0.
\]
Here (a) follows from the violation of~\eqref{eq:simple-curvature-test},
and (b) from $1-e^{-s}\ge3s/4$ for $0\le s\le1/2$.
Since $f(x)-f^\ast\le\Phi/\mathcal A\le R^2/\mathcal A$,
\[
 \ip{\nabla f(y)}{y-x}\by{\ge}{\text{\ref{asm:components}}}f(y)-f(x)
 \ge f(y)-f^\ast-\frac{R^2}{\mathcal A}
 \by{\ge}{(b)}\frac{\ell(y)-\frac32L_0-2L_1^2(\delta+R^2/\mathcal A)}{2L_1^2}.
\]
Step (b) follows from Corollary~\ref{cor:curvature-lower-gap}.
The right-hand side is positive; consequently, $\nabla f(y)\ne0$. Using
$v-y=(1-\alpha)(y-x)/\alpha$ and $\norm{\nabla f(y)}\le\ell(y)/L_1$,
we obtain
\begin{align*}
 \dist(v,\mathcal H_y)&=\frac1{\norm{\nabla f(y)}}\left(\ip{\nabla f(y)}{v-y}
   +\frac{\ell(y)-\frac32L_0-2L_1^2\delta}{2L_1^2}\right)
 \\&\ge\frac{\ell(y)-\frac32L_0-2L_1^2(\delta+R^2/\mathcal A)}
          {2L_1^2\,\alpha\norm{\nabla f(y)}}
 \ge\frac{3r}{8\alpha}.
\end{align*}
The projection onto $\mathcal H_y$ decreases the squared distance to $x^\ast$
by at least the squared distance from $v$ to $\mathcal H_y$, and the projection onto $B(0,R)$ does not increase it,
by Assumption~\ref{asm:minimizer}.
This gives \eqref{eq:cut-descent}. The argument does not require $x$
to belong to the halfspace.
\item Suppose the triple is accepted by the second condition in~\eqref{eq:certificate-test}.
Then $g$ is a subgradient of $f+I_{B(0,R)}$ at the point $z$, and
\begin{align*}
 \Phi^+-\Phi&\by{\le}{(a)}
 \gamma\mathcal A\ip g{z-x}+\gamma a_\lambda\ip g{z-v}
       +\tfrac12\gamma^2a_\lambda^2\norm g^2 =-\gamma(\mathcal A+a_\lambda)\ip g{y-z}
       +\tfrac12\gamma^2a_\lambda^2\norm g^2
 \\&=\tfrac12\gamma^2a_\lambda^2
       \bigl(\norm e^2-(\lambda+\beta)^2\norm{y-z}^2\bigr) \by{\le}{(b)}
 -\frac9{50}\gamma^2a_\lambda^2(\lambda+\beta)^2\norm{y-z}^2
 =-\frac9{50\alpha^2}\norm{y-z}^2.
\end{align*}
In (a) we used the convexity of $f$ (Assumption~\ref{asm:components}),
the combination~\eqref{eq:damped-update}, and the nonexpansiveness of the projection,
which preserves $x^\ast$ by Assumption~\ref{asm:minimizer};
(b) follows from $\norm e\le\frac45(\lambda+\beta)\norm{y-z}$.
In the last equality we used the identities $\gamma(\lambda+\beta)=\lambda$ and
$\lambda a_\lambda=1/\alpha$.
If $\beta>0$, the displacement equals $r$. All new points belong to
$B(0,R)$, and $\mathcal A$ does not decrease. On an up-step,
$h^+=2\lambda\ge\lambda_\downarrow(\mathcal A^+)/2$;
on a down-step, $h^+=\lambda/2\ge\lambda_\downarrow(\mathcal A^+)/2$.
Therefore, conditions~\eqref{eq:outer-invariant} are preserved.
\item If a triple with $\norm g\le\varepsilon/(4R)$ is accepted, then for
$z,x^\ast\in B(0,R)$,
\[
 f(z)-f^\ast
 \by{\le}{(a)}\ip g{z-x^\ast}
 \by{\le}{(b)}2R\norm g\le\varepsilon/2,
\]
where (a) follows from $g\in\partial(f+I_{B(0,R)})(z)$, and (b) from
the Cauchy--Schwarz inequality
and Assumption~\ref{asm:minimizer}.
Stopping due to the growth of $\mathcal A$ gives
$f(x)-f^\ast\le R^2/\mathcal A\le\varepsilon/2$.
By~\eqref{eq:history-full-lower} and Proposition~\ref{prop:mixed-history},
we have $\widetilde f_{\rm best}\le f^\ast$. Therefore, when stopping
by criterion~\eqref{eq:history-stop},
\[
 f(w_{\rm best})-f^\ast
 \le f(w_{\rm best})-\widetilde f_{\rm best}\le\varepsilon/2.
\]
Finally, a triple accepted by the second condition with $\beta=0$ and
$\lambda=\varepsilon/(8rR)$ already satisfies
\[
 \norm g\le\frac95\cdot\frac{\varepsilon}{8rR}\,r=\frac{9\varepsilon}{40R}<\frac{\varepsilon}{4R}.
\]
Therefore, no separate stopping rule is needed for $\lambda=\varepsilon/(8rR)$.
\end{itemize}
\end{proof}

\subsection{Number of main iterations}\label{app:count-proof}

\begin{lemma}[number of main iterations]\label{lem:outer-count}
For any choice of the heuristic in Algorithm~\ref{alg:cheap},
the method almost surely stops
after at most $K_{\rm out}$ main iterations, where
$K_{\rm out}$ is a deterministic quantity satisfying the bound
\begin{equation}\label{eq:outer-budget}
 \begin{aligned}
 K_{\rm out}\lesssim{}
 &\left(\frac Rr\right)^{2/3}\left(1+\log\frac\Delta\varepsilon\right)
 +\log\left(2+\frac{rR\lambda_{\uparrow}}\varepsilon\right)\\
 &+\frac{L_1R}{\sqrt n}\log\frac\Delta\varepsilon
 +\sqrt{\frac{(L_0+L_1^2\delta)R^2}{n\varepsilon}}.
 \end{aligned}
\end{equation}
\end{lemma}

\begin{proof}
We count only the iterations after which the method has not yet stopped;
on them, $\mathcal A<2R^2/\varepsilon$.
Accepting the output of the heuristic does not change $\mathcal A,h$
and does not increase the Lyapunov function $\Phi$.

\begin{itemize}
\item Let $n_U$ be the number of iterations that increase the regularization.
By Lemma~\ref{lem:outer-potential}, on each of them
\[
 \Phi-\Phi^+\ge\frac{9r^2}{128\alpha^2}
 \ge\frac9{32}R^2\left(\frac rR\right)^{2/3},
 \qquad \alpha\le\frac12\left(\frac rR\right)^{2/3}.
\]
Since $0\le\Phi\le R^2$, we have $n_U\le\frac{32}{9}(R/r)^{2/3}$.
Each such iteration increases $h$ by at most a factor of four.
An iteration with $\beta=0$ and $\lambda>\lambda_\downarrow(\mathcal A)$
decreases $h$ by at least a factor of two, and for
$\lambda=\lambda_\downarrow(\mathcal A)$ it does not increase $h$. From
$h\ge\lambda_\downarrow(\mathcal A)/2\ge\varepsilon/(16rR)$
it follows that the number $n_D$ of iterations with $\beta=0$ and
$\lambda>\lambda_\downarrow(\mathcal A)$ satisfies
\[
 h\,2^{n_D}\le\lambda_{\uparrow}4^{n_U},\qquad
 n_D\le2n_U+\log_2\frac{16rR\lambda_{\uparrow}}\varepsilon.
\]

\item Consider the remaining iterations with $\beta=0$ and
$\lambda=\lambda_\downarrow(\mathcal A)$.
On them, $\mathcal A^+=\mathcal A+a_\lambda$.
If the maximum in~\eqref{eq:outer-parameters} is attained at the first term,
then $\alpha=\frac12(r/R)^{2/3}$ and
\[
 \log\frac{\mathcal A^+}{\mathcal A}
 =-\log(1-\alpha)\ge\alpha=\frac12(r/R)^{2/3}.
\]
Since $\mathcal A_0=R^2/(2\Delta)$ and, before stopping,
$\mathcal A<2R^2/\varepsilon$, the total growth of $\log\mathcal A$
is less than $\log(4\Delta/\varepsilon)$. Therefore, there are at most
$2(R/r)^{2/3}\log(4\Delta/\varepsilon)$ such iterations.

\item If the maximum is attained at the second term, then
\[
 \lambda=\frac1n\left(\frac32L_0+2L_1^2\delta
                       +\frac{2L_1^2R^2}{\mathcal A}\right).
\]
If $\frac32L_0+2L_1^2\delta\ge2L_1^2R^2/\mathcal A$, then
$\lambda\le(3L_0+4L_1^2\delta)/n$.
From $a_\lambda\ge\sqrt{\mathcal A/\lambda}$ and
$\mathcal A^+\le2\mathcal A$ we obtain
\[
 \begin{aligned}
 \sqrt{\mathcal A^+}-\sqrt{\mathcal A}
 &=\frac{a_\lambda}{\sqrt{\mathcal A+a_\lambda}+\sqrt{\mathcal A}}\\
 &\ge\frac1{(1+\sqrt2)\sqrt\lambda}
 \ge\frac{\sqrt n}{(1+\sqrt2)\sqrt{3L_0+4L_1^2\delta}}.
 \end{aligned}
\]
Since $\sqrt{\mathcal A}<\sqrt{2R^2/\varepsilon}$, the number of these iterations
does not exceed
\[
 (1+\sqrt2)\sqrt{\frac{2(3L_0+4L_1^2\delta)R^2}{n\varepsilon}}.
\]
In the opposite case, $\lambda\mathcal A<4(L_1R)^2/n$.
Since $\alpha\le1/2$, we have
\[
 \alpha^2=\frac{1-\alpha}{\lambda\mathcal A}
 \ge\frac n{8(L_1R)^2},\qquad
 \log\frac{\mathcal A^+}{\mathcal A}\ge\alpha
 \ge\frac{\sqrt n}{2\sqrt2\,L_1R}.
\]
Therefore, the number of such iterations is
$O((L_1R/\sqrt n)\log(\Delta/\varepsilon))$.

\item If the maximum is attained at the last term, that is,
$\lambda=\varepsilon/(8rR)$, then the method stops
by Lemma~\ref{lem:outer-potential}.
If several terms are equal, the iteration is assigned to any one of them.
\end{itemize}
Summing the obtained bounds and accounting for the last iteration, we obtain~\eqref{eq:outer-budget}.
Each solver call almost surely terminates after finitely many steps,
and each heuristic block is finite by definition
(Section~\ref{sec:heuristic-block}).
Consequently, the method
almost surely stops.
\end{proof}

\subsection{Number of heuristic queries}\label{app:budget-proof}

\begin{proof}[Proof of Theorem~\ref{thm:adaptive-cheap-budget}\proofup{thm:adaptive-cheap-budget}]
\leavevmode
\begin{itemize}
\item
Initially, $G=H=0$. The inequality $H\le B_0+\omega G$ implies
that the number of allowed queries in~\eqref{eq:adaptive-cheap-budget} is nonnegative.
A single block spends at most $B_0+\omega G-H$ queries, so after the block
the inequality is preserved. During a main iteration, only $G$ increases,
so the inequality is preserved as well. The limit of each block is at most $B_0$, which
proves the second bound $H\le kB_0$.

\item
Blocks with limit $B_k\le B_0$ satisfy the requirements of
Section~\ref{sec:heuristic-block};
therefore, Lemmas~\ref{lem:outer-potential}
and~\ref{lem:outer-count} apply to them.
All conditional bounds for the main iterations and the bound on the number of oracle calls
in each heuristic block remain valid. Consequently, the bound of
Theorem~\ref{thm:cheap-upper} is preserved. Finally, for every run,
\[
 T=G+H\by{\le}{(a)}
 B_0+(1+\omega)G,
\]
where (a) follows from~\eqref{eq:budget-prefix-invariant}.
Taking expectations, we obtain~\eqref{eq:adaptive-budget-total}.
\end{itemize}
\end{proof}

\section{Proof of the upper bound}\label{app:upper-constants}
We use $N$ and $T_0$ from Section~\ref{sec:setting}.
All bounds below are uniform in the problem parameters and in the history of the method.

\subsection{Number of main iterations}
By the choice of the radius $r$, we have
\[
 \frac Rr=\max\{1,2L_1R\}\le2\sqrt N.
\]
From~\eqref{eq:outer-parameters} it follows that
\begin{equation}\label{eq:lambda-range-bound}
 \begin{aligned}
 \frac{\lambda_{\uparrow}rR}{\varepsilon}
 &\lesssim
 \frac\Delta\varepsilon\left(\frac Rr\right)^{1/3}
 +\frac{(L_0+L_1^2\delta)R^2}{n\varepsilon}
 +\frac\Delta\varepsilon\frac{L_1R}{n}+1
 \lesssim N^{3/2}.
 \end{aligned}
\end{equation}
Therefore, the logarithms in~\eqref{eq:outer-budget} are of order $O(\log N)$.

\begin{lemma}\label{lem:outer-numerical}
The number of main iterations satisfies
\begin{equation}\label{eq:outer-13}
 nK_{\rm out}\lesssim T_0\log N.
\end{equation}
\end{lemma}
\begin{proof}
Substituting~\eqref{eq:lambda-range-bound} into~\eqref{eq:outer-budget}, we obtain
\[
 nK_{\rm out}\lesssim\log N\left[
 n\left(\frac Rr\right)^{2/3}
 +\sqrt n\,L_1R
 +\sqrt{\frac{n(L_0+L_1^2\delta)R^2}{\varepsilon}}\right].
\]
It remains to use $(R/r)^{2/3}\sim1+(L_1R)^{2/3}$
and $a+b\sim\sqrt{a^2+b^2}$ for $a,b\ge0$.
\end{proof}

\subsection{Number of local solver queries}
The curvature test~\eqref{eq:simple-curvature-test} and the lower bound
on $\lambda$ give
\begin{equation}\label{eq:condition-number-simple}
 \kappa:=\frac{\overline L}{\lambda}
 =1+\frac{e^{3L_1r}\ell(y)}\lambda
 \le1+ne^{4L_1r}\le1+ne^2\le8n,
\end{equation}
since $L_1r\le1/2$ and $n\ge2$. Here we used the fact that the same quantity
$\frac32L_0+2L_1^2(\delta+R^2/\mathcal A)$ enters both the curvature test
and the lower bound on $n\lambda$.

\begin{lemma}\label{lem:epoch-numerical}
For any history preceding a local call,
\begin{equation}\label{eq:epoch-37}
 \E[T_{\rm ep}\mid\text{history}]\le2+\lceil\log_2\Gamma\rceil,
 \qquad n+M\le36n,
\end{equation}
where $\Gamma$ is defined in Lemma~\ref{lem:local-certificate} and
$\Gamma\le9(8n)^3$.
In particular, the conditional expectation of the number of local solver
queries is of order $O(n\log(n+2))$.
\end{lemma}
\begin{proof}
By~\eqref{eq:condition-number-simple}, the epoch length is
$M=\lceil21\kappa/5\rceil\le1+168n/5\le35n$.
For the quantity $\Gamma$ from~\eqref{eq:epoch-tail} we have
\[
 \Gamma=\left(1+\frac{2\kappa}{3}\right)
         \left(1+\frac{5\kappa}{4}\right)^2
 \le\frac{5\kappa}{3}\left(\frac{9\kappa}{4}\right)^2
 \le9\kappa^3\le9(8n)^3,
\]
since $\kappa\ge1$. Set $j_0=\lceil\log_2\Gamma\rceil$.
Summing the tail bound~\eqref{eq:epoch-tail}, we obtain
\[
 \E[T_{\rm ep}\mid\text{history}]
 \le1+\sum_{j\ge1}\min\{1,\Gamma2^{-j}\}
 \le1+j_0+\sum_{j>j_0}\Gamma2^{-j}\le j_0+2.
\]
Each epoch, together with the verification of its result, requires $O(n+M)=O(n)$ queries.
\end{proof}

\subsection{Total number of queries}\label{app:upper-total}
\begin{proof}[Proof of Theorem~\ref{thm:cheap-upper}]
By Lemmas~\ref{lem:outer-potential} and~\ref{lem:outer-count},
the method almost surely stops with an optimality gap of at most
$\varepsilon/2$. This proves~\eqref{eq:cheap-upper-final}.

If the curvature test fails, a main iteration requires $n$ queries.
Otherwise, by Lemma~\ref{lem:epoch-numerical}, the conditional expected number of queries
is $O(n\log(n+2))$. This bound includes the evaluations at the center,
at the anchor points, and in the verification of the solver's output.
The history-based stopping can only decrease the number of queries:
for comparison, the solver can be continued with fresh random indices
until its own stopping criterion is met. By the definition in
Section~\ref{sec:heuristic-block}, each heuristic block
makes $O(n\log N)$ queries.
Since $n+2\le N$, the same bound applies to a main iteration as well.

The number of main iterations and of the blocks preceding them does not exceed
$K_{\rm out}$. Summing the conditional bounds and using~\eqref{eq:outer-13}, we obtain
\[
 \E T\lesssim nK_{\rm out}\log N\lesssim T_0\log^2N.
\]
The bound is uniform over all heuristics admissible
in the block from Section~\ref{sec:heuristic-block}.
\end{proof}

\section{Constructions for the lower bound}\label{app:lower-proofs}

Throughout this appendix, Assumptions~\ref{asm:components}--\ref{asm:accuracy} are assumed to hold.

\subsection{Extension to randomized methods}

\begin{lemma}[averaging over a finite family]\label{lem:lower-finite}
Let $\mathcal F$ be a nonempty finite family of admissible inputs,
and let $c_\tau>0$ and $c_{\mathcal E}\ge0$ be weight coefficients.
Denote by $\tau_F$ the number of oracle calls made by the algorithm on input $F$,
and by $\mathcal E_F$ the optimality gap of its output on this input,
that is, $f_F(\hat x_F)-f_F^\ast$.
If for every deterministic algorithm
\[
 \frac1{|\mathcal F|}\sum_{F\in\mathcal F}
 \bigl(c_\tau\tau_F+c_{\mathcal E}\mathcal E_F\bigr)\ge c,
\]
then $\mathcal T^\ast\ge(c-c_{\mathcal E}\varepsilon)/c_\tau$.
Inputs in the family may repeat: the sum is taken over a finite set
of parameters with equal weights.
\end{lemma}

\begin{proof}[Proof of Lemma~\ref{lem:lower-finite}]
If all random variables used by the algorithm are fixed,
it becomes deterministic. Measurability of the query selection rules
in the adopted oracle model allows us to average the bound
over the random variables of the algorithm.
Let $\mathsf M$ be a valid algorithm (that is, its expected
optimality gap on every input does not exceed $\varepsilon$) with $\mathcal T(\mathsf M)<+\infty$.
Then
\[
 c \by{\le}{(a)}
 \mathbb E\!\left[
 \frac1{|\mathcal F|}\sum_{F\in\mathcal F}
 (c_\tau\tau_F+c_{\mathcal E}\mathcal E_F)\right] =
 \frac1{|\mathcal F|}\sum_{F\in\mathcal F}
 \bigl(c_\tau\mathbb E\tau_F+
 c_{\mathcal E}\mathbb E\mathcal E_F\bigr) \by{\le}{(b)}
 c_\tau\mathcal T(\mathsf M)+c_{\mathcal E}\varepsilon.
\]
Here inequality (a) follows from the hypothesis of the lemma applied for every value of the random variables of the algorithm,
and (b) uses the inequalities $\mathbb E\tau_F\le\mathcal T(\mathsf M)$
and $\mathbb E\mathcal E_F\le\varepsilon$.
Hence $\mathcal T(\mathsf M)\ge(c-c_{\mathcal E}\varepsilon)/c_\tau$;
for $\mathcal T(\mathsf M)=+\infty$ this inequality is obvious.
Taking the infimum over valid algorithms, we obtain the claim
for $\mathcal T^\ast$. This is a version of Yao's minimax principle~\citep{yao1977} for a finite family.
\end{proof}

Let $c:\{1,\ldots,s\}\to\{1,\ldots,d\}$ be chosen uniformly,
and let coordinate $0$ in $\R^{d+1}$ be distinguished. The revelation of coordinates and the oracle answers
depend only on the current query and on the already revealed values,
while the queries and the returned point depend on the history of answers.
When bounding the probability of revealing coordinates, we also take into account
the returned point; the number of queries does not increase in doing so.

\Needspace{14\baselineskip}
\begin{lemma}[averaging over a hidden coordinate]\label{lem:lower-resample}
Let $\alpha,L>0$, and let $N_{\rm qry}\ge1$ be an arbitrary integer.
The fraction of maps $c$ for which
among the first $N_{\rm qry}$ query points and the point returned
before this deadline there is a point $x$ with
$\|x\|\le L$ and $|x_{c(j)}|\ge\alpha$ for some
not yet revealed coordinate $c(j)$ does not exceed
$(N_{\rm qry}+1)sL^2/(d\alpha^2)$. Taking non-injective maps into account, the fraction
of exceptions does not exceed
\begin{equation}\label{eq:lower-pbad}
 \frac{s^2+(N_{\rm qry}+1)sL^2/\alpha^2}{d}.
\end{equation}
\end{lemma}

Bound~\eqref{eq:lower-pbad} refers to the uniform distribution
on the set of all maps, not to the conditional distribution of injections.
To each non-injective map we assign the input corresponding to
one fixed injection. Thus, every input of the finite
distribution is admissible. On the exceptions we use only the nonnegativity
of the number of queries and of the optimality gap; no properties of the inadmissible construction
are required there. The dimension is chosen after $N_{\rm qry},L,\alpha$, but before the choice of
the algorithm. Therefore, one and the same dimension yields the required bound
for every algorithm and every value of the random variables it uses.

\begin{proof}[Proof of Lemma~\ref{lem:lower-resample}]
Fix a position $j$ and all values of $c$ except $c(j)$. Until $j$ is
revealed, all corresponding runs have the same history: this is verified by induction on time,
since the answers before the revelation do not depend on $c(j)$. Therefore, at any
step before the revelation, the query point or the returned point is the same for all $d$ values
of $c(j)$. If its norm is at most $L$, then
\[
 \alpha^2\#\{v:|x_v|\ge\alpha\}
 \le\sum_{v:|x_v|\ge\alpha}|x_v|^2
 \le\|x\|^2\le L^2.
\]
The fraction of values $c(j)$ for which $|x_{c(j)}|\ge\alpha$ is at most $L^2/(d\alpha^2)$.
It remains to average over the values at the other positions and to sum over the $s$ positions and
the $N_{\rm qry}+1$ points. The probability that $c$ takes equal values at two distinct positions
is at most $\binom{s}{2}/d\le s^2/d$.
\end{proof}

\begin{lemma}[lower bound via the number of revealed links]\label{lem:lower-revelations}
Suppose that for a finite distribution of admissible inputs an idealized process
of revealing $s$ links is given. For every deterministic algorithm
and every deadline $N_{\rm qry}$, denote by $Q$ the number of queries
in this process before stopping or before the given deadline, and by $P$
the number of revealed links. Assume that
$\E Q\ge c_\tau\E P$ for some $c_\tau>0$.

Assume, moreover, that outside the exceptional set the real and idealized histories coincide,
and that the optimality gap of the point returned before this deadline is at least $A-BP$, where
\[
 A\ge7\varepsilon/4,\qquad B>0,\qquad sB\ge A.
\]
If the fraction of exceptions is bounded by~\eqref{eq:lower-pbad}
with fixed $L,\alpha>0$ independent of the algorithm and of the dimension,
then in some finite dimension
\[
 \mathcal T^\ast\ge\frac{c_\tau A}{4B}.
\]
\end{lemma}

\begin{proof}
Set $N_{\rm qry}=\lceil16Ac_\tau/B\rceil$,
$q=BN_{\rm qry}/(Ac_\tau)\ge16$ and
\[
 Z=\frac{BQ}{Ac_\tau}+1-\frac BA P.
\]
By assumption, $\E Z\ge1$, and $Q\le N_{\rm qry}$ implies $Z\le q+1$.
By~\eqref{eq:lower-pbad}, the choice
\[
 d\ge16q\bigl(s^2+(N_{\rm qry}+1)sL^2/\alpha^2\bigr)
\]
ensures $q\Pr(\mathrm{bad})\le1/16$ simultaneously for all algorithms.

Let $\tau$ be the actual number of queries and $\mathcal E$ the optimality gap.
Outside the exceptions, for $\tau<N_{\rm qry}$ we have $Q\le\tau$
and $\mathcal E\ge\max\{0,A-BP\}$; for $\tau\ge N_{\rm qry}$
we use $Z\le q+1$. In all cases,
\[
 \frac{B\tau}{Ac_\tau}+\frac{\mathcal E}{A}
 \ge\frac{q}{q+1}Z-q\mathbf1_{\mathrm{bad}}.
\]
Indeed, for $Z<0$ the nonnegativity of the left-hand side suffices,
and on the exceptions the right-hand side is nonpositive.
Averaging over the inputs and over the random variables of a valid algorithm,
as in Lemma~\ref{lem:lower-finite}, we obtain
\[
 \frac{B\mathcal T(\mathsf M)}{Ac_\tau}
 \ge\frac{q}{q+1}-\frac1{16}-\frac\varepsilon A
 \ge\frac{16}{17}-\frac1{16}-\frac47>\frac14.
\]
It remains to take the infimum over valid algorithms.
\end{proof}

\subsection{Construction of the chain}

Define the symmetric $C^1$ transition function
\[
 h(t)=\begin{cases}
 0,&t\le0,\\
 3t^2-2t^3,&0\le t\le1,\\
 1,&t\ge1.
 \end{cases}
\]
It belongs to $C^1$, is nondecreasing, and satisfies $h(1-t)=1-h(t)$.
On $[0,1]$ we have $h'(t)=6t(1-t)$; hence $0\le h'\le M=3/2$.
Set
\begin{equation}\label{eq:lower-profile}
 \begin{aligned}
 &p(t)=h(4|t|-1)\bigl(1-h(|t|-2)\bigr),\\
 &\widehat\psi'(t)=\int_0^t p(u)\,du,\\
 &\widehat\psi(t)=\int_0^t\widehat\psi'(u)\,du,\\
 &G_H=\widehat\psi'(1),\\
 &c_H=G_H-\widehat\psi(1).
 \end{aligned}
\end{equation}
On $[1/4,1/2]$, the substitution $u=4t-1$ and integration of the polynomial $h$
give
\begin{equation}\label{eq:lower-profile-numbers}
 G_H=\frac58,\qquad \widehat\psi(1)=\frac{63}{320},\qquad
 c_H=\frac{137}{320}.
\end{equation}
The function $p$ is even, $0\le p\le1$, vanishes for $|t|\le1/4$
and $|t|\ge3$, and equals one on $[1/2,2]$.
Consequently, the function $\widehat\psi$ belongs to $C^2$, is convex and nonnegative, vanishes
on $[-1/4,1/4]$, and satisfies the bounds
\begin{equation}
\label{eq:lower-profile-bounds}
 \begin{gathered}
 |\widehat\psi'|\le3,\\
 0\le\widehat\psi''\le1,\\
 \widehat\psi(v)\ge\tfrac32(v-2).
 \end{gathered}
\end{equation}
For $v\ge2$, the last bound follows from the tangent at the point $2$:
\[
 \widehat\psi'(2)\ge\int_{1/2}^2p=3/2.
\]
For arbitrary $c_{\psi},\sigma>0$, consider the family of one-dimensional functions
\[
 \psi(u)=c_{\psi}\widehat\psi(u/\sigma).
\]
The parameters $c_{\psi}$ and $\sigma$ will be chosen separately for each family.

For $w>0$, to smooth the maximum exactly, we define
\begin{gather*}
 S_w'(u)=2h\bigl((u+w)/(2w)\bigr)-1,\\
 S_w(u)=w-\int_0^w S_w'(v)\,dv+\int_0^u S_w'(v)\,dv,\\
 \operatorname{smax}_w(a,b)=(a+b)/2+S_w((a-b)/2).
\end{gather*}
The symmetry of $h$ implies that $S_w$ is even; moreover, $S_w(w)=w$.
From $|S_w'|\le1$ and $0\le S_w''\le M/w$ it follows that
$S_w\ge|\cdot|$, with equality for $|u|\ge w$.
Therefore, the function $\operatorname{smax}_w$ belongs to $C^2$, is convex, nondecreasing in each
argument, not smaller than the maximum of its arguments, and coincides with it for
$|a-b|\ge2w$. For convex $C^2$ functions $a(x),b(x)$,
\begin{equation}\label{eq:lower-smax-hessian}
 \begin{aligned}
 &\nabla\operatorname{smax}_w(a,b)=\frac{\nabla a+\nabla b}{2}
 +\frac12S_w'((a-b)/2)(\nabla a-\nabla b),\\
 &\|\nabla^2\operatorname{smax}_w(a,b)\|\le\|\nabla^2a\|+\|\nabla^2b\|
 +\frac{M}{4w}\|\nabla a-\nabla b\|^2.
 \end{aligned}
\end{equation}
These formulas are obtained by differentiating twice; the additional term
in the Hessian equals
$\frac14S_w''((a-b)/2)(\nabla a-\nabla b)(\nabla a-\nabla b)^\top$.

For an arbitrary integer $T\ge1$, set $y_j=x_{c(j)}$ for $j<T$
and $y_T=0$. Define the chain
\begin{equation}
\label{eq:lower-chain}
 C(x)=-\frac{c_{\psi}G_H}{\sigma}y_0+\sum_{j=0}^{T-1}\psi(y_j-y_{j+1}).
\end{equation}
The convexity of $\psi$ and the tangent at $\sigma$ give, for $0\le p\le T$,
\begin{equation}\label{eq:lower-fenchel}
 C(x) \by{\ge}{(a)}
 -\frac{c_{\psi}G_H}{\sigma}y_0+\sum_{j=0}^{p-1}
 \Bigl(\frac{c_{\psi}G_H}{\sigma}(y_j-y_{j+1})-c_{\psi}c_H\Bigr) \by{=}{(b)}-\frac{c_{\psi}G_H}{\sigma}y_p-pc_{\psi}c_H.
\end{equation}
Here (a) uses the tangent to the function $\psi$ at the point $\sigma$, and (b) uses the cancellation of adjacent terms.

For $p=T$ we obtain $C\ge-Tc_{\psi}c_H$, with equality
at $y_j=(T-j)\sigma$. For injective $c$, the operator
$y\mapsto(y_j-y_{j+1})_{j<T}$ has norm at most $2$:
\[
 \sum_{j<T}(y_j-y_{j+1})^2
 \le2\sum_{j<T}(y_j^2+y_{j+1}^2)
 \le4\sum_{j<T}y_j^2.
\]
Hence, for the full chain and for any sum of its links,
\begin{equation}
\label{eq:lower-chain-bounds}
 \begin{gathered}
 \|\nabla C\|\le \frac{c_{\psi}}{\sigma}(G_H+6\sqrt T),\\
 \|\nabla^2C\|\le\frac{4c_{\psi}}{\sigma^2}.
 \end{gathered}
\end{equation}
If $|y_j|<\alpha\le\sigma/8$ for all $j\ge p$, then the links
with $j\ge p$ vanish in some neighborhood of $x$.
Therefore, the value and the gradient of the chain at the point $x$ coincide with the answer for the function
\[
 C^{[p]}(x)=-\frac{c_{\psi}G_H}{\sigma}y_0+\sum_{0\le j<\min(p,T)}\psi(y_j-y_{j+1}).
\]
The answer for $C^{[p]}$ depends only on the coordinates with indices at most $p$;
before the query some of them may still be hidden, but after this query
at most one new coordinate is revealed.

To handle distant queries, we introduce a radial extension.
It does not change the function in a ball containing its minimizer,
and it hides the chain coordinates outside a larger ball.

\begin{lemma}[radial extension]\label{lem:lower-radial}
Let $r,G,L_{\rm aux}>0$, and let a convex $C^2$ function $C$
attain its minimum $C_\ast$ in $B(0,r)$ and satisfy
\[
 \|\nabla C\|\le G,\qquad
 \|\nabla^2C\|\le L_{\rm aux}/3,\qquad
 4G/r\le L_{\rm aux}/3.
\]
Set
\begin{equation}\label{eq:lower-mask-S}
 \begin{aligned}
 w&=\frac{27MG^2}{4L_{\rm aux}},\\
 \rho(x)&=2G\sqrt{\|x\|^2+r^2}-4Gr+C_\ast-2w,\\
 H(x)&=\operatorname{smax}_w(C(x),\rho(x)).
 \end{aligned}
\end{equation}
Then $H$ is convex, belongs to $C^2$, coincides with $C$ on $B(0,r)$,
has the same minimum, and satisfies
\[
 \|\nabla H\|\le2G,\qquad \|\nabla^2H\|\le L_{\rm aux}.
\]
For $\|x\|\ge4r+(C(0)-C_\ast+4w)/G$ we have $H(x)=\rho(x)\ge C(0)$.
\end{lemma}

\begin{proof}
On $B(0,r)$ we have $\rho\le C_\ast-2w\le C-2w$; hence $H=C$.
Everywhere $H\ge C\ge C_\ast$, so the minimum is preserved.
For the points $x$ in the last statement, we use $C(x)\le C(0)+G\|x\|$:
\[
 \rho(x)\ge2G\|x\|-4Gr+C_\ast-2w
 \ge C(0)+G\|x\|+2w\ge C(x)+2w.
\]
Hence, there $H=\rho\ge C(0)$.

We have $\|\nabla\rho\|\le2G$ and $\|\nabla^2\rho\|\le2G/r$.
The gradient of $H$ is a convex combination of the gradients of $C$ and $\rho$.
By~\eqref{eq:lower-smax-hessian},
\[
 \|\nabla^2H\|\le\frac{L_{\rm aux}}3+\frac{2G}{r}
 +\frac{9MG^2}{4w}\le L_{\rm aux}.
\]
\end{proof}

For the chain~\eqref{eq:lower-chain} we take
$G=c_\psi(G_H+6\sqrt T)/\sigma$ and $C_\ast=-Tc_\psi c_H$.
The conditions of Lemma~\ref{lem:lower-radial} take the form
\begin{equation}\label{eq:lower-mask-conditions}
 \frac{4c_\psi}{\sigma^2}\le\frac{L_{\rm aux}}3,
 \qquad \frac{4c_\psi(G_H+6\sqrt T)}{\sigma r}\le\frac{L_{\rm aux}}3.
\end{equation}
The truncated chains $C^{[p]}$ satisfy the same gradient bound
and vanish at the origin. We extend them using
the same function $\rho$: outside the indicated ball, all answers coincide with $\rho$.
If the full and truncated chains coincide in a neighborhood of a point,
applying the same smoothed maximum preserves this property.
The idea of hiding information from distant queries by means of a radial
extension is used in lower bounds on oracle complexity~\citep{carmon2020lower}.

\subsection{Searching for a hidden pair: the \texorpdfstring{$n$}{n} term}

Choose distinct indices $a,b$ and a sign $s\in\{-1,1\}$
uniformly among all ordered triples. In the one-dimensional space, define
\[
 \begin{gathered}
 f_a(u)=\frac{n\min\left\{\delta,\frac\Delta{\cosh(L_1R)-1}\right\}}{2}e^{L_1(u-sR)},\\
 f_b(u)=\frac{n\min\left\{\delta,\frac\Delta{\cosh(L_1R)-1}\right\}}{2}e^{-L_1(u-sR)},\\
 f_i=0\ (i\ne a,b).
 \end{gathered}
\]

\begin{lemma}\label{lem:lower-pair}
Under Assumption~\ref{asm:accuracy},
$\mathcal T^\ast\ge(n+1)/9$.
\end{lemma}

\begin{proof}[Proof of Lemma~\ref{lem:lower-pair}]
\leavevmode
\begin{itemize}
\item The components of the hidden pair satisfy
Assumptions~\ref{asm:components} and~\ref{asm:smoothness}: they are convex,
smooth, bounded below, and $|f_i''|=L_1|f_i'|$. Their average equals
\[
 f(u)=\min\left\{\delta,\frac\Delta{\cosh(L_1R)-1}\right\}\cosh(L_1(u-sR)),
\]
and the minimum is attained at $sR$, so Assumption~\ref{asm:minimizer} holds.
The initial optimality gap equals $\min\{\Delta,(\cosh(L_1R)-1)\delta\}\le\Delta$
(Assumption~\ref{asm:initial-gap}), and Assumption~\ref{asm:heterogeneity}
follows from
\[
 \begin{gathered}
 f^\ast=\min\left\{\delta,\frac\Delta{\cosh(L_1R)-1}\right\},\\
 f^\ast-n^{-1}\sum_i\inf f_i=\min\left\{\delta,\frac\Delta{\cosh(L_1R)-1}\right\}\le\delta.
 \end{gathered}
\]

\item Fix an algorithm and trace its run on the zero input.
Until the first query to $a$ or $b$, all answers coincide with the answers on the zero input.
If this pair is not touched before stopping, the returned point $u$ does not depend on $s$, and
\begin{align*}
 \frac12\sum_{s=\pm1}(f(u)-f^\ast)&=\min\left\{\delta,\frac\Delta{\cosh(L_1R)-1}\right\}\bigl(\cosh(L_1u)\cosh(L_1R)-1\bigr)
 \\&\by{\ge}{(a)}\min\{\Delta,(\cosh(L_1R)-1)\delta\}.
\end{align*}
Here (a) uses the inequality $\cosh(L_1u)\ge1$.

\item Suppose that the run on the zero input examines $r\le n-1$ distinct
indices (after $n-1$ examined indices, every pair has already been touched).
Repeated queries only increase the number of queries.
The probability of touching the pair equals
\[
 1-\frac{(n-r)(n-r-1)}{n(n-1)}.
\]
The expected number of queries before the first hit or stopping is at least
\begin{align*}
 \sum_{j=0}^{r-1}\frac{(n-j)(n-j-1)}{n(n-1)}
 &=\frac{n(n^2-1)-(n-r)((n-r)^2-1)}{3n(n-1)},\\
 \sum_{j=0}^{r-1}\frac{(n-j)(n-j-1)}{n(n-1)}
 -\frac{n+1}{3}\left(1-\frac{(n-r)(n-r-1)}{n(n-1)}\right)
 &=\frac{(n-r)(n-r-1)r}{3n(n-1)}\ge0.
\end{align*}
The first equality uses the identity
$\sum_{u=1}^v u(u-1)=v(v^2-1)/3$.

If on the zero input the algorithm makes infinitely many queries but accesses fewer than $n-1$ distinct
indices, then on some input of the family the number of queries is infinite,
and the required averaged bound also holds. Denote the number of queries on the chosen
input by $\tau$, and the optimality gap of the output by $\mathcal E$. Then
\[
 \frac{\mathbb E_{a,b,s}\tau}{(n+1)/3}
 +\frac{\mathbb E_{a,b,s}\mathcal E}{\min\{\Delta,(\cosh(L_1R)-1)\delta\}}\ge1.
\]
\item By Lemma~\ref{lem:lower-finite},
\[
 \mathcal T^\ast\ge\frac{n+1}{3}\left(1-
 \frac\varepsilon{\min\{\Delta,(\cosh(L_1R)-1)\delta\}}\right).
\]
From $\cosh u-1\ge u^2/2$ and Assumption~\ref{asm:accuracy}
it follows that $(\cosh(L_1R)-1)\delta\ge3\varepsilon/2$.
By Assumption~\ref{asm:initial-gap}, $\Delta\ge2\varepsilon$;
hence the expression in parentheses is at least $1/3$.
Consequently, $\mathcal T^\ast\ge(n+1)/9$.
\end{itemize}
\end{proof}

\subsection{Family S: independent smooth chains}

Under the conditions of the next lemma, we define $n$ chains~\eqref{eq:lower-chain}
in disjoint coordinates $c(i,j)$, $1\le i\le n$, $0\le j<T$, where
\[
 T=\left\lfloor\sqrt{\frac{c_HL_0R^2}{24n\varepsilon}}\right\rfloor,
 \qquad c_\psi=\frac{2\varepsilon}{Tc_H},
 \qquad \sigma=\sqrt{\frac{24\varepsilon}{Tc_HL_0}}.
\]
Set $f_i=H_i$, where $H_i$ is the extension of the chain $C_i$ by Lemma~\ref{lem:lower-radial}
with radius $R$ and Hessian bound $L_0$.
The common minimizer has coordinates
$x^\ast_{c(i,j)}=(T-j)\sigma$.

After $q_i$ queries to component $i$, we regard
$p_i=\min\{q_i,T\}$ positions of its chain as revealed. The next idealized answer
is given by the extension of $C_i^{[p_i]}$ and reveals at most one position.
For the hidden coordinates we use the threshold
\[
 \alpha=\sigma\min\{1/8,c_H/(8G_H)\}.
\]
The dimension is chosen by Lemma~\ref{lem:lower-revelations} with $s=nT$.

\begin{lemma}\label{lem:lower-S}
If
\begin{equation}\label{eq:lower-S-regime}
 \frac{L_0R^2}{\varepsilon}
 \ge\frac{24(G_H+6)^2}{c_H}\,n,
\end{equation}
then
\begin{equation}\label{eq:lower-S}
 \mathcal T^\ast\gtrsim\sqrt{\frac{nL_0R^2}{\varepsilon}}.
\end{equation}
\end{lemma}

\begin{proof}
Condition~\eqref{eq:lower-S-regime} ensures $T\ge1$ and leaves
a sufficient margin of smoothness for the radial extension.
We verify Assumption~\ref{asm:minimizer}:
\[
 \|x^\ast\|^2=n\sigma^2\sum_{j=1}^Tj^2
 \le\frac{24n\varepsilon}{c_HL_0}T^2\le R^2.
\]
By~\eqref{eq:lower-chain-bounds} we have
\[
 \|\nabla^2C_i\|\le\frac{4c_\psi}{\sigma^2}=\frac{L_0}{3},\qquad
 G^2=\frac{L_0\varepsilon}{6c_H}
       \left(\frac{G_H}{\sqrt T}+6\right)^2
 \le\frac{L_0\varepsilon(G_H+6)^2}{6c_H}
 \le\left(\frac{L_0R}{12}\right)^2.
\]
The last inequality follows from~\eqref{eq:lower-S-regime}.
Hence $4G/R\le L_0/3$, and Lemma~\ref{lem:lower-radial}
yields convex, bounded below, $L_0$-smooth components,
that is, Assumptions~\ref{asm:components} and~\ref{asm:smoothness} hold.
We have $f(0)=0$ and $f^\ast=-2\varepsilon$;
the initial optimality gap equals $2\varepsilon\le\Delta$, and the heterogeneity is zero.
Thus, Assumptions~\ref{asm:initial-gap} and~\ref{asm:heterogeneity} hold.

Outside the exceptional set, the real and idealized histories coincide:
for interior queries this follows from the presence of a flat region of $\psi$,
and for exterior ones from Lemma~\ref{lem:lower-radial}.
The same holds for the returned point $\hat x$.
Inside the extension ball, by~\eqref{eq:lower-fenchel},
\[
 f_i(\hat x)\ge-\frac{2\varepsilon}{T}\left(p_i+\frac18\right).
\]
For $p_i=T$ it suffices that $f_i\ge-2\varepsilon$.
Therefore, with $P=\sum_i p_i$, the optimality gap satisfies
\[
 \mathcal E\ge2\varepsilon-\frac{\varepsilon}{4T}
 -\frac{2\varepsilon}{nT}P
 \ge\frac{7\varepsilon}{4}-\frac{2\varepsilon}{nT}P.
\]
Outside the ball it is at least $2\varepsilon$.
Each query reveals at most one link.
We apply Lemma~\ref{lem:lower-revelations} with
$c_\tau=1$, $A=7\varepsilon/4$, $B=2\varepsilon/(nT)$.
The condition $sB\ge A$ holds, and
\[
 \mathcal T^\ast\ge\frac{7nT}{32}
 \ge\frac7{64}\sqrt{\frac{c_HnL_0R^2}{24\varepsilon}}.
\]
Here we used the inequality $\lfloor u\rfloor\ge u/2$ for $u\ge1$.
\end{proof}

\subsection{Family X: exponential transformation of chains}

\begin{lemma}[exponential transformation]\label{lem:lower-lift}
Let a function $m\ge0$ be convex, belong to $C^2$, and satisfy
$\|\nabla m\|\le L_1/7$, $\|\nabla^2m\|\le4L_1^2/49$.
Then, for $\delta'>0$, both functions
\[
 v(x)=\delta'\bigl(e^{\pm2L_1x_0/7+m(x)}-1\bigr)
\]
are convex, belong to $C^2$, have infimum $-\delta'$,
and satisfy $\|\nabla^2v\|\le L_1\|\nabla v\|$.
\end{lemma}

\begin{proof}
Set $u=\pm2L_1x_0/7+m(x)$. The function $m$ may depend on $x_0$.
We have
\begin{equation}\label{eq:lower-lift}
 \begin{aligned}
 L_1/7&\le\|\nabla u\|\le3L_1/7,\\
 \|\nabla^2u\|+\|\nabla u\|^2
 &\le\left(\frac{4L_1}{7}+\frac{3L_1}{7}\right)\|\nabla u\|
 =L_1\|\nabla u\|,\\
 \|\nabla^2v\|
 &\le\delta'e^u\bigl(\|\nabla^2u\|+\|\nabla u\|^2\bigr)
 \le L_1\|\nabla v\|.
 \end{aligned}
\end{equation}
The convexity of $v$ follows from the convexity of $u$.
Along $x=\mp te_0$ we have $u(x)\le m(0)-L_1t/7\to-\infty$;
hence $\inf v=-\delta'$.
\end{proof}

Thus, $v$ belongs to the required class for every $L_0\ge0$.
Verifying membership in the class via a Hessian bound corresponds to the characterization
from the first part~\citep[Lemma~2.2]{rudomanov2026everest}.

Set
\begin{equation}\label{eq:lower-X-parameters}
 \begin{gathered}
 \eta=\log\left(1+\frac{2\varepsilon}{\delta}\right),\qquad
 \rho=\left(1+\frac{\delta}{2\varepsilon}\right)\eta,\\
 D_\delta=\max\left\{\sqrt{\frac{12\eta}{c_H}},
                   \frac{4(G_H+6)\eta}{c_H}\right\},
 \qquad U=\frac{2L_1R}{7D_\delta}.
 \end{gathered}
\end{equation}
Here $D_\delta\sim\max\{\sqrt\eta,\eta\}$.
First consider arbitrary integers $g,T\ge1$ with $2g\le n$.
We construct $g$ chains~\eqref{eq:lower-chain} in disjoint coordinates
$c(b,j)$, $1\le b\le g$, $0\le j<T$, with parameters
\[
 c_\psi=\frac{\eta}{Tc_H},\qquad
 \sigma=\frac7{2L_1}\max\left\{
 \sqrt{\frac{12\eta}{Tc_H}},\frac{4(G_H+6\sqrt T)\eta}{Tc_H}\right\}.
\]
Denote by $\sigma_1$ the value of $\sigma$ for $T=1$.
Then $\sigma\le\sigma_1/\sqrt T$ and $U=R/\sigma_1$.
We extend each chain by Lemma~\ref{lem:lower-radial}
with radius $r_X=\max\{R,21/(2L_1)\}$ and Hessian bound $4L_1^2/49$.
The resulting functions are denoted by $H_b$.
Their common minimum is attained at the point
$x^\ast_{c(b,j)}=(T-j)\sigma$.

We define $2g$ active components:
\[
 \begin{aligned}
 m_b(x)&=\eta+H_b(x),\\
 f_{b,\pm}(x)&=\frac{n\delta}{2g}
 \bigl(e^{\pm2L_1x_0/7+m_b(x)}-1\bigr).
 \end{aligned}
\]
The remaining components are zero. The hidden coordinates are different from $0$;
hence $x^\ast_0=0$.

A query to either component of a pair reveals at most one link of its chain.
For the hidden coordinates we use the threshold
\[
 \alpha=\sigma\min\{1/8,c_H/(8G_H\rho)\}.
\]
The dimension is chosen by Lemma~\ref{lem:lower-revelations} with $s=gT$.
For $U\ge1$ we take
\[
 g=\min\{\lfloor n/2\rfloor,\lfloor U^2\rfloor\},
 \qquad T=\lfloor U/\sqrt g\rfloor.
\]
Both numbers are positive, and $2g\le n$.

\begin{lemma}\label{lem:lower-X}
If $U\ge1$, then
\begin{equation}\label{eq:lower-X}
 \mathcal T^\ast\gtrsim
 \frac{U\min\{\sqrt n,U\}}{\rho}.
\end{equation}
\end{lemma}

\begin{proof}
From the choice of $\sigma$ and~\eqref{eq:lower-chain-bounds} we obtain
\[
 \|\nabla C_b\|\le L_1/14,\qquad
 \|\nabla^2C_b\|\le4L_1^2/147.
\]
Together with $r_X\ge21/(2L_1)$, this means that the conditions
of Lemma~\ref{lem:lower-radial} hold. Moreover,
\[
 \|x^\ast\|^2=g\sigma^2\sum_{j=1}^Tj^2
 \le g\sigma_1^2T^2\le R^2.
\]
Consequently, $m_b\ge0$,
$\|\nabla m_b\|\le L_1/7$ and $\|\nabla^2m_b\|\le4L_1^2/49$.
By Lemma~\ref{lem:lower-lift}, Assumptions~\ref{asm:components}
and~\ref{asm:smoothness} hold, and
\[
 f(x)=\frac\delta g\sum_{b=1}^g
 \bigl(\cosh(2L_1x_0/7)e^{m_b(x)}-1\bigr)\ge0.
\]
At the point $x^\ast$ all $m_b=0$; hence $f(x^\ast)=f^\ast=0$,
and Assumption~\ref{asm:minimizer} holds.
At zero $m_b=\eta$, whence $f(0)=2\varepsilon\le\Delta$
(Assumption~\ref{asm:initial-gap}).
Since $\inf f_{b,\pm}=-n\delta/(2g)$, the heterogeneity equals $\delta$
(Assumption~\ref{asm:heterogeneity}).

The coincidence of the real and idealized histories is proved in the same way as for family S:
the exponential transformation preserves local equality of functions.
Let $p_b\le T$ be the number of revealed links of pair $b$, and let $P=\sum_b p_b$.
Inside the extension ball and outside the exceptional set,
bound~\eqref{eq:lower-fenchel} gives
\[
 m_b(\hat x)\ge\eta\left(1-\frac{p_b}{T}\right)
 -\frac{\eta}{8T\rho}.
\]
For $p_b=T$ it suffices that $m_b\ge0$.
The tangent to the exponential at $\eta$ and the equality
$e^\eta\eta\delta=2\varepsilon\rho$ give
\[
 \frac{\mathcal E}{2\varepsilon}
 \ge1-\frac{\rho}{gT}P-\frac1{8T}.
\]
Outside the ball, the optimality gap is at least $2\varepsilon$.
We apply Lemma~\ref{lem:lower-revelations} with
$c_\tau=1$, $A=7\varepsilon/4$, $B=2\varepsilon\rho/(gT)$.
The inequality $(1+1/z)\log(1+z)\ge1$ for $z>0$ gives
$\rho\ge1$; hence $sB\ge A$.
We obtain $\mathcal T^\ast\ge7gT/(32\rho)$.
From $n\ge2$, $U\ge1$ and the definition of $g,T$ it follows that
\[
 g\ge\frac13\min\{n,U^2\},\qquad
 T\ge\frac{U}{2\sqrt g},\qquad
 gT\ge\frac{U\min\{\sqrt n,U\}}{2\sqrt3}.
\]
These inequalities imply~\eqref{eq:lower-X}.
\end{proof}

\subsection{Number of queries needed to find a link}

\begin{lemma}[number of queries needed to find a link]\label{lem:lower-search}
Let $T$ links be placed by an injection
$\pi:\{0,\ldots,T-1\}\to\{1,\ldots,m\}$, $T\le m$.
At level $j$, the idealized answer contains only the links with indices
at most $j$, and the level increases by one only upon a query
to the pair $\pi(j)$. Let $Q_j$ be the number of queries at level
$j$ before stopping or before time $N_{\rm qry}$, and let $R_j$ be the indicator
that this level has been passed by time $N_{\rm qry}$. Then, for every
fixed coordinate map,
\begin{equation}\label{eq:lower-search}
 \mathbb E_\pi Q_j\ge\frac{m-T+2}{2}\mathbb E_\pi R_j.
\end{equation}
\end{lemma}

\begin{proof}[Proof of Lemma~\ref{lem:lower-search}]
\leavevmode
\begin{itemize}
\item Fix the placement of all links except $j$.
For the placement of link $j$, there remain $m-T+1$ free pairs.
Before level $j$, the history does not depend on the placement of this link.
At level $j$ itself, we compare it with the virtual history in which
link $j$ is placed in no pair. The runs coincide until the first
query to the chosen free pair. Suppose that in the virtual history
$r$ distinct free pairs are queried before stopping or before time $N_{\rm qry}$.
\item Exactly $r$ placement options lead to passing the level;
for the $k$-th of them at least $k$ queries are made at this level,
and for each of the remaining $(m-T+1)-r$ at least $r$.
Therefore,
\begin{align*}
 \sum_{v=1}^{(m-T+1)}Q_j(v)&\by{\ge}{(a)}
 \sum_{k=1}^r k+((m-T+1)-r)r =r\left((m-T+1)+\frac{1-r}{2}\right)
 \\&\by{\ge}{(b)}
 \frac{(m-T+1)+1}{2}\sum_{v=1}^{(m-T+1)}R_j(v).
\end{align*}
Here (a) takes into account the order of the first visits to the free pairs, and (b) uses the relations $r\le (m-T+1),\;\sum_vR_j(v)=r$.

\item If the level is not reached at all, both right-hand quantities are zero.
Repeated queries and queries to occupied pairs only increase the left-hand side.
Summing over all fixed placements of the other links,
we obtain~\eqref{eq:lower-search}. This is precisely averaging over all injections:
each of them enters exactly one layer with weight $1/(m-T+1)$.
\end{itemize}
\end{proof}

\subsection{Family Q: searching for the pair with the next link}

The exponent $2/3$ already appeared in the first part in the lower bound for a single
function with bounded initial optimality gap
\citep[Theorem~B.3]{rudomanov2026everest}.
The family Q used here is defined by exponentials of auxiliary functions
and by a random placement of links over pairs of components.

First we construct the family for arbitrary integers $m,T\ge1$
such that $2m\le n$, $T\le m/2$ and
\[
 \begin{gathered}
 2L_1R/7\ge3,\\
 (2L_1R/7)^{2/3}\le T\le2(2L_1R/7)^{2/3},\\
 24\log\!\left(1+\frac{2\varepsilon m}{\delta}\right)\le\widehat\psi(1)(2L_1R/7)^{2/3}.
 \end{gathered}
\]
The specific numbers of pairs and links will be set below.

We define the auxiliary functions
\[
 \begin{aligned}
 \Psi_0(u)&=\frac{\log\!\left(1+\frac{2\varepsilon m}{\delta}\right)}{\widehat\psi(1)}
             \widehat\psi(u\sqrt T/R),\\
 \Psi(u)&=\frac1{48}\widehat\psi(4L_1u/7).
 \end{aligned}
\]
The links have the form $\Psi_0(R/\sqrt T-x_{c(0)})$ and
$\Psi(x_{c(j-1)}-x_{c(j)})$, $1\le j<T$. The coordinates $c(j)$ are pairwise distinct and different from $0$.
The links are placed over $m$ pairs by a uniformly chosen injection $\pi$, so that each pair contains
at most one link. Denote by $\Psi_b$ the sum of the links of pair $b$, and set
\[
 \begin{aligned}
 \rho(x)&=\frac{L_1}7\sqrt{\|x\|^2+R^2}-\frac{2L_1R}7-\frac{27M}{32},\\
 m_b(x)&=\operatorname{smax}_{27M/64}(\Psi_b(x),\rho(x)).
 \end{aligned}
\]

Define the point
\[
 x^\ast=\frac R{\sqrt T}\sum_{j<T}e_{c(j)}.
\]

Define the components
\[
 f_{b,\pm}(x)=\frac{n\delta}{2m}
 \left(e^{\pm2L_1x_0/7+m_b(x)}-1\right).
\]
The remaining components are zero.
Define the functions
\[
 \begin{gathered}
 \phi_0(u)=e^{\Psi_0(u)}-1,\\
 \phi(u)=e^{\Psi(u)}-1.
 \end{gathered}
\]

At level $p$, the idealized oracle keeps only the links
$j\le p$ and increases the level after a query to the pair $\pi(p)$.
Before the query, the coordinates with indices $j<p$ are revealed, and after it at most one new coordinate is revealed.
For the hidden coordinates we use the threshold
\[
 \min\left\{\frac7{32L_1},
 \frac{\varepsilon m}{4\phi_0'(R/\sqrt T)\delta}\right\}.
\]

For the next lemma, set
\[
 t=(2L_1R/7)^{2/3},\qquad m=\lfloor n/2\rfloor,\qquad
 T=\lceil t\rceil,\qquad \eta_m=\log(1+2\varepsilon m/\delta).
\]

\begin{lemma}\label{lem:lower-Q}
If
\begin{equation}\label{eq:lower-Q-regime}
 L_1R\ge\frac{21}{2},\qquad 4t\le m,\qquad
 24\eta_m\le\widehat\psi(1)t,
\end{equation}
then
\begin{equation}\label{eq:lower-Q}
 \mathcal T^\ast\gtrsim
 \frac{n(L_1R)^{2/3}}{(1+\eta_m)^2}.
\end{equation}
\end{lemma}

\begin{proof}[Proof of Lemma~\ref{lem:lower-Q}]
The first condition in~\eqref{eq:lower-Q-regime} allows us to extend the links
beyond the ball $B(0,R)$. The second gives $t\le T\le2t\le m/2$:
when searching for the next link, at least half of the pairs remain free.
The third bounds the curvature of the first link.
\begin{itemize}
\item For an ordinary link, \eqref{eq:lower-profile-bounds} implies the bounds
$L_1/14$ for the gradient norm and $4L_1^2/147$ for the Hessian norm.
For the first link, the same bounds follow from
\[
 \begin{aligned}
 \frac{3\sqrt T\,\eta_m}{\widehat\psi(1)R}
 &\le\frac{6\sqrt t\,\eta_m}{\widehat\psi(1)R}
 \le\frac{t^{3/2}}{4R}=\frac{L_1}{14},\\
 \frac{T\eta_m}{\widehat\psi(1)R^2}
 &\le\frac{t^2}{12R^2}\le\frac{4L_1^2}{147}.
 \end{aligned}
\]
For $2L_1R/7\ge3$, Lemma~\ref{lem:lower-radial} applies
with $G=L_1/14$, $L_{\rm aux}=4L_1^2/49$, $r=R$, $C_\ast=0$.
Therefore,
\[
 \begin{gathered}
 \|\nabla m_b\|\le L_1/7,\\
 \|\nabla^2m_b\|\le4L_1^2/49.
 \end{gathered}
\]
All potentials are nonnegative. The point $x^\ast$ has norm $R$,
and all links and potentials vanish at it.
At zero, only the potential of the pair containing link $0$ is nonzero; it equals
 $\eta_m$.
As for family X, the minimum of the average function equals
zero, the heterogeneity equals $\delta$, and
\[
 f(0)=2\varepsilon\le\Delta.
\]
Thus, Assumptions~\ref{asm:minimizer}--\ref{asm:heterogeneity} hold.
Assumptions~\ref{asm:components} and~\ref{asm:smoothness}
follow from Lemma~\ref{lem:lower-lift}.

\item We prove a lower bound on the optimality gap. We have
\[
 \begin{aligned}
 \phi_0(R/\sqrt T)&=e^{\eta_m}-1=\frac{2\varepsilon m}{\delta},\\
 \phi_0'(R/\sqrt T)&=
 \frac{G_H\sqrt T}{R\widehat\psi(1)}e^{\eta_m}\eta_m.
 \end{aligned}
\]
Denote $s_0=\phi_0'(R/\sqrt T)>0$,
$\log_+u=\max\{0,\log u\}$ and
\[
 J=1+16\log_+(56s_0/L_1),\qquad K=\frac{7s_0}{2L_1}J.
\]
The quantity $K$ bounds the loss in the lower bound when one link is added:
\begin{equation}\label{eq:lower-conjugate}
 \sup_{u\in\mathbb R}\{s_0u-\phi(u)\}\le K.
\end{equation}
For $u\le7/(2L_1)$ it suffices that $\phi\ge0$.
For $u>7/(2L_1)$, set $v=4L_1u/7-2>0$.
By~\eqref{eq:lower-profile-bounds}, $\Psi(u)\ge v/32$.
For arbitrary $\mu\ge0$, $a>0$, $v\ge0$, the tangent to the exponential at
$\log(\mu/a)$ for $\mu>a$ and $e^{av}\ge1+av$ for $\mu\le a$ give
\[
 \mu v-(e^{av}-1)\le\frac\mu a\max\{0,\log(\mu/a)\}.
\]
Substituting $\mu=7s_0/(4L_1)$ and $a=1/32$, we obtain
\[
 s_0u-\phi(u)
 \le\frac{7s_0}{2L_1}
     +\frac{56s_0}{L_1}\log_+(56s_0/L_1)=K.
\]

If the first $p<T$ coordinates are revealed and the absolute value of the next coordinate
is less than the chosen threshold, then the sum of the differences
$R/\sqrt T-x_{c(0)}$ and $x_{c(j-1)}-x_{c(j)}$, $1\le j\le p$,
equals $R/\sqrt T-x_{c(p)}$.
The tangent to $\phi_0$ and~\eqref{eq:lower-conjugate} give
\begin{equation}\label{eq:lower-Q-gap}
 \begin{aligned}
 &\phi_0(R/\sqrt T-x_{c(0)})+\sum_{j=1}^p\phi(x_{c(j-1)}-x_{c(j)})\\
 &\quad\ge\frac{2\varepsilon m}{\delta}-s_0x_{c(p)}-pK.
 \end{aligned}
\end{equation}
The injectivity of $\pi$ allows us to bound the average function from below by the sum
of these exponential links: an empty pair makes a nonnegative contribution,
and a pair with link $j$ contributes at least $(\delta/m)\phi(x_{c(j-1)}-x_{c(j)})$.

\item The threshold for the hidden coordinates chosen in the construction of Q
does not exceed $7/(32L_1)$; therefore,
the distant links vanish in a neighborhood of the query point.
For
\[
 \|x\|\ge4R+\frac{14}{L_1}\log\!\left(1+\frac{2\varepsilon m}{\delta}\right)+\frac{189M}{8L_1},
\]
by Lemma~\ref{lem:lower-radial}, the full and truncated components are determined only by the radial extension
and the coordinate $x_0$. All potentials are at least $\log\!\left(1+\frac{2\varepsilon m}{\delta}\right)$, and the optimality gap
is at least the initial one. Thus, the histories coincide
for arbitrary distant queries as well.
By~\eqref{eq:lower-Q-gap} and the choice of the threshold, the optimality gap after
passing $p$ levels is at least $A-Bp$, where one can take
\[
 A=7\varepsilon/4,\qquad
 B=\frac{K\delta}{m}=\frac{7s_0\delta}{2mL_1}J.
\]
Indeed, after multiplication by $\delta/m$, the term with the hidden coordinate
in~\eqref{eq:lower-Q-gap} causes a loss of at most $\varepsilon/4$.
To apply Lemma~\ref{lem:lower-revelations}, we verify two relations.
From
\[
 s_0=\frac{G_H}{\widehat\psi(1)}\frac{\sqrt T}{R}e^{\eta_m}\eta_m,
 \qquad \frac{\delta}{m}=\frac{2\varepsilon}{e^{\eta_m}-1}
\]
and $\sqrt T\le2\sqrt t$ it follows that
\begin{equation}\label{eq:lower-Q-ratio}
 \frac AB
 =\frac{\widehat\psi(1)}{4G_H}
   \frac{L_1R}{\sqrt T}\frac{1-e^{-\eta_m}}{\eta_m J}
 \ge\frac{7\widehat\psi(1)}{16G_H}
      \frac{t}{(1+\eta_m)J}.
\end{equation}
Here we used the equality $L_1R=7t^{3/2}/2$ and the inequality
$1-e^{-u}\ge u/(1+u)$ for $u\ge0$.
Moreover, $T^{3/2}\ge2L_1R/7$, $J\ge1$ and
$G_H/\widehat\psi(1)\ge1$ give
\[
 TB=\frac{7\delta}{2mL_1}\frac{G_H}{\widehat\psi(1)}
       \frac{T^{3/2}}R e^{\eta_m}\eta_m J
 \ge\frac\delta m e^{\eta_m}\eta_m
 \ge2\varepsilon>A.
\]
For $p=T$, the bound on the optimality gap follows from $f\ge0$.

\item Summing~\eqref{eq:lower-search} over the levels, we obtain
$\E Q\ge c_\tau\E P$ with $c_\tau=(m-T+2)/2\ge m/4$.
Lemma~\ref{lem:lower-revelations} with $s=T$ gives
$\mathcal T^\ast\ge mA/(16B)$.
It remains to bound $J$. By the third condition in~\eqref{eq:lower-Q-regime},
\[
 \frac{56s_0}{L_1}
 \le\frac{32G_H}{\widehat\psi(1)}\frac{e^{\eta_m}\eta_m}{t}
 \le\frac{4G_H}{3}e^{\eta_m}
 =\frac56e^{\eta_m}.
\]
Consequently, $J\le1+16\eta_m$. Substituting this bound
into~\eqref{eq:lower-Q-ratio} and using $m\ge n/3$ and
$t=(2/7)^{2/3}(L_1R)^{2/3}$, we obtain~\eqref{eq:lower-Q}.
\end{itemize}
\end{proof}

\section{Proof of the lower bound}\label{app:lower-constants}

\begin{proof}[Proof of Theorem~\ref{thm:lower}\proofup{thm:lower}]
It suffices to bound the four terms $B_n,B_S,B_X,B_Q$
from~\eqref{eq:lower-four-terms}: their sum is of order $T_0$.
Set $a=L_1R$; we use the quantity $N_0$ from~\eqref{eq:complexity-parameters}.
The numbers $G_H$, $\widehat\psi(1)$ and $c_H$ are computed
in~\eqref{eq:lower-profile-numbers}.

\subsection{The \texorpdfstring{$L_0$}{L0} term}
By Lemma~\ref{lem:lower-pair},
\begin{equation}\label{eq:lower-n-number}
 B_n=n\le9\mathcal T^\ast.
\end{equation}
If~\eqref{eq:lower-S-regime} holds, then
Lemma~\ref{lem:lower-S} applies. Otherwise,
\[
 B_S<\sqrt{24/c_H}\,(G_H+6)n<50n.
\]
In both cases,
\begin{equation}\label{eq:lower-S-number}
 B_S=\sqrt{\frac{nL_0R^2}{\varepsilon}}\lesssim\mathcal T^\ast.
\end{equation}

\subsection{Terms obtained from family X}
We use the notation $\eta,\rho,D_\delta,U$ from~\eqref{eq:lower-X-parameters}.
By Assumption~\ref{asm:accuracy},
$\varepsilon/\delta\le a^2/3$, whence
\begin{equation}\label{eq:lower-rho-number}
 \eta\le\log(1+2a^2/3)\le\log N_0,\qquad
 1\le\rho\le1+\eta\le1+\log N_0.
\end{equation}
For the second bound, the inequalities
$z/(1+z)\le\log(1+z)\le z$ with $z=2\varepsilon/\delta$ suffice.
Moreover,
\begin{equation}\label{eq:lower-D-bound}
 D_\delta\sqrt{1+\delta/\varepsilon}\lesssim\log N_0.
\end{equation}
Indeed, $\eta\delta/\varepsilon\le2$ by the inequality
$\log(1+z)\le z$, and
$\sqrt{1+\delta/\varepsilon}\le1+\delta/(2\varepsilon)$. Therefore,
\[
 D_\delta\sqrt{1+\delta/\varepsilon}
 \le\max\left\{\sqrt{\frac{12(\eta+2)}{c_H}},
                 \frac{4(G_H+6)}{c_H}\rho\right\}
 \lesssim\log N_0.
\]

We show that
\begin{equation}\label{eq:lower-Y-number}
 \sqrt n\,U\lesssim\log N_0\,\mathcal T^\ast.
\end{equation}
For $U\ge\sqrt n$, this follows from Lemma~\ref{lem:lower-X}
and~\eqref{eq:lower-rho-number}.
For $1\le U<\sqrt n$, the same lemma gives $U^2\lesssim\rho\mathcal T^\ast$.
Multiplying this inequality and~\eqref{eq:lower-n-number} side by side, we obtain
$\sqrt n\,U\lesssim\sqrt\rho\,\mathcal T^\ast\lesssim\log N_0\,\mathcal T^\ast$.
Finally, for $U<1$ it suffices that $\sqrt n\,U\le n$.
From $a=7UD_\delta/2$, \eqref{eq:lower-D-bound}
and~\eqref{eq:lower-Y-number} it follows that
\begin{equation}\label{eq:lower-X-number}
 B_X=\sqrt n\,a\sqrt{1+\delta/\varepsilon}
 \lesssim\log^2N_0\,\mathcal T^\ast.
\end{equation}

\subsection{The term with exponent \texorpdfstring{$2/3$}{2/3}}
Set $t=(2a/7)^{2/3}$, $m=\lfloor n/2\rfloor$ and
$\eta_m=\log(1+2\varepsilon m/\delta)$.
By Assumption~\ref{asm:accuracy},
\[
 \eta_m\le\log(1+na^2/3)\le2\log N_0.
\]
If conditions~\eqref{eq:lower-Q-regime} hold,
Lemma~\ref{lem:lower-Q} gives
\begin{equation}\label{eq:lower-Q-number}
 B_Q=na^{2/3}\lesssim(1+\eta_m)^2\mathcal T^\ast
 \lesssim\log^2N_0\,\mathcal T^\ast.
\end{equation}
If at least one of the conditions is violated, the bound follows from those already proved:
\begin{itemize}
\item If $a<21/2$, then $B_Q<5n$.
\item If $a\ge21/2$ but $m<4t$, then $t>1$ and
$n\le2m+1<9t$. Therefore,
\[
 \frac{B_Q}{\sqrt n\,a}=\frac{\sqrt n}{a^{1/3}}
 <3(2/7)^{1/3}<2,
 \qquad B_Q<2B_X.
\]
\item If $24\eta_m>\widehat\psi(1)t$, then
\[
 B_Q<(7/2)^{2/3}\frac{24}{\widehat\psi(1)}n\eta_m
 \lesssim n\log N_0.
\]
\end{itemize}
Thus, bound~\eqref{eq:lower-Q-number} holds in all cases.
Adding the bounds for the four terms, we obtain
\[
 T_0\sim B_n+B_S+B_X+B_Q
 \lesssim\log^2N_0\,\mathcal T^\ast.
\]
\end{proof}

\Needspace{12\baselineskip}
\section{Variants of YETI}
\label{app:variants}

\subsection{Why the local curvature must be bounded}

The condition $\norm{y-x^\ast}\le R$ by itself does not yield a bound
on the ratio of local smoothness constants that is polynomial in $L_1R$.
Consider the one-dimensional function
\[
 \varphi(t)=\cosh(L_1t).
\]
It is convex, attains its minimum at zero, and satisfies
\[
 \varphi''(t)=L_1^2\cosh(L_1t)
 \le L_1^2+L_1|\varphi'(t)|.
\]
However,
\[
 \frac{\varphi''(R)}{\varphi''(0)}=\cosh(L_1R)
\]
grows exponentially. Therefore, the mere fact that the centers lie in the ball
is not sufficient to drop the curvature test while preserving the same
bound on the number of oracle calls per local call. This example concerns only the solver
used here and does not rule out the existence of other algorithms.

\subsection{Choice of the stochastic gradient}

Let the functions $\varphi_i$ have local smoothness constants $L_i>0$
on the enlarged ball, and let the points $x,w$ lie in the original ball.
For a random index $I$ with probabilities $q_i>0$,
Lemma~\ref{lem:local-bregman} gives
\begin{equation}\label{eq:importance-variance}
 \begin{aligned}
 \mathbb E\left\|
 \frac{\nabla\varphi_I(x)-\nabla\varphi_I(w)}{nq_I}
 \right\|^2&=\frac1{n^2}\sum_i
       \frac{\norm{\nabla\varphi_i(x)-\nabla\varphi_i(w)}^2}{q_i}
 \\&\le\frac2{n^2}\sum_i\frac{L_i}{q_i}
       \bigl(\varphi_i(x)-\varphi_i(w)
                  -\ip{\nabla\varphi_i(w)}{x-w}\bigr).
 \end{aligned}
\end{equation}
If
\[
 q_i=\frac{L_i}{\sum_jL_j},
\]
then the last line equals
\[
 \left(\frac2n\sum_iL_i\right)\frac1n\sum_i
 \bigl(\varphi_i(x)-\varphi_i(w)-\ip{\nabla\varphi_i(w)}{x-w}\bigr).
\]
Under uniform sampling, the standard upper bound contains
$2\max_iL_i$ instead of $2n^{-1}\sum_iL_i$.

The ratio of these quantities can be of order $n$ even within the class under consideration.
In the one-dimensional case with $L_0>0$, take
\[
 \begin{gathered}
 f_1(t)=\frac{L_0}{2}(t-R)^2,\\
 f_i(t)=0\quad(i>1).
 \end{gathered}
\]
Set the known upper bounds $L_{0,i}=L_0$ for all $i$.
Choose the parameters
\[
 \begin{gathered}
 y=0,\\
 L_1R=n,\\
 r=\frac1{2L_1},\\
 \lambda=L_0.
 \end{gathered}
\]
Then the constants from~\eqref{eq:local-smoothness}, that is,
\[
 L_i=e^{3L_1r}\bigl(L_{0,i}+L_1\norm{\nabla f_i(y)}\bigr)+\lambda,
\]
satisfy
\[
 \frac{\max_iL_i}{n^{-1}\sum_iL_i}
 =\frac{1+e^{3/2}(1+n)}{1+2e^{3/2}}.
\]
This ratio grows linearly in $n$. Hence, for this bound of the
inner solver, non-uniform sampling is essential. It does not follow
that every method with uniform sampling is necessarily worse.

The YETI solver uses the probabilities
\[
 q_i=\frac{L_i}{\sum_jL_j}.
\]
If another distribution guarantees only
\[
 \frac{L_i}{nq_i}\le\frac cn\sum_jL_j,
\]
then a factor $c$ appears on the right-hand side of
bound~\eqref{eq:importance-variance}. Its effect on the complexity bound
must be taken into account together with the choice of the epoch parameters.

Storing the anchor gradient is also essential. For $a>0$ and the functions
\[
 \begin{gathered}
 f_1(t)=\frac{(t-a)^2}{2},\\
 f_2(t)=\frac{(t+a)^2}{2}
 \end{gathered}
\]
the minimum of their average is attained at zero, but
\[
 \frac12\sum_{i=1}^2|f_i'(0)|^2=a^2>0.
\]
Therefore, the gradient of a randomly chosen component has nonzero
variance at the solution. In bound~\eqref{eq:importance-variance}, by contrast,
the gradient differences vanish at $x=w$.
Without the correction based on the stored gradients, the variance at the solution
may remain positive, and the above proof of
geometric convergence does not apply.

\subsection{Componentwise secants for generalized linear models}

Let
\[
 f_i(x)=\phi_i(a_i^{\mathsf T}x),
\]
where $\phi_i$ is convex and twice continuously differentiable.
Using the already computed gradients at the points $x_0,x_1$, define
\[
 \begin{gathered}
 s=x_1-x_0,\\
 \Delta g_i=\nabla f_i(x_1)-\nabla f_i(x_0).
 \end{gathered}
\]

\begin{proposition}[recovering the average Hessian]\label{prop:glm-secant}
If $\ip{\Delta g_i}s>0$, then
\begin{equation}\label{eq:glm-secant}
 \frac{\Delta g_i\Delta g_i^{\mathsf T}}{\ip{\Delta g_i}s}
 =\int_0^1\nabla^2f_i(x_0+ts)\,dt.
\end{equation}
If $a_i^{\mathsf T}s=0$, then a zero secant may correspond to a
nonzero matrix on the right-hand side. If, however, the average scalar curvature
is zero and the secant is explicitly defined as the zero matrix (without division
by zero), then both sides of the equality are zero.
\end{proposition}

\begin{proof}
By the Newton--Leibniz formula,
\[
 \begin{aligned}
 \Delta g_i&=(a_i^{\mathsf T}s)a_i
     \int_0^1\phi_i''(a_i^{\mathsf T}(x_0+ts))\,dt,\\
 \ip{\Delta g_i}s&=(a_i^{\mathsf T}s)^2
     \int_0^1\phi_i''(a_i^{\mathsf T}(x_0+ts))\,dt,\\
 \int_0^1\nabla^2f_i(x_0+ts)\,dt
 &=a_i a_i^{\mathsf T}\int_0^1\phi_i''(a_i^{\mathsf T}(x_0+ts))\,dt.
 \end{aligned}
\]
When the denominator is positive, cancelling the scalar factors
gives~\eqref{eq:glm-secant}. When $a_i^{\mathsf T}s=0$,
the equality $\Delta g_i=0$ does not determine the integral of $\phi_i''$.
If, however, this integral is zero, then the average Hessian and the gradient difference are zero.
\end{proof}

Thus, a positive componentwise secant recovers the
average Hessian along the segment. Obtaining the Hessian at the endpoint of the segment
requires an additional assumption.

\begin{proposition}[approximating the Hessian of the Poisson model]\label{prop:poisson-secants}
Let
\[
 f_i(x)=e^{a_i^{\mathsf T}x}-y_i a_i^{\mathsf T}x,
\]
where $y_i\ge0$ and $\norm{a_i}\le L_1$.
These components satisfy the condition of the class with
\[
 L_{0,i}=y_i\norm{a_i}^2.
\]
If $a_i^{\mathsf T}s\ne0$ for every nonzero $a_i$,
then the matrix $M$ obtained by averaging the secants~\eqref{eq:glm-secant}
satisfies
\begin{equation}\label{eq:poisson-secant-comparison}
 e^{-L_1\norm s}\nabla^2f(x_1)
 \preceq M\preceq e^{L_1\norm s}\nabla^2f(x_1).
\end{equation}
Zero $a_i$ contribute zero matrices to the sum.
\end{proposition}

\begin{proof}
Directly from the form of the component we obtain
\[
 \norm{\nabla^2f_i(x)}
 =e^{a_i^{\mathsf T}x}\norm{a_i}^2
 \le y_i\norm{a_i}^2+|e^{a_i^{\mathsf T}x}-y_i|\norm{a_i}^2
 \le L_{0,i}+L_1\norm{\nabla f_i(x)}.
\]
Moreover, for $0\le t\le1$, by the Cauchy--Schwarz inequality,
\[
 |a_i^{\mathsf T}(x_0+ts-x_1)|\le L_1\norm s.
\]
Therefore,
\[
 e^{-L_1\norm s}e^{a_i^{\mathsf T}x_1}
 \le e^{a_i^{\mathsf T}(x_0+ts)}
 \le e^{L_1\norm s}e^{a_i^{\mathsf T}x_1}.
\]
Multiplying by the positive semidefinite matrix
$a_i a_i^{\mathsf T}$, integrating, and averaging over $i$,
we obtain~\eqref{eq:poisson-secant-comparison}.
\end{proof}

Applying the main theorem also requires
Assumptions~\ref{asm:minimizer}--\ref{asm:heterogeneity}.
The condition that no nonzero directions are missed is essential:
if $a_i^{\mathsf T}s=0$ with $a_i\ne0$, then the gradient secant is zero,
whereas the Hessian of the Poisson component is nonzero.
Moreover, two observed parallel derivatives do not prove that
the component is globally representable in the form $\phi_i(a_i^{\mathsf T}x)$.
Therefore, the recovered matrix can be used to choose the
heuristic proposal of YETI: the checks and the main iterations of the method
preserve its complexity bound.
Replacing the full sum of secants by a sample of a fixed number of components
would require a separate theorem on spectral approximation.

\subsection{Limiting the number of oracle queries made by the heuristic}
\label{sec:budget}

Let $B_0=n\lceil\log N\rceil$, let $G$ be the total number of oracle queries
in the main iterations, and let $H$ be the total number in the heuristic blocks.
Choose $\omega\ge0$. Before the next heuristic block, we allow
\begin{equation}\label{eq:adaptive-cheap-budget}
 B_k=n\left\lfloor
       \frac{\min\{B_0,B_0+\omega G-H\}}n
       \right\rfloor
\end{equation}
queries. If $B_k<n$, the block is skipped. Before each evaluation of the function and the full gradient,
$n$ queries are reserved. All oracle calls made by the heuristic, including
checks of the proposed points, are counted in $H$.

The proof of the following theorem reduces to a step-by-step verification of an invariant
and is given in Appendix~\ref{app:budget-proof}.

\proofdown{thm:adaptive-cheap-budget}
\begin{theorem}[number of heuristic queries]
\label{thm:adaptive-cheap-budget}
If $k$ heuristic blocks have been started by some moment of the run,
then at this moment
\begin{equation}\label{eq:budget-prefix-invariant}
 H\le\min\{kB_0,B_0+\omega G\}.
\end{equation}
Rule~\eqref{eq:adaptive-cheap-budget} preserves the upper bound of
Theorem~\ref{thm:cheap-upper}. Moreover,
\begin{equation}\label{eq:adaptive-budget-total}
 \mathbb E T\le B_0+(1+\omega)\mathbb E G.
\end{equation}
Here all final evaluations of the function and the full gradient are included in $G$ or $H$.
\end{theorem}

For $\omega=1$, beyond the first $B_0$ queries, the heuristic cannot make
more oracle calls than the main iterations have already made. The first block has all $B_0$ queries available,
so for the same values of the random variables it makes
the same queries as a block with the constant limit $B_0$.
This does not require any assumption on the usefulness of the heuristic.
Limiting the number of queries is essential: repeating the same point without
stopping and without passing to a main iteration can make
the number of oracle calls infinite, even if the function value does not increase.

\clearpage
\section{Detailed description of YETI}\label{app:implementation}
\subsection{History-based stopping criterion}

Evaluating the function and the full gradient at a point $p\in\mathbb R^d$ yields the affine minorant
$u\mapsto f(p)+\ip{\nabla f(p)}{u-p}$. Its minimum over $u\in B(0,R)$
equals
\begin{equation}\label{eq:history-full-lower}
 \widetilde f(p)
 =f(p)-\ip{\nabla f(p)}p-R\norm{\nabla f(p)}.
\end{equation}
This is a lower bound on the optimal value:
\[
 f^\ast\by{\ge}{\text{\ref{asm:components}}} f(p)+\ip{\nabla f(p)}{x^\ast-p}
 \by{\ge}{\text{\ref{asm:minimizer}}}
 f(p)-\ip{\nabla f(p)}p-R\norm{\nabla f(p)}.
\]
The second step also uses the Cauchy--Schwarz inequality.
If $p\in B(0,R)$, then the difference $f(p)-\widetilde f(p)$ is the
Frank--Wolfe gap~\citep{jaggi2013}; the bound $\widetilde f(p)\le f^\ast$
also holds for $p\notin B(0,R)$.
Denote by $\widetilde f_{\rm best}$ the maximum of all stored
lower bounds, and by $w_{\rm best}$ the best of the points of the ball $B(0,R)$
at which the function value is known. The method stops and returns
$w_{\rm best}$ if
\begin{equation}\label{eq:history-stop}
 f(w_{\rm best})-\widetilde f_{\rm best}\le\varepsilon/2.
\end{equation}
The point $w_{\rm best}$ is stored separately from $x$: under the accelerated update of $x$,
the function value does not necessarily decrease monotonically.
All evaluations of the function and the full gradient are performed by the subroutine \textsc{Evaluate}
(Algorithm~\ref{alg:pass}): at the points proposed by the heuristic,
at the centers, and at the anchor and test points of the solver. If the values and gradients of all components at a given point are already
known, repeated queries are not needed. To check~\eqref{eq:history-stop},
it suffices to store $\widetilde f_{\rm best}$, $w_{\rm best}$,
and $f(w_{\rm best})$; the heuristic may also store additional data for choosing points.

\begin{algorithm}[H]
\caption{\textsc{Evaluate}$(p)$: evaluation of the function and gradient, and accuracy check}\label{alg:pass}
\begin{algorithmic}[1]
\Require point $p$
\State query $\bigl(f_i(p),\nabla f_i(p)\bigr)$ for all $i$ whose answers at $p$ have not yet been obtained
\State $\widetilde f_{\rm best}\gets\max\bigl\{\widetilde f_{\rm best},\ f(p)-\ip{\nabla f(p)}p-R\norm{\nabla f(p)}\bigr\}$
\If{$\norm p\le R$ and $f(p)<f(w_{\rm best})$} \Comment{$f(w_{\rm best})=+\infty$ until the first such point}
 \State $w_{\rm best}\gets p$
\EndIf
\If{$f(w_{\rm best})-\widetilde f_{\rm best}\le\varepsilon/2$}
 \State terminate YETI with output $w_{\rm best}$
\EndIf
\end{algorithmic}
\end{algorithm}

\proofdown{lem:history-small-gradient}
\begin{lemma}[stopping at a small subgradient]
\label{lem:history-small-gradient}
Let $z\in B(0,R)$, let the function and the full gradient at $z$ be already computed, and let
$g=\nabla f(z)+v_z$, where $v_z\in N_{B(0,R)}(z)$. Then
\[
 f(z)-\widetilde f(z)\le2R\norm g.
\]
In particular, if $\norm g\le\varepsilon/(4R)$, then after $z$ is recorded
in the history, \eqref{eq:history-stop} holds.
\end{lemma}

The stopping test can also use component answers
obtained at different points. The corresponding lower bound and the way to
update it without new queries are given in
Appendix~\ref{app:component-history}. This option is used by the
heuristic in the experiments, but it is not needed for the complexity bound of YETI.

\subsection{Estimating the local curvature}

After the function and the full gradient at $y$ are evaluated, the quantity
\begin{equation}\label{eq:observed-curvature}
 \ell(y)=\frac1n\sum_{i=1}^n
       \bigl(L_{0,i}+L_1\norm{\nabla f_i(y)}\bigr)
\end{equation}
is known.
\proofdown{lem:observed-curvature}
\begin{lemma}[curvature estimate via the optimality gap]\label{lem:observed-curvature}
Under Assumptions~\ref{asm:components}--\ref{asm:minimizer}
and~\ref{asm:heterogeneity}, for every point $y$,
\begin{equation}\label{eq:curvature-gap}
 \begin{gathered}
 \norm{\nabla f(y)}
 \le L_1\bigl(\delta+f(y)-f^\ast\bigr)
 +\sqrt{2L_0\bigl(\delta+f(y)-f^\ast\bigr)},\\
 L_1\norm{\nabla f(y)}\le\ell(y),\\
 \ell(y)\le\tfrac32L_0+2L_1^2\bigl(\delta+f(y)-f^\ast\bigr).
 \end{gathered}
\end{equation}
Moreover, for all $x,y$,
\begin{equation}\label{eq:curvature-transfer}
 \begin{gathered}
 \ell(x)\le e^{L_1\norm{x-y}}\ell(y),\\
 \norm{\nabla^2f(x)}\le\ell(x).
 \end{gathered}
\end{equation}
In particular, for every $r>0$, the norm of the Hessian of the component $f_i$ on $B(y,3r)$
does not exceed
\begin{equation}\label{eq:component-local-smoothness}
 e^{3L_1r}\bigl(L_{0,i}+L_1\norm{\nabla f_i(y)}\bigr).
\end{equation}
\end{lemma}

\begin{corollary}[lower bound on the optimality gap via curvature]\label{cor:curvature-lower-gap}
Evaluating the function and the full gradient at $y$ yields the lower bound on the optimality gap
\[
 f(y)-f^\ast\ge
 \frac{\ell(y)-\tfrac32L_0-2L_1^2\delta}{2L_1^2}.
\]
If the right-hand side is negative, the bound is trivial.
\end{corollary}

\subsection{Parameters of YETI}
The state of the method consists of points $x,v\in B(0,R)$, a number
$\mathcal A>0$, and a parameter $h>0$. The number $\mathcal A$ plays the same role
as the sum of the coefficients $A_k$ in accelerated methods with estimate
sequences: we prove below that always
$f(x)-f^\ast\le R^2/\mathcal A$, and that $\mathcal A$ does not decrease
and increases under the update~\eqref{eq:damped-update}.
Unlike the sum $A_k$ introduced in the notation, the initial value of
$\mathcal A$ is positive. The radius $r$ is chosen so that on the small ball the curvature changes
by at most a constant factor. The proximal parameter
is bounded from below for three reasons, related to the weight of the
accelerated step, the conditioning of the local subproblem, and stopping.
We use the following values:
\begin{equation}\label{eq:outer-parameters}
 \begin{aligned}
 &r=\min\left\{R,\frac1{2L_1}\right\},\\
 &\lambda_\downarrow(\mathcal A)
 =\max\left\{\frac{1-\frac12(r/R)^{2/3}}{\frac14(r/R)^{4/3}\mathcal A},\
       \frac1n\Bigl(\frac32L_0+2L_1^2\Bigl(\delta+\frac{R^2}{\mathcal A}\Bigr)\Bigr),\
       \frac\varepsilon{8rR}\right\},\\
 &\lambda_{\uparrow}=\lambda_\downarrow\Bigl(\frac{R^2}{2\Delta}\Bigr).
 \end{aligned}
\end{equation}
Here $r$ is the radius of the small ball, $\lambda_\downarrow(\mathcal A)$ is the
lower limit of the proximal parameter, and $\lambda_{\uparrow}$ is the upper one.
The three expressions under the maximum correspond to these reasons: the first prevents
the weight $\alpha$ of the accelerated step from exceeding $\frac12(r/R)^{2/3}$, the second makes
the local subproblem well conditioned, and the third is needed for stopping.

The method starts with $x=v=0$, $\mathcal A=R^2/(2\Delta)$, $h=\lambda_{\uparrow}$, and an empty
query history.

\subsection{Main iteration of YETI}
In the state $(x,v,\mathcal A,h)$, we compute
\begin{equation}\label{eq:outer-step-setup}
 \begin{aligned}
 &\lambda=\min\bigl\{\max\{h,\lambda_\downarrow(\mathcal A)\},\lambda_{\uparrow}\bigr\},\\
 &a_\lambda=\frac{1+\sqrt{1+4\lambda\mathcal A}}{2\lambda},\\
 &\alpha=\frac{a_\lambda}{\mathcal A+a_\lambda},\\
 &y=(1-\alpha)x+\alpha v.
 \end{aligned}
\end{equation}
The number $a_\lambda$ is the positive root of the equation
$\lambda a_\lambda^2=\mathcal A+a_\lambda$; therefore,
\[
 \begin{gathered}
 \lambda\mathcal A\alpha^2=1-\alpha,\\
 \lambda a_\lambda=1/\alpha.
 \end{gathered}
\]
Using the values of the component gradients at $y$, we check the curvature condition
\begin{equation}\label{eq:simple-curvature-test}
 \ell(y)\le e^{L_1r}\Bigl(\frac32L_0+2L_1^2\Bigl(\delta+\frac{R^2}{\mathcal A}\Bigr)\Bigr).
\end{equation}
If this condition is violated, then by Corollary~\ref{cor:curvature-lower-gap}
the optimality gap at $y$ is large, and the halfspace
\begin{equation}\label{eq:cut-update}
 \mathcal H_y=\left\{w:\ \ip{\nabla f(y)}{w-y}
 \le-\frac{\ell(y)-\frac32L_0-2L_1^2\delta}{2L_1^2}\right\}
\end{equation}
contains the minimizer. In this case, we perform the update
\[
 \begin{aligned}
 v^+&=\Pi_{B(0,R)}\Pi_{\mathcal H_y}v,\\
 h^+&=2\lambda,
 \end{aligned}
\]
under which $x$ and $\mathcal A$ do not change. Such a step is called
an up-step.

If condition~\eqref{eq:simple-curvature-test} holds, we run the
local solver (Algorithm~\ref{alg:local-solver}), which uses the already
obtained values and gradients of the components at $y$. It returns a point $z$, a subgradient $g$
of the problem on the ball $B(0,R)$ at $z$, and a number $\beta\ge0$.
If $g$ is small, the method returns $z$. Otherwise, setting $\gamma=\lambda/(\lambda+\beta)$,
we reduce the weight of the accelerated update by multiplying it by $\gamma$,
following the scheme of~\citep{carmon2022adaptive}:
\begin{equation}\label{eq:damped-update}
 \begin{aligned}
 \mathcal A^+&=\mathcal A+\gamma a_\lambda,\\
 x^+&=\frac{(1-\gamma)\mathcal A x+
                 \gamma(\mathcal A+a_\lambda)z}{\mathcal A^+},\\
 v^+&=\Pi_{B(0,R)}(v-\gamma a_\lambda g),\\
 h^+&=\begin{cases}2\lambda,&\beta>0,\\\lambda/2,&\beta=0.\end{cases}
 \end{aligned}
\end{equation}
The step with $\beta>0$ is also called an up-step, and the step with $\beta=0$
a down-step. If $\mathcal A^+\ge2R^2/\varepsilon$, the method returns $x^+$.

\subsection{Heuristic block}\label{sec:heuristic-block}

The heuristic uses the known parameters of the problem, the history of oracle
answers, and its own random variables. Each block returns
control after finitely many computational steps. The limit $n\lceil\log N\rceil$
applies to all queries, including queries at rejected points and the initial
evaluation of the function and the full gradient at $x$. Before the heuristic is called,
$n$ queries are reserved for checking the proposal.

A point $z$ is accepted if $z\in B(0,R)$ and $f(z)\le f(x)$;
in this case only $x$ changes. Stopping is possible by criterion~\eqref{eq:history-stop}.
The decisions made by the heuristic are measurable. The indices of the local solver use
fresh independent random variables; given the history, the indices have
the distribution from Section~\ref{sec:local-solver}.

\begin{algorithm}[H]
\caption{\textsc{Block}$(x)$: heuristic block}\label{alg:block}
\begin{algorithmic}[1]
\Require point $x\in B(0,R)$, heuristic
\State $b=n\lceil\log N\rceil$ \Comment{each query decreases $b$ by one}
\State \Call{Evaluate}{$x$}
\While{$b\ge n$ and the heuristic has not terminated the block}
 \State the heuristic makes at most $b-n$ queries and proposes $z\in B(0,R)$
 \State \Call{Evaluate}{$z$}
 \If{$f(z)\le f(x)$}
  \State $x\gets z$
 \EndIf
\EndWhile
\State \Return $x$
\end{algorithmic}
\end{algorithm}

\section{Local solution by SVRG}\label{app:solver-implementation}
\subsection{Local subproblem}
Let $y\in B(0,R)$, $\lambda>0$, and $r>0$. Consider the functions
\begin{equation}\label{eq:local-problem}
 \begin{aligned}
 P(z)&=f(z)+\frac\lambda2\norm{z-y}^2,\\
 P_i(z)&=f_i(z)+\frac\lambda2\norm{z-y}^2,
 \end{aligned}
\end{equation}
the set $\mathcal K=B(0,R)\cap B(y,r)$, and the point
\[
 u=\argmin_{z\in\mathcal K}P(z).
\]
The set $\mathcal K$ is compact and convex, and $P$ is strongly convex with parameter
$\lambda$; hence $u$ exists and is unique. Evaluating the function and the full gradient at $y$
allows us to find the numbers
\begin{equation}\label{eq:local-smoothness}
 L_i=e^{3L_1r}\bigl(L_{0,i}+L_1\norm{\nabla f_i(y)}\bigr)+\lambda
\end{equation}
and their average
\[
 \overline L=\frac1n\sum_iL_i=e^{3L_1r}\ell(y)+\lambda.
\]
By Lemma~\ref{lem:observed-curvature}, the norm of the Hessian of $P_i$ on $B(y,3r)$
does not exceed $L_i$.

The following lemma describes the projection onto $\mathcal K$ and the local
cocoercivity of the components; its proof is given
in Appendix~\ref{app:local-proofs}.

\proofdown{lem:local-projection}
\begin{lemma}[projection onto the intersection of balls]\label{lem:local-projection}
If $z=\Pi_{\mathcal K}(b)$, then there exist $\xi,\zeta\ge0$ such that
\begin{equation}\label{eq:two-ball-normal}
 \begin{gathered}
 b-z=\xi z+\zeta(z-y),\\
 \xi>0\Longrightarrow\norm z=R,\\
 \zeta>0\Longrightarrow\norm{z-y}=r.
 \end{gathered}
\end{equation}
For $x,t\in B(y,r)$,
\begin{equation}\label{eq:local-bregman}
 \norm{\nabla P_i(x)-\nabla P_i(t)}^2
 \le2L_i\bigl(P_i(x)-P_i(t)-\ip{\nabla P_i(t)}{x-t}\bigr).
\end{equation}
In particular, the gradient of $P$ on $B(y,r)$ is $1/\overline L$-cocoercive.
\end{lemma}

\subsection{SVRG epoch}
Within a single call of the solver, the quantities $\lambda$, $y$, and $L_i$ do not change. The index $i$ is chosen with probability $L_i/(n\overline L)$.
The step size and the epoch length are
\begin{equation}\label{eq:svrg-parameters}
 \begin{aligned}
 h_s&=\frac1{6\overline L},\\
 M&=\left\lceil\frac{21\overline L}{5\lambda}\right\rceil.
 \end{aligned}
\end{equation}
Since $L_i\ge\lambda>0$, all probabilities $L_i/(n\overline L)$ are positive and sum to one.
At the anchor point $s\in\mathcal K$, the values and gradients of all components are known. From the current point
$x_j$, a step is made with a new index $i_j$, conditionally independent of the previous ones:
\begin{equation}\label{eq:svrg-step}
 \begin{aligned}
 &v_{i_j}=\nabla P(s)+\frac{\nabla P_{i_j}(x_j)-\nabla P_{i_j}(s)}{L_{i_j}/\overline L},\\
 &x_{j+1}=\Pi_{\mathcal K}(x_j-h_sv_{i_j}).
 \end{aligned}
\end{equation}
The next epoch starts at the point $w^+=x_M$, and the new anchor point is
\begin{equation}\label{eq:svrg-average}
 s^+=\frac{\sum_{j=0}^{M-1}(1-h_s\lambda)^{M-1-j}x_j}
          {\sum_{j=0}^{M-1}(1-h_s\lambda)^j}.
\end{equation}
The first epoch of a call starts with $w=s=y$. The component gradients
at the anchor point are stored, so one inner step requires
one query.

\begin{algorithm}[H]
\caption{\textsc{Epoch}$(w,s)$: SVRG epoch for the local subproblem}\label{alg:svrg-epoch}
\begin{algorithmic}[1]
\Require points $w,s\in\mathcal K$, answers of all components at $s$,
$P_i$, $\mathcal K$, $L_i$, $\overline L$, $h_s$, $M$ from Algorithm~\ref{alg:local-solver}
\State $x_0=w$
\For{$j=0,\dots,M-1$}
 \State choose the index $i_j=i$ with probability $L_i/(n\overline L)$ independently of the past
 \State query $\nabla f_{i_j}(x_j)$
 \State $v_j=\nabla P(s)+\dfrac{\overline L}{L_{i_j}}\bigl(\nabla P_{i_j}(x_j)-\nabla P_{i_j}(s)\bigr)$
 \State $x_{j+1}=\Pi_{\mathcal K}(x_j-h_sv_j)$
\EndFor
\State \Return $\Bigl(x_M,\ \dfrac{\sum_{j=0}^{M-1}(1-h_s\lambda)^{M-1-j}x_j}{\sum_{j=0}^{M-1}(1-h_s\lambda)^j}\Bigr)$
\end{algorithmic}
\end{algorithm}

\subsection{Stopping the local solver}
After an epoch, the full gradient at $w$ is computed. First, the
tuple
\[
 \begin{gathered}
 z=w,\\
 \beta=0,\\
 g=\nabla f(w),\\
 e=\nabla P(w)
 \end{gathered}
\]
is checked.
If it is not accepted, we compute the point and the vectors
\begin{equation}\label{eq:certificate-projection}
 \begin{aligned}
 z&=\Pi_{\mathcal K}\bigl(w-\nabla P(w)/\overline L\bigr),\\
 e&=\nabla P(z)-\nabla P(w)+\overline L(w-z),\\
 g&=\nabla f(z)+\overline L\xi z,
 \end{aligned}
\end{equation}
and the number
\[
 \beta=\overline L\zeta,
\]
where $\xi,\zeta$ are the coefficients from~\eqref{eq:two-ball-normal} for precisely
the projection in the first line. Computing $e$ and $g$ requires the full gradient
at $z$. In both cases,
\begin{equation}\label{eq:certificate-properties}
 \begin{gathered}
 g\in\partial(f+I_{B(0,R)})(z),\\
 e=g-(\lambda+\beta)(y-z),\\
 \beta>0\Longrightarrow\norm{y-z}=r,
 \end{gathered}
\end{equation}
where $I_{B(0,R)}$ is the indicator function of the convex set; its
subdifferential is the normal cone. The tuple is accepted if
at least one of the conditions
\begin{equation}\label{eq:certificate-test}
 \begin{gathered}
 \norm g\le\frac\varepsilon{4R},\\
 \norm e\le\frac45(\lambda+\beta)\norm{y-z}
 \end{gathered}
\end{equation}
holds. The first condition is checked before the second.

The complete solver is given in Algorithm~\ref{alg:local-solver}.
If no tuple is accepted after an epoch, the function and the full gradient are evaluated
at the new anchor point, and the next epoch begins.

\begin{algorithm}[H]
\caption{\textsc{Solver}$(y,\lambda)$: approximate solution of the local subproblem}\label{alg:local-solver}
\begin{algorithmic}[1]
\Require center $y\in B(0,R)$ at which the answers of all components are known, parameter $\lambda>0$
\State $P_i(z)=f_i(z)+\frac\lambda2\norm{z-y}^2$, $i=1,\dots,n$
\State $P(z)=f(z)+\frac\lambda2\norm{z-y}^2$
\State $\mathcal K=B(0,R)\cap B(y,r)$
\State $L_i=e^{3L_1r}\bigl(L_{0,i}+L_1\norm{\nabla f_i(y)}\bigr)+\lambda$, $i=1,\dots,n$
\State $\overline L=\frac1n\sum_iL_i$
\State $h_s=1/(6\overline L)$
\State $M=\lceil21\overline L/(5\lambda)\rceil$
\State a tuple $(z,g,\beta)$ is \emph{accepted} if $\norm g\le\dfrac\varepsilon{4R}$ or $\norm{g-(\lambda+\beta)(y-z)}\le\dfrac45(\lambda+\beta)\norm{y-z}$
\State $w=s=y$
\Loop
 \State $(w,s)\gets\Call{Epoch}{w,s}$
 \State \Call{Evaluate}{$w$}
 \If{the tuple $\bigl(w,\nabla f(w),0\bigr)$ is accepted}
  \State \Return $\bigl(w,\nabla f(w),0\bigr)$
 \EndIf
 \State $z=\Pi_{\mathcal K}\bigl(w-\nabla P(w)/\overline L\bigr)$
 \State $\xi,\zeta\ge0$ are the multipliers of the normals to the two balls at $z$:
 \Statex \hspace{\algorithmicindent}\quad$w-\nabla P(w)/\overline L-z=\xi z+\zeta(z-y)$,
 \Statex \hspace{\algorithmicindent}\quad$\xi>0\Longrightarrow\norm z=R$,
 \Statex \hspace{\algorithmicindent}\quad$\zeta>0\Longrightarrow\norm{z-y}=r$
 \State \Call{Evaluate}{$z$}
 \If{the tuple $\bigl(z,\nabla f(z)+\overline L\xi z,\overline L\zeta\bigr)$ is accepted}
  \State \Return $\bigl(z,\nabla f(z)+\overline L\xi z,\overline L\zeta\bigr)$
 \EndIf
 \State \Call{Evaluate}{$s$}
\EndLoop
\end{algorithmic}
\end{algorithm}

\proofdown{lem:local-certificate}
\begin{lemma}[number of epochs until the solver stops]\label{lem:local-certificate}
For the tuple \eqref{eq:certificate-projection}, the bounds
\begin{equation}\label{eq:sharp-residual}
 \begin{gathered}
 \norm{z-u}\le\norm{w-u},\\
 \norm e\le\overline L\norm{w-u}
 \end{gathered}
\end{equation}
hold. Therefore, if
$\norm{w-u}\le\frac{4\lambda}{5\overline L+4\lambda}\norm{y-u}$,
the tuple is accepted.
Let $T_{\rm ep}$ be the number of epochs until the first accepted tuple. Then for $j\ge1$,
\begin{equation}\label{eq:epoch-tail}
 \Pr(T_{\rm ep}>j\mid\text{history of the call})
 \le\min\{1,\Gamma 2^{-j}\},
\end{equation}
where
\[
 \Gamma=\left(1+\frac{2\overline L}{3\lambda}\right)
 \left(1+\frac{5\overline L}{4\lambda}\right)^2.
\]
In particular, the solver terminates almost surely and makes a finite
expected number of oracle calls.
\end{lemma}

\clearpage
\section{Experimental details}\label{app:experiments}

\subsection{Experimental setup}\label{sec:exp-protocol}

\paragraph{Problems.}
We use three groups of problems:
\begin{itemize}
\item Chains S and X: $n=16$, $T\in\{8,16,32,64\}$.
These are numerical versions of the families from Section~\ref{sec:lower}
with tabulated integrals, without hidden coordinates and without radial
extension. They allow us to compare methods as $T$ increases;
these experiments do not test the lower bound for arbitrary methods.
\item Poisson models $f_i(x)=w_i e^{a_i^{\mathsf T}x}-y_i a_i^{\mathsf T}x$,
where $w_i>0$, $y_i\ge0$, $L_{0,i}=y_i\norm{a_i}^2$,
$L_1=\max_i\norm{a_i}$.
In the synthetic problems, $d=50$, $n\in\{2000,20000,100000\}$,
and the condition number of the feature covariance is $10^4$.
On the Bike Sharing data, $n=17379$, $d=48$;
on freMTPL2freq, $n=678013$, $d=68$.
\item 65 problems: 13 families with five random instances each.
This series includes rotated, ill-conditioned (up to $10^8$),
nonlinear, and sparse Poisson models, high-dimensional problems,
sums of quartic functions, and regularized logistic regression.
The implementation settings were fixed before these experiments were run;
the problems used during its development are not included in the series.
\end{itemize}

\paragraph{Metric.}
We count all queries of a method, including queries for constructing
the initial metric, computing full gradients, stochastic corrections,
and line search.
One pass equals $n$ queries.
The minimum on the chains is known from the construction; on the other problems,
$f_{\rm ref}$ is found numerically. The reference solution is also used to compute
the bounds $R,\Delta,\delta$ passed to YETI as problem parameters.
Queries for computing the reference and for external diagnostics are not counted.
The number of queries to reach the tolerance is determined by the first checked
point satisfying the tolerance; on the 65 problems, we consider all points
with a complete set of component answers, including trial points of the line search.
The plots show the smallest optimality gap attained up to a given moment.
The triggering of a method's own stopping criterion is measured separately.

\paragraph{Methods.}
On the 65 problems, we compare YETI, L-BFGS with memory 50,
Cov-LBFGS $1/4$ and $1/2$, BFGS, Nesterov's method with line search and restart,
GRAAL with stabilization of the curvature estimate, and SVRG
\citep{liu1989lbfgs,nocedal2006,borodich2026graal,johnson2013svrg}.
Cov-LBFGS uses the coordinate transformation $C^{-1/4}$ or $C^{-1/2}$,
where $C$ is the centered covariance of the initial component gradients;
on the cut-off subspace, the pseudoinverse is zero.
On the main test set, we additionally consider nonlinear CG,
AdGD, clipped gradient descent, SAGA, and Katyusha
\citep{malitsky2020adgd,zhang2019gradient,defazio2014saga,allenzhu2018katyusha}.
For the variance reduction methods, the local smoothness estimates
are updated using the available component gradients.

\paragraph{Run limits.}
On the chains, the baseline methods are limited to 20000 passes,
on the Poisson models to 1000, and on the 65 problems to
3000 passes and 1800 seconds.
The main tables report results for one instance of each problem configuration.
In the extended series, all runs are counted, including failed ones;
the family median is reported only if all five instances succeed.

\paragraph{Tested variant of YETI.}
For collinear component gradients, the heuristic uses the secants
from Proposition~\ref{prop:glm-secant}; in the remaining cases, it uses
L-BFGS in the covariance metric with the power chosen from the observed secant.
When the history is insufficient, standard L-BFGS is used.
On large samples, the initial secant matrix is built from a subsample of
up to 4096 rows. After 32 evaluations of the full gradient or after an unsuccessful block,
a switch to the full history of component answers is possible.
The dense covariance requires $O(nd^2)$ operations,
and its spectral decomposition $O(d^3)$.
The implementation limits the fraction of the accelerated step and the number of inner
calls, and changes the settings of the local solver.
Theorem~\ref{thm:main} is proved for Algorithm~\ref{alg:cheap};
these modifications require a separate analysis.

\subsection{Results by family}

Figure~\ref{fig:exp-profiles} compares the methods on all 65 problems,
and Table~\ref{tab:exp-series} shows the medians over the five instances of each family.
In particular, on the quartic problems with $n<d$, Cov-LBFGS does not reach
the tolerance on all instances because directions are lost when inverting
the degenerate covariance matrix.

\begin{figure}[!htbp]
 \centering
 \input{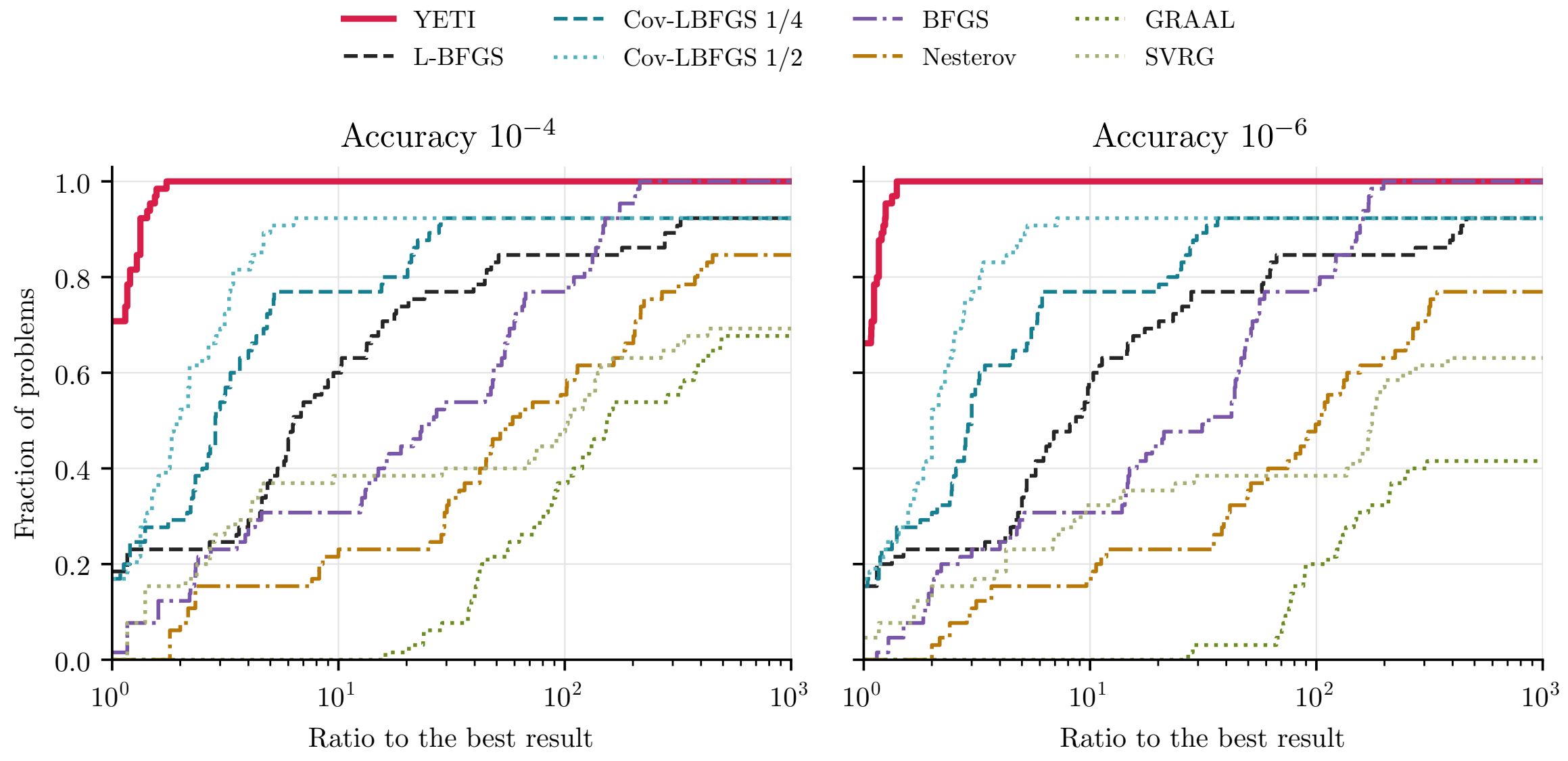}
 \caption{Fraction of the 65 problems on which the number of queries does not exceed
 the best result by more than a factor of $\tau$ \citep{dolan2002profiles}. At $\tau=1$,
 all first places are counted, including shared ones; failed runs remain
 in the denominator.}
 \label{fig:exp-profiles}
\end{figure}

\begin{table}[!htbp]
 \centering\small
 \setlength{\tabcolsep}{2.3pt}
 \begin{tabular}{@{}lrrrrrrrr@{}}
 \toprule
 Family & YETI & L-BFGS & Cov-LBFGS & Cov-LBFGS & BFGS & Nesterov & GRAAL & SVRG\\
 & & (50) & $1/4$ & $1/2$ & & & &\\
 \midrule
Rotated, $10^4$ & \textbf{6.0} & 163 & 33 & 15 & 337 & 1969 & --- & --- \\
Rotated, $10^8$ & \textbf{6.0} & 2617 & 155 & 15 & 725 & --- & --- & --- \\
High dimension & \textbf{7.0} & 440 & 42 & 37 & 1123 & --- & --- & --- \\
Strong nonlinearity & \textbf{8.0} & 45 & 24 & 27 & 40 & 323 & 1080 & 34 \\
Rare features & \textbf{11} & 57 & 27 & 19 & 229 & 670 & 1614 & 94 \\
Rare, $n=5\cdot10^4$ & 27 & 104 & 30 & \textbf{25} & 464 & 1230 & --- & 190 \\
$n=5\cdot10^4$, cond.\ $10^8$ & \textbf{8.0} & --- & 237 & 10 & 1364 & --- & --- & --- \\
Quartics $32\times48$ & \textbf{101} & \textbf{101} & --- & --- & 211 & 1121 & --- & --- \\
Quartics $256\times32$ & 13 & 73 & \textbf{12} & 29 & 177 & 1091 & --- & --- \\
Quartics, $\nu=0.01$ & \textbf{7.0} & 8.0 & \textbf{7.0} & 11 & 9.0 & 22 & 474 & \textbf{7.0} \\
Quartics, $\nu=0.1$ & 7.0 & \textbf{6.0} & 7.0 & 11 & 12 & 12 & 444 & 10 \\
Logistic, rotated & 10 & 88 & 20 & \textbf{9.0} & 392 & 1011 & --- & 1519 \\
Logistic, rare & 10 & 132 & 28 & \textbf{9.0} & 470 & 2403 & --- & 1591 \\
 
 \bottomrule
 \end{tabular}
 \caption{Median number of passes to relative accuracy $10^{-6}$
 over five instances. ``---'' indicates failure on at least one instance;
 the minimum before rounding is shown in bold. $\kappa$ is the condition number
 of the feature covariance, and $n\times d$ are the number of components and the dimension.}
 \label{tab:exp-series}
\end{table}

Results on the main test set are given in Tables~\ref{tab:exp-chains}
and~\ref{tab:exp-glm}. In additional runs of YETI on five
synthetic Poisson instances, the median numbers of passes are
4, 5, and 5 for $n=2000$, $20000$, and $100000$;
the ranges are 4--5, 5--7, and 5--8, respectively.

\begin{table}[!htbp]
 \centering\small
 \setlength{\tabcolsep}{4pt}
 \begin{tabular}{@{}lrrrrrrrr@{}}
 \toprule
 & \multicolumn{4}{c}{S, length $T$} & \multicolumn{4}{c}{X, length $T$}\\
 \cmidrule(lr){2-5}\cmidrule(l){6-9}
 Method & 8 & 16 & 32 & 64 & 8 & 16 & 32 & 64\\
 \midrule
YETI & \textbf{22} & 44 & \textbf{80} & \textbf{158} & \textbf{21} & \textbf{39} & \textbf{81} & \textbf{151} \\
L-BFGS & \textbf{22} & \textbf{43} & 87 & 181 & \textbf{21} & 41 & 88 & 175 \\
BFGS & 39 & 79 & 157 & 313 & 27 & 55 & 111 & 217 \\
Katyusha & 223 & 436 & 868 & 1726 & 265 & 514 & 1009 & 2005 \\
Nesterov & 82 & 182 & 377 & 763 & 159 & 416 & 873 & 1850 \\
GRAAL & 756 & 1472 & 2906 & 5814 & 496 & 1096 & 2610 & 5504 \\
AdGD & 149 & 156 & 1005 & 4640 & 392 & 1650 & 6538 & --- \\
SVRG & 3658 & 13699 & --- & --- & 4222 & 15676 & --- & --- \\
 
 \bottomrule
 \end{tabular}
 \caption{Number of passes to relative accuracy $10^{-2}$
 on the chains. The minimum before rounding is shown in bold;
 ``---'' indicates that the tolerance was not reached within the allotted number of passes.}
 \label{tab:exp-chains}
\end{table}

\begin{table}[!htbp]
 \centering\small
 \setlength{\tabcolsep}{4pt}
 \begin{tabular}{@{}crrrrrrr@{}}
 \toprule
 & & \multicolumn{3}{c}{To $10^{-2}\widehat\Delta$}
 & \multicolumn{3}{c}{To the stopping criterion}\\
 \cmidrule(lr){3-5}\cmidrule(l){6-8}
 Family & $T$ & YETI & No & Random & YETI & No & Random\\
 & & & proposals & proposals & & proposals & proposals\\
 \midrule
S & 8 & 22 & 898 & 19290 & 27 & 2343 & 46145 \\
S & 16 & 44 & 2299 & 43590 & 47 & 6329 & 113107 \\
X & 8 & 21 & 1218 & 63896 & 33 & 3837 & 177593 \\
X & 16 & 39 & 1968 & 106768 & 55 & 7418 & --- \\
 
 \bottomrule
 \end{tabular}
 \caption{Number of passes for three variants of the heuristic.
 The method's own stopping criterion uses $\varepsilon=10^{-3}\widehat\Delta$.
 ``---'' indicates that it was not satisfied within 200000 passes.}
 \label{tab:exp-insurance}
\end{table}

\begin{table}[!htbp]
 \centering\small
 \begin{tabular}{@{}lrrrrr@{}}
 \toprule
 & \multicolumn{3}{c}{Synthetic problems, $n$}
 & Bike & freMTPL2\\
 \cmidrule(lr){2-4}
 Method & $2000$ & $20000$ & $100000$ & Sharing & freq\\
 \midrule
YETI & \textbf{4.0} & \textbf{5.0} & \textbf{5.0} & \textbf{9.0} & 7.0 \\
L-BFGS & 197 & 168 & 99 & 44 & 22 \\
Cov-LBFGS $1/4$ & 22 & 23 & 16 & 21 & 12 \\
BFGS & 353 & 325 & 253 & 26 & 116 \\
Nonlinear CG & 705 & 707 & 359 & 108 & 102 \\
Nesterov & --- & 775 & 508 & 242 & 77 \\
GRAAL & --- & --- & 984 & 642 & 396 \\
AdGD & --- & --- & --- & 937 & 86 \\
Clipped GD & --- & --- & --- & --- & --- \\
SVRG & --- & 145 & 16 & 10 & 10 \\
SAGA & --- & 906 & 96 & 32 & \textbf{5.0} \\
Katyusha & 394 & 133 & 43 & 22 & 10 \\
 
 \bottomrule
 \end{tabular}
 \caption{Number of passes to relative accuracy $10^{-4}$
 on the Poisson models. Notation as in Table~\ref{tab:exp-chains}.}
 \label{tab:exp-glm}
\end{table}
\FloatBarrier

\subsection{Stopping test}

On the 65 problems with $\varepsilon=10^{-6}\widehat\Delta$,
at least 90\% of the runs terminate within 37 passes.
In 63 of the 65 runs, YETI stops already in the first heuristic
block. The main proximal step was needed only on two instances
with rare features at $n=50000$: it was performed 4 and 8 times.
On the main test set, all 45 runs (24 on the chains and 21 on regression)
terminated before the first such step.
Figure~\ref{fig:exp-certificate} shows two additional runs
with a stricter tolerance.

\begin{figure}[!htbp]
 \centering
 \input{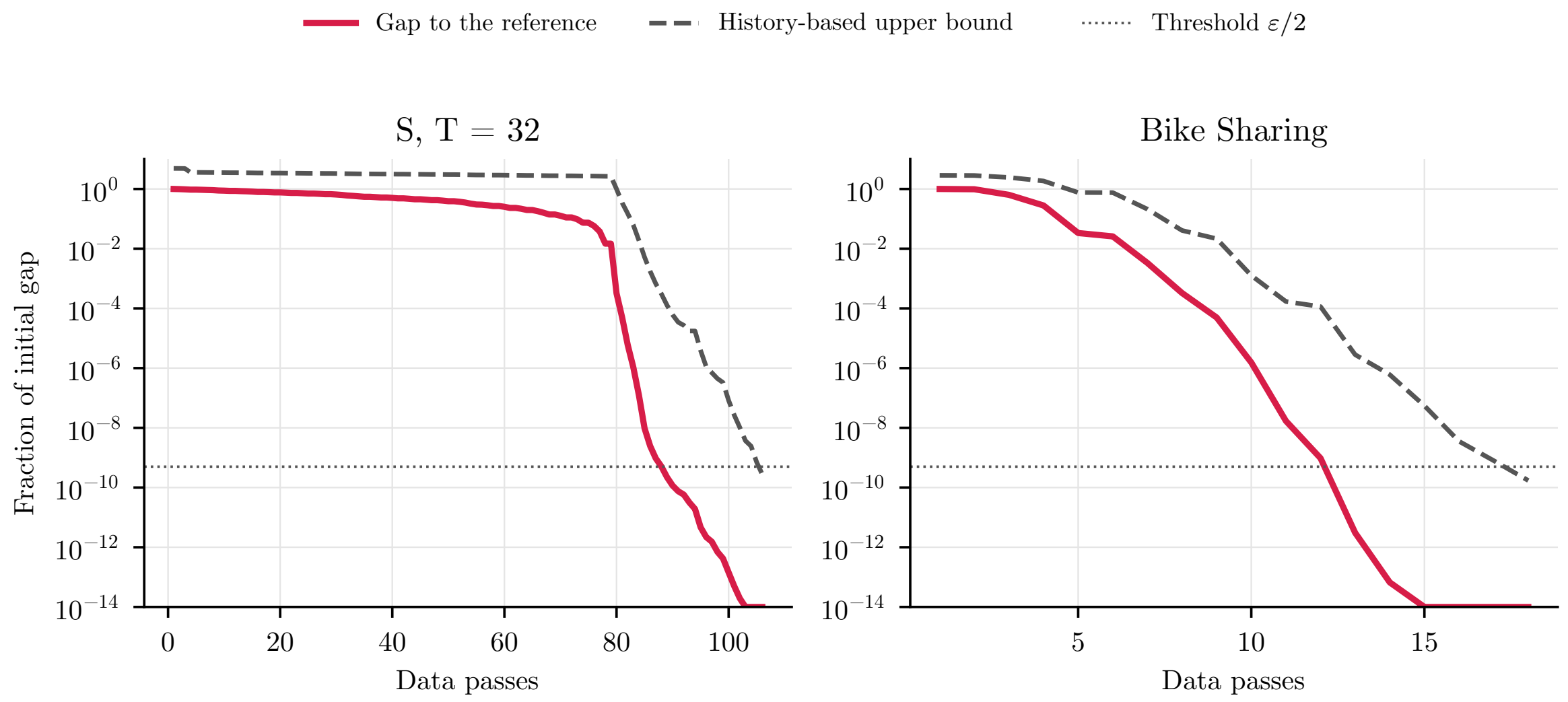}
 \caption{Optimality gap relative to the reference value and its history-based upper bound
 at $\varepsilon=10^{-9}\widehat\Delta$. All quantities are divided by
 $\widehat\Delta$; the dotted line indicates the threshold $\varepsilon/2$.
 Values below $10^{-14}$ are shown at the lower boundary of the axis.}
 \label{fig:exp-certificate}
\end{figure}

Table~\ref{tab:exp-insurance} compares the standard heuristic,
a block consisting of a single checking pass without proposals, and a block that proposes
random points of the feasible ball. In the last variant, each block spends
the entire allowed limit: 243--264 passes.
The overall limit is 60000 passes in the mode without proposals
and 200000 with random proposals.

The tolerance $10^{-2}\widehat\Delta$ is reached in all runs.
Disabling the heuristic proposals increases the number of passes to the tolerance
by a factor of about 41--58, and random proposals increase it by a further factor of about
19--54 relative to the mode without proposals. Thus,
in these experiments, the main gain comes from the heuristic. However, its matrix
operations can be expensive, so a gain in the number of queries
does not necessarily imply the same gain in time.
On X with $T=16$, the method's own stopping criterion is not satisfied within the allotted
number of passes. Limiting the block length does not eliminate the cost
of useless proposals.
\FloatBarrier

\clearpage
\section{Summary of known bounds}\label{app:literature}
\enlargethispage{2\baselineskip}
\relpenalty=500 \binoppenalty=700
\providecommand\yes{\textcolor{green!55!black}{\ding{51}}}
\providecommand\no{\textcolor{red!80!black}{\ding{55}}}

\begin{center}
\footnotesize
\captionof{table}{Finite-sum minimization: number of calls to the component oracle.}
\label{tab:lit-sum}
\renewcommand\arraystretch{1.0}\setlength\tabcolsep{4pt}
\begin{tabularx}{\linewidth}{@{}P{.18\linewidth}P{.13\linewidth}cP{.14\linewidth}Y@{}}
\toprule
Method and work & Smoothness & Convex & Criterion & Bound\\
\midrule
\multicolumn{5}{@{}l}{\textit{Finite sum, $L$-smooth components}}\\
SAGA\newline\citet{defazio2014saga} & $L$ & \yes & $\E f(\hat x)-f^\ast\le\varepsilon$
 & $\displaystyle O\Bigl(n+\frac{nF_0+LR^2}{\varepsilon}\Bigr)$\\
Katyusha with AdaptReg\textsuperscript{\dag}\newline\citet{allenzhu2018katyusha} & $L$ & \yes & $\E f(\hat x)-f^\ast\le\varepsilon$
 & $\displaystyle O\Bigl(n\log\frac{F_0}{\varepsilon}+\sqrt{\frac{nLR^2}{\varepsilon}}\Bigr)$\\
Varag\textsuperscript{\ddag}\newline\citet{lan2019varag} & $L$ & \yes & $\E f(\hat x)-f^\ast\le\varepsilon$
 & $\displaystyle O\Bigl(n+n\log\min\Bigl\{n,\frac{D_0}{\varepsilon}\Bigr\}+\sqrt{\frac{nD_0}{\varepsilon}}\Bigr)$\\
Lower bound\textsuperscript{\S}\newline\citet{woodworth2016finite} & $L$ & \yes & $\E f(\hat x)-f^\ast\le\varepsilon$
 & $\displaystyle \Omega\Bigl(n+\sqrt{\frac{nLR^2}{\varepsilon}}\Bigr)$\\
\midrule
\multicolumn{5}{@{}l}{\textit{One function ($n=1$), generalized smoothness, full-gradient oracle}}\\
Accelerated method\textsuperscript{\P}\newline\citet{tyurin2026accelerated} & sym.\ $(L_0,L_1)$ of $f$ & \yes & $f(\hat x)-f^\ast\le\varepsilon$
 & $\displaystyle O\Bigl(\sqrt{\frac{L_0R^2}{\varepsilon}}+L_1R\log_+\min\Bigl\{\frac{2L_1^2F_0}{L_0},\frac{2F_0}{\varepsilon}\Bigr\}$\newline
 \hspace*{1em}$\displaystyle{}+\log_+\frac{2F_0}{L_0R^2}\Bigr)$\\
Upper and lower bounds\textsuperscript{\dag\dag}\newline\citet{rudomanov2026everest} & sym.\ $(L_0,L_1)$ of $f$ & \yes & $f(\hat x)-f^\ast\le\varepsilon$
 & $\displaystyle\widetilde\Theta\Bigl(\sqrt{\frac{L_0R^2}{\varepsilon}}+L_1R\Bigr)$;\newline
 if $F_0\le4\varepsilon$: $\displaystyle\widetilde\Theta\Bigl(\sqrt{\frac{L_0R^2}{\varepsilon}}+(L_1R)^{2/3}\Bigr)$\\
\midrule
\multicolumn{5}{@{}l}{\textit{Finite sum, generalized smoothness}}\\
$(L_0,L_1)$-SGD, SGD-PS\textsuperscript{\ddag\ddag}\newline\citet{gorbunov2024methods} & sym.\ $(L_0,L_1)$ of each $f_i$ & \yes & see note~\textsuperscript{\ddag\ddag}
 & $\displaystyle O\Bigl(\frac{L_0R^2}{\varepsilon}\Bigr)$\newline
 if $\delta=0$ (all $f_i$ share a minimizer)\\
$(L_0,L_1)$-SPIDER\textsuperscript{\S\S}\newline\citet{reisizadeh2023variance} & avg.\ $(L_0,L_1)$ & \no & $\norm{\nabla f(\hat x)}\le\varepsilon$ with probability ${}\ge1/2$
 & $\displaystyle O\Bigl(n+\frac{\sqrt n\,L_0F_0}{\varepsilon^2}\Bigr)$\\
\midrule
\textbf{YETI} and lower bound\textsuperscript{\P\P}\newline\textbf{This paper}, Theorem~\ref{thm:main} & sym.\ $(L_{0,i},L_1)$ of each $f_i$ & \yes & $\E f(\hat x)-f^\ast\le\varepsilon$
 & $\displaystyle\widetilde\Theta\Bigl(n\bigl(1+(L_1R)^{2/3}\bigr)+\sqrt{\frac{n(L_0+L_1^2\delta)R^2}{\varepsilon}+n(L_1R)^2}\Bigr)$\\
\bottomrule
\end{tabularx}
\par\smallskip\raggedright\scriptsize
\emph{Notation.} $R=\norm{x_0-x^\ast}$, $F_0=f(x_0)-f^\ast$, $\log_+=\max\{\log,0\}$;
$O(\cdot)$ and $\Omega(\cdot)$ hide absolute constants, and $\widetilde\Theta(\cdot)$ also hides
logarithmic factors. Column ``Convex'': \yes{} --- the result is obtained
for convex $f_i$, \no{} --- convexity is not assumed.\par
\emph{Smoothness.} $L$: the gradient of each $f_i$ is Lipschitz with constant~$L$.
sym.: symmetric $(L_0,L_1)$-smoothness \citep[Assumption~1.3]{gorbunov2024methods},
$\norm{\nabla g(x)-\nabla g(y)}\le\bigl(L_0+L_1\sup_{u\in[x,y]}\norm{\nabla g(u)}\bigr)\norm{x-y}$
for all $x,y$; for $C^2$ functions this is equivalent to $\norm{\nabla^2g(x)}\le L_0+L_1\norm{\nabla g(x)}$.
Avg.: averaged $(L_0,L_1)$-smoothness \citep[Assumption~1(ii)]{reisizadeh2023variance},
$\bigl(\frac1n\sum_{i=1}^n\norm{\nabla f_i(x)-\nabla f_i(y)}^2\bigr)^{1/2}\le\bigl(L_0+L_1\norm{\nabla f(x)}\bigr)\norm{x-y}$
for $\norm{x-y}\le1/L_1$: the components are bounded only in the mean-square sense,
and the right-hand side contains the gradient of the sum $f$ rather than the gradients of the components.\par\smallskip
\textsuperscript{\dag}~For $0<\varepsilon\le F_0/2$.\quad
\textsuperscript{\ddag}~$D_0=2F_0+\frac32LR^2$; the bound is stated for $D_0\ge\varepsilon$.\quad
\textsuperscript{\dag\dag}~For deterministic methods with $L_0R^2\ge\varepsilon$.\par
\textsuperscript{\S}~For randomized algorithms in sufficiently high dimension, on the ball $\norm x\le R$.\par
\textsuperscript{\P}~The result is obtained for $\ell$-smooth functions,
$\norm{\nabla^2f(x)}\le\ell(\norm{\nabla f(x)})$: for a general $\ell$ the bound is
$O\bigl(\sqrt{\ell(0)}\,R/\sqrt\varepsilon+T\bigr)$ with an implicitly defined $T$
\citep[Corollary~5.3]{tyurin2026accelerated}. Here it is reduced to $\ell(s)=L_0+L_1s$
\citep[Theorem~4.3]{tyurin2026accelerated} with $\bar R=R$, $\Gamma_0=2F_0/R^2$.\par
\textsuperscript{\ddag\ddag}~Number of stochastic gradients until
$\min_k\E\min\bigl\{\frac{\nu L_0}{4nL_1^2},f(x_k)-f^\ast\bigr\}\le\varepsilon$, where
$\nu=e^{-\nu}$ \citep[Theorems~7.2, 7.3]{gorbunov2024methods}; SGD-PS additionally
requires knowing all $\min f_i$.\par
\textsuperscript{\S\S}~For $\varepsilon<L_0/(20L_1)$; here $F_0=f(x_0)-\inf f$
\citep[Theorem~4]{reisizadeh2023variance}.\par
\textsuperscript{\P\P}~The $L_{0,i}$ may differ, $\frac1n\sum_iL_{0,i}\le L_0$.
Complexity in the sense of Definition~\ref{def:oracle}; the hidden factors are $\log^2N_0$ (lower bound)
and $\log^2N$ (upper bound). Assumptions~\ref{asm:components}--\ref{asm:accuracy}, including $n\ge2$,
$\Delta\ge2\varepsilon$, $L_1^2R^2\delta/\varepsilon\ge3$; the case $\delta=0$ is not covered.\par
\end{center}

\endgroup
\end{document}